\documentclass{amsart}
\usepackage{latexsym}
\usepackage{amsopn,amscd}
\usepackage{amsfonts}
\usepackage{amssymb}
\usepackage{amsthm}
\usepackage{amsmath}
\usepackage{amssymb}
\usepackage{mathrsfs}
\usepackage{graphicx}
\usepackage[all]{xy}

\usepackage{a4wide}
\usepackage{framed}
\usepackage{stmaryrd}
\usepackage{enumerate,framed,color}
\usepackage{eucal}

\newtheorem{Theorem}{Theorem}[section]
\newtheorem{Proposition}[Theorem]{Proposition}
\newtheorem{Lemma}[Theorem]{Lemma}
\newtheorem{Corollary}[Theorem]{Corollary}
\newtheorem{Remark}[Theorem]{Remark}
\newtheorem{Example}[Theorem]{Example}

\numberwithin{equation}{section}

\def\momen{\delta}
\def\Mo{V}
\def\ML{M}
\def\cotan{\tau^*}
\def\BB{B}
\def\QQ{Q}
\def\empha{\em}
\def\emphas{}

\def\pemphas{}

\def\Sigm{\mathcal S}
\def\prin{\xi}
\numberwithin{equation}{section}

\begin{document}

\title[Extensions of Lie-Rinehart algebras and cotangent bundle reduction]
{Extensions of Lie-Rinehart algebras and cotangent bundle reduction}
\author{J. Huebschmann}
\address{(JH)USTL, UFR de Math\'ematiques, CNRS-UMR 8524, 59655 Villeneuve d'Ascq Cedex, France}
\email{Johannes.Huebschmann@math.univ-lille1.fr}
\author{M. Perlmutter}
\address{(MP)Departamento de Matematica, Universidade Federal de Minas Gerais, Brasil}
\email{matthew@mat.ufmg.br}
\author{T. S. Ratiu}
\address{(TSR)Section de Math\'ematiques and Bernoulli Center, \'{E}cole
Polytechnique F\'{e}d\'{e}rale de Lausanne, CH-1015 Lausanne, Switzerland}
\email{tudor.ratiu@epfl.ch}

\maketitle
\begin{abstract}
{Let $Q$ denote a smooth manifold acted upon smoothly by a Lie group
$G$. The $G$-action lifts to
an action on the total space $\mathrm{T}^*Q$ of the cotangent
bundle of $Q$ and hence on the standard symplectic Poisson algebra
of smooth functions on $\mathrm{T}^*Q$. The Poisson algebra of
$G$-invariant functions on $\mathrm{T}^*Q$ yields a Poisson
structure on the space $(\mathrm{T}^*Q)\big/G$ of $G$-orbits. We
relate this Poisson algebra with extensions of Lie-Rinehart algebras
and derive an explicit formula for this Poisson structure
in terms of differentials.
We then show, for the particular case where the $G$-action on $Q$ is
principal, how an explicit description of the Poisson algebra
derived in the literature by an ad hoc construction 
is essentially a special case of the formula
for the corresponding extension
of Lie-Rinehart algebras.  By means of various examples, we also 
show that this kind 
of description breaks down when 
the $G$-action does not define a principal bundle.}
\end{abstract}


\section{Introduction}
Let $Q$ denote a smooth manifold acted upon smoothly by a Lie group $G$. The $G$-action lifts
to a symplectic action on the total space $\mathrm{T}^*Q$ of the cotangent
bundle of $Q$, and hence $G $ acts on the ordinary symplectic Poisson
algebra of smooth functions on $\mathrm{T}^*Q$. Collective
symplectic reduction then leads to the Poisson algebra of
$G$-invariant functions on $\mathrm{T}^*Q$, and the question
arises of determining this Poisson algebra in terms of the space of $G$-orbits in $Q$
and the requisite additional data. 

The simplest non-trivial case one could perhaps think of
arises when we take  $Q$ to be the group $G$ itself, with $G$-action
via left translation. In this case, the Poisson algebra of
$G$-invariant functions on $\mathrm{T}^*Q$  simply comes down to
the familiar Lie-Poisson algebra of smooth functions on the dual
$\mathfrak{g}^*$ of the Lie algebra $\mathfrak{g}$ of $G$. A
non-degenerate $G$-invariant scalar product on $\mathfrak{g}$ then
defines a $G$-invariant Riemannian metric or, equivalently,
a kinetic energy on $Q=G$, and the associated equations of motion, that is, the equations defining the geodesics, are determined
by the Euler-Poincar\'e equations on $\mathfrak{g}$. For
$G=\mathrm{SO}(3,\mathbb R)$, viewed as the configuration space of
a rigid body moving freely about its center of mass, the inertia tensor of the body is precisely such a
quadratic form and the Euler-Poincar\'e equations are the
familiar Euler equations for the motion of a free rigid body
\cite{abramars}, \cite{arnobook}, \cite{marssche}.

More generally, when $Q \rightarrow B$ is a principal
$G$-bundle, the manifold $(\mathrm{T}^\ast Q)/G$ carries
a natural Poisson structure. Decomposing it in
terms of the canonical Poisson structures on the total space
 $\mathrm{T}^\ast (Q/G)$ 
of the cotangent bundle of the \textit{shape space} $Q/G$
and the Lie-Poisson structure on
the dual $\mathfrak{g}^\ast$ of the Lie algebra 
$\mathfrak{g}$ of $G$ requires as 
an additional ingredient 
a principal connection. Under that hypothesis, the
Poisson bracket has been determined in \cite{Mont},\cite{MMR1984}
 and used to study the motion of a particle
in a Yang-Mills field (the Wong equations) and various
continuum mechanical problems. A 
coordinate free proof of the resulting
formula 
was given in \cite{zaalaone}. A characterization of the symplectic
leaves of this Poisson structure was found 
in
\cite{MP}, based on previous results in \cite{S77},
\cite{Weinstein1978}. 

For the case where the $G$-action is not free
(and in particular not principal), 
the problem has been 
studied, is the subject of ongoing research,
and partial results are known. 
The singular symplectic reduced space at zero, see \cite{SjLe1991},
has additional structure induced by the cotangent bundle
character of $\mathrm{T}^\ast Q$; this structure has been
completely described in \cite{PeROSD2007}. For the particular case
where $Q$ has
a single $G$-orbit type, 
see \cite{or}, the same program has been carried out both
in the symplectic and Poisson categories for
arbitrary momentum values in \cite{hochrain}, \cite{PeRO2009}. 
In the completely general case of
a proper non-free action of $G$ on $Q$, the description of the Poisson
geometry of the stratified space $(\mathrm{T}^\ast Q)/G$
is built on the previously quoted results and elaborated on in
work in progress.

In the present paper we will show that the
constructions in the quoted references admit a simple explanation in
terms of what we refer to as the {\em tautological Poisson structure
of a Lie-Rinehart algebra\/} and, furthermore, in terms of
extensions of Lie-Rinehart algebras. One of us has already explored this
tautological Poisson structure in \cite{poiscoho}
(without having introduced that terminology). At the present stage,
suffice it to say the following:
Take $(A,L)$ to be the
Lie-Rinehart algebra 
$(C^{\infty}(M),\mathfrak{X}(M))$ which consists of the smooth
functions $C^{\infty}(M)$ and the smooth vector fields $
\mathfrak{X}(M)$ on a  smooth manifold $M$; now the tautological
Poisson algebra is the algebra of smooth functions on the
total space $\mathrm{T}^*M$ of the cotangent bundle of $M$ that are
polynomial 
on the fibers, endowed with the 
(negative of the)
ordinary cotangent bundle
Poisson structure.
For clarity we note that 
the tautological Poisson structure associated to a general Lie-Rinehart 
algebra is a purely algebraic construction; the algebraicity of the
construction is here reflected in the fact
that in the particular case under discussion 
we take functions that are polynomial on 
the fibers. 
We hope 
that the terminology 
\lq tautological Poisson structure\rq\ 
does not 
conceal from the browsing reader that
this concept is a substantial generalization
of the Poisson structure on the dual of a Lie algebroid;
we are indebted to the referee for having observed this risk.

We mention in passing that Lie-Rinehart algebras yield
substantial new insight elsewhere,
see \cite{poiscoho}--\cite{quasi} and, in particular,
are a crucial tool in a systematic approach to singular K\"ahler reduction \cite{kaehler}--\cite{adjoint}.
In particular Lie-Rinehart algebras
provide the appropriate algebraic tool
for a satisfactory description of
the behavior of the theory at singularities. 
The corresponding singular K\"ahler 
quantization was developed in \cite{qr}
and a particular gauge model for quantum mechanics on a stratified
space was developed in \cite{hurusch}; an unforeseen physical phenomenon
(a tunneling) arose out of that model.

At various stages in the paper, 
to make the paper somewhat more easily accessible to a wider audience,
we make explicit the necessary transition
from the algebra to the geometry that underlies
our initial problem.
In this vein, in Section \ref{derifo},
we reproduce the relationship between formal differentials,
derivations, and vector fields.
A good understanding of this relationship
is indispensable for the rest of the paper.
In Section \ref{sec_2}
we explain the tautological Poisson algebra of a Lie-Rinehart
algebra. In Section \ref{sec_extensions} 
we develop the behavior of the 
tautological Poisson algebra relative to an extension of Lie-Rinehart 
algebras. Here the general algebraic result is Theorem 
\ref{general}; it spells out the corresponding Poisson algebra
in terms of data that characterize the underlying extension
of Lie-Rinehart algebras. In Section \ref{desc} 
we then proceed to describe the
Poisson structure 
arising from an extension of Lie-Rinehart algebras
in terms of the formal differentials resulting from the kernel,
of those resulting from the quotient  and, furthermore,
in terms of the requisite additional data;
see Theorem \ref{thm1} for details for the case where
the quotient is a general Lie-Rinehart algebra
and 
 Theorem \ref{thm2} for the particular case where
the quotient is the Lie-Rinehart algebra of derivations of a commutative
algebra $A$.
These are structure results in the 
purely algebraic
theory of Lie-Rinehart algebras
that are interesting in their own right.

As for our initial problem, that of describing
the Poisson structure of $G$-invariant functions on $\mathrm{T}^*Q$
in terms of the space of $G$-orbits in $Q$ and the requisite additional 
data, when the action of $G$ on $Q$ is principal, the Poisson algebra of
$G$-invariant functions on $\mathrm{T}^*Q$ can be completely
characterized in terms of the tautological Poisson structure of a
suitable Lie-Rinehart algebra and  therefore in terms of extensions
of Lie-Rinehart algebras, cf. Proposition \ref{zaal}.
This observation recovers the familiar expressions for the gauged
Lie-Poisson bracket.
This kind
of description generalizes to Lie algebroids with constant rank
structure map, see Theorem \ref{alg}, in particular, to
transitive Lie algebroids, as explained in Corollary \ref{cor1}. 

However, in the situation 
when the $G$-action on $Q$ is not free, 
this kind of description and hence the characterizations
given in the quoted references can only partially recover the
Poisson algebra of $G$-invariant functions on $\mathrm T^*Q$. 
In Section \ref{nonreg} below, we
justify this observation with a number of examples. 
The description 
of the Poisson algebra
in terms of the tautological Poisson algebra 
associated to a Lie-Rinehart algebra
that we elaborate upon in the present paper
opens the door
to a systematic treatment of singular cotangent bundle reduction
by means of techniques from real algebraic geometry;
the examples in Section \ref{nonreg} only hint at a huge unexplored 
territory here, similar to that arising in singular K\"ahler reduction
and singular  K\"ahler quantization.

Our approach is essentially algebraic: rather than working with the Poisson 2-tensor we exploit the
corresponding Poisson 2-form, and the requisite algebra provides a road 
map through the jungle of covariant derivatives,
see, e.g., Theorems  \ref{thm1} and \ref{thm2}.
Our approach immediately translates to other categories, for example,
it applies to 
the orbit space of the
total space of the cotangent bundle on a non-singular algebraic variety 
$Q$,
endowed with the associated Poisson structure on the structure sheaf,
when an algebraic group acts principally on $Q$.
We spare the reader and ourselves these added troubles here.
We plan to explore, at another occasion, the significance of our
approach for the equations of motion.

Part of the present apparatus was first developed
for the smooth case, phrased in the language of Lie algebroids,
in papers of T. Courant, A. Weinstein, and others;
see e.~g. 
 \cite{CW}, \cite{courant0}, \cite{courant},
\cite{mackbtwo},
 and the references there.
\medskip

\noindent\textbf{Conventions.} 

We will write the de Rham operator on a smooth manifold 
 as $\mathbf{d}$. 
In the standard formalism, given a smooth symplectic 
manifold $(M, \omega)$,
the Hamiltonian vector field $X_h$ associated to the function $h \in C ^{\infty}(M)$
is defined by the identity
\[
\mathbf{i}_{X_h}\omega:=\omega(X_h,\,\cdot \,)
=  \mathbf{d}h;
\]
given  another smooth function $h$ on $M$, the associated Poisson
bracket of $f$ and $h$ is defined by 
\[
\{f,h\}:= \omega(X_f,X_h)= X_h(f) = -X_f(h).
\]
In particular, when $M$ is the total space $\mathrm T^* Q$ of the
cotangent bundle an a smooth manifold $Q$
(the configuration space of a mechanical 
system)  and $\vartheta$  the tautological 1-form on  $\mathrm T^* Q$,
the standard cotangent bundle symplectic form
is defined by $\omega:=-\mathbf d\vartheta$ 
(see, e.g., \cite{abramars, arnobook, marsden81, MMOPR, MP, MaRa1994, or}), 
so that for $h \in 
C^{\infty}(\mathrm{T}^\ast Q)$, the Hamiltonian
vector field $X_h \in \mathfrak{X}( \mathrm{T}^\ast Q)$ 
is given by the identity $ -\mathbf{i}_{X_h}\mathbf{d} \vartheta = \mathbf{d} h.
$
This convention yields the standard Poisson brackets
$\{q,p\}=1$ etc. and Hamilton's equation in the form
$\dot{f}=\{f,h\}$; in particular, with $h= {p^2 \over 2}$, we get the standard relationship
$
\dot q = \{q,h\} = p,
$
and the Poisson bracket of the two functions $f$ and $h$ on
$\mathrm T^* Q$ 
is given by $\{ f,h\}=- \mathbf{d}\vartheta(X_f, X_h)$.

However the Poisson bracket $\{\, \cdot \, , \,\cdot \,\}$  on  $C^{\infty}(\mathrm T^*Q)$
we shall use in Section \ref{sec_2} 
assigns to $f,h\in C ^{\infty}(\mathrm{T}^* Q)$ the function
$\{ f,h\}= \mathbf{d}\vartheta(X_f, X_h)$.
In standard cotangent bundle coordinates $(q^i, p _i)$ on $\mathrm{T}^\ast Q$ 
the value $\{f, h\}$ 
of the Poisson bracket of $f$ and $h$
is then given by the expression
\begin{equation}
\{f, h\}= \sum_{i=1}^{\dim M}\left(
\frac{\partial f}{\partial p_i}
\frac{\partial h}{\partial q^i} - 
\frac{\partial h}{\partial p_i}
\frac{\partial f}{\partial q^i}
\right)
\label{poisexp}
\end{equation}
and Hamilton's equations defined by $h \in 
C^{\infty}(\mathrm{T}^\ast Q)$ then take the form
$\dot{f}=\{h,f\}$. 

Throughout, we use the symbol
$\blacklozenge$ to indicate the end of a remark.

\section{Derivations and formal differentials}
\label{derifo}
Let $R$ be a commutative ring with identity,
$A$ an  $R$-algebra with identity,
at first not necessarily commutative,
and $M$ an $A$-bimodule.
We recall that an $M$-valued $R$-{\em derivation\/} 
or $M$-valued {\em derivation over\/} $R$
(relative Hochschild 1-cocycle) 
$d\colon A \to M$ on $A$
is a morphism of $R$-modules satisfying the additional requirement
\begin{equation}
d(ab)=(da)b+a(db), \ a,b \in A.
\label{der1}
\end{equation}
When $A$ is commutative and $M$ an ordinary $A$-module,
the requisite left and right $A$-module structures on $M$
being identified, 
the identity \eqref{der1} simplifies to
\begin{equation}
d(ab)=b(da)+a(db), \ a,b \in A.
\label{der2}
\end{equation}
The alerted reader will notice that this is a purely algebraic
identity and not just
beginning calculus.

Let now $A$ be a commutative
$R$-algebra with unit and $\mathcal A$ a commutative $A$-algebra
with unit.
We denote the 
$R$-module
of derivations of $A$ over $R$ by
$\mathrm{Der}(A)$, that of
derivations of $\mathcal A$ over $R$ by
$\mathrm{Der}(\mathcal A)$, and the $A$-module 
of
derivations of $\mathcal A$ over $A$ by  
$\mathrm{Der}(\mathcal A\big|A)$;
actually 
$\mathrm{Der}(A)$ 
acquires an $A$-module 
and an $R$-Lie algebra structure
and $\mathrm{Der}(\mathcal A)$ 
an $\mathcal A$-module 
and an $A$-Lie algebra structure.
More generally given an $A$-module $M$, we denote the
$A$-module of $M$-valued $R$-derivations on $A$
by $\mathrm{Der}(A,M)$ and, likewise,
given an $\mathcal A$-module $\mathcal M$, we denote the
$\mathcal A$-module of $\mathcal M$-valued $R$-derivations on $\mathcal A$
by $\mathrm{Der}(\mathcal A,\mathcal M)$
and the
$\mathcal A$-module of $\mathcal M$-valued $A$-derivations on $\mathcal A$
by $\mathrm{Der}(\mathcal A\big|A,\mathcal M)$.

Let 
$\mathrm D_A$ denote the  $A$-module of
formal differentials of $A$ over $R$,
likewise $\mathrm D_{\mathcal A}$  the $\mathcal A$-module of
formal differentials of $\mathcal A$ over $R$,
and $\mathrm D_{\mathcal A\big|A}$ the
$\mathcal A$-module of
formal differentials of $\mathcal A$ over $A$.
We will write the elements of $\mathrm D_A$ as $da$ ($a\in A$) etc.
By definition,
the $A$-module $\mathrm D_A$ represents the derivation functor
$\operatorname{Der}(A,\,\cdot\,)$
on the category of $A$-modules:
The universal derivation $d\colon A \to \mathrm D_A$ 
that assigns $da\in \mathrm D_A$ to $a \in A$
has the property
that,
given an $A$-module $M$, the $A$-module morphism
\begin{equation}
\mathrm{Hom}_A(\mathrm D_A, M)
\longrightarrow \mathrm{Der}(A,M),
\ 
\varphi \mapsto \varphi \circ d,\ \varphi \in \mathrm{Hom}_A(\mathrm D_A, M),
\label{univ}
\end{equation}
is an isomorphism of $A$-modules, plainly natural in $M$.

We will now recall a number of facts that we will subsequently use.
To this end,
let $\Mo$ be a finitely generated 
projective
$A$-module and let $\Sigm=\Sigm_A[\Mo]$, the symmetric
$A$-algebra on $\Mo$.
Without a hypothesis of the kind that $\Mo$ be finitely generated
and projective, a number of sequences that we exploit below
need no longer be exact; 
moreover the theory to be developed below involves
the $A$-dual of $\Mo$ and, without that hypothesis
it would be hard to handle that $A$-dual effectively.
In fact, without a hypothesis of that kind, many of the basic constructions
in the present paper break down; 
this is, perhaps, the reason
why the presence of singularities spoils the description of the
Poisson algebra we are looking for. Indeed, in geometry,
the projectivity hypothesis (of the space of
derivations of the algebra of functions, viewed as a module over
that algebra of functions) corresponds to the smoothness of the
underlying space.

\subsection{}
The sequence of inclusions 
$R \subseteq A \subseteq \mathcal A$ determines a canonical exact 
sequence
\begin{equation}
\mathcal A \otimes _A \mathrm D_A \longrightarrow \mathrm
D_{\mathcal A} \longrightarrow \mathrm D_{\mathcal A|A}
\longrightarrow 0 \label{dif}
\end{equation}
of $\mathcal A$-modules.

\subsection{} The exact sequence \eqref{dif} dualizes to a canonical 
exact sequence
\begin{equation}
0
\longrightarrow
\mathrm {Der}{(\mathcal A|A)}
\longrightarrow
\mathrm {Der}{(\mathcal A)}
\longrightarrow
\mathrm {Der}{(A,\mathcal A)}
\label{ddif}
\end{equation}
of $\mathcal A$-modules.

\subsection{}
The canonical $\mathcal S$-module morphism
\begin{equation}
\mathrm{Hom}_A(\Mo,\Sigm)
\longrightarrow
\mathrm {Der}({\mathcal S|A}),
\quad
 \phi \mapsto d_{\phi}, \quad 
d_{\phi}(\alpha) :=\phi(\alpha),\ \alpha \in \Mo,
\label{dfordif}
\end{equation}
is an isomorphism. 
Further, the
canonical morphism 
$\Sigm \otimes_A \mathrm{Hom}_A(\Mo,A) \to \mathrm{Hom}_A(\Mo,\Sigm)$ 
is an isomorphism of $\Sigm$-modules, 
and therefore
the canonical $\mathcal S$-module morphism
\begin{equation}
\Sigm \otimes_A \mathrm{Hom}_A(\Mo,A)
\longrightarrow
\mathrm {Der}({\mathcal S|A})
\label{ddfordif}
\end{equation}
is an isomorphism of $\Sigm$-modules as well. 
Consequently the assignment to $\alpha \in \Mo$,
viewed as a member of $\Sigm =\Sigm_A[\Mo]$, of the differential
$d\alpha \in \mathrm D_{\Sigm|A}$ induces an isomorphism
\begin{equation}
\Sigm \otimes_A \Mo
\longrightarrow 
\mathrm D_{\Sigm|A}
\label{fordif}
\end{equation}
of $\Sigm$-modules.
The isomorphism \eqref{dfordif} is inverse to the dual of \eqref{fordif}.

\subsection{}
With $\mathcal S$ substituted for $\mathcal A$, 
the exact sequence \eqref{dif}
of $\mathcal S$-modules  refines to
\begin{equation}
0
\longrightarrow
\Sigm\otimes_A \mathrm{D}_A
\longrightarrow
\mathrm \mathrm{D}_\Sigm
\longrightarrow
\Sigm\otimes_A \Mo
\longrightarrow
\label{dif2} 
0.
\end{equation}
Here the second unlabeled arrow sends
a member of $\Sigm\otimes_A \mathrm{D}_A$ of the kind  $s \otimes_A da$
($s \in \Sigm$ and $a \in A$) 
to the formal differential $sda$ of $\Sigm$ over $R$;
the injectivity of this arrow is, perhaps, not entirely obvious,
and we justify it in the appendix.
Further, given $\alpha \in \Mo$,
the third arrow sends the formal differential $d\alpha$
over $R$ to the differential of $\Sigm$ over $A$ induced by $\alpha$ 
relative to the $A$-module structure on $\Mo$, the $\Sigm$-module
$\mathrm D_{\Sigm|A}$ being identified with $\Sigm \otimes_A\Mo$ via
\eqref{fordif}. 

\subsection{} 
With $\Sigm$ substituted for $\mathcal A$, 
the exact sequence \eqref{ddif}
of $\Sigm$-modules  refines to
\begin{equation}
0
\longrightarrow
\mathrm {Hom}_A(\Mo,\Sigm)
\longrightarrow
\mathrm {Der}{(\Sigm)}
\longrightarrow
\Sigm\otimes_A\mathrm {Der}{(A)}
\longrightarrow
0;
\label{dddif}
\end{equation}
in terms of
the canonical isomorphisms
$\mathrm {Hom}_A(\Mo,\Sigm)\to \mathrm {Der}{(\Sigm\big|A)}$
and
$\Sigm\otimes_A\mathrm {Der}{(A)}\to \mathrm {Der}{(A,\Sigm )}$,
the two unlabeled arrows in  \eqref{dddif}
are given by the obvious associations:
$\mathrm {Der}{(\Sigm\big|A)}\longrightarrow
\mathrm {Der}{(\Sigm)}$
is the obvious inclusion and
$\mathrm {Der}{(\Sigm)}
\longrightarrow\mathrm {Der}{(A,\Sigm )}$
the obvious restriction.
We note that the sequence \eqref{dddif} arises from
\eqref{dif2} by dualization.
However, when we dualize the sequence \eqref{dddif}, 
we do not in general get back the sequence 
\eqref{dif2} unless the $A$-module $\mathrm D_A$ of formal differentials
of $A$ is reflexive.

\subsection{}\label{2.6} Let $B$ be a smooth manifold, and let $A=C^{\infty}(B)$.
A smooth version $\mathrm D^{\mathrm{smooth}}_A$ of 
the $A$-module
$\mathrm D_A$ of formal differentials 
is the $A$-module $\mathrm D^{\mathrm{smooth}}_A: 
= \Gamma(\cotan_B)$ 
of differential forms on $B$ or, equivalently,
smooth sections of the cotangent 
bundle $\cotan_B\colon 
\mathrm T ^\ast B \rightarrow B$ of $B$.
However the canonical morphism
$\mathrm D_A\to \mathrm D^{\mathrm{smooth}}_A$
of $A$-modules is surjective but not an isomorphism.
For example, over the real line, the formal differential
$d\sin (x) -\cos(x)dx$ is non-zero but is zero as a differential form.

\subsection{} 
To spell out a geometric interpretation of the exact sequences \eqref{dif2}
and \eqref{dddif},
let $B$ be a smooth manifold and 
$\lambda \colon  \mathcal L \to B$ a smooth vector bundle  
on $B$.  Let $\rho\colon  \mathcal L^* \to B$ 
denote the vector bundle on $B$
that 
is dual to $\lambda$ and, for $b \in B$,  
let $\mathcal L^*_b: =\rho^{-1}(b)$
be the fiber over $b\in B$.
When we view
 the total space $\mathcal L^*$ of the 
vector bundle $\rho$ merely as a smooth manifold
we denote that total space by $N$. 

The injective smooth vector bundle
map $\operatorname{Ver}\colon N \times _B \mathcal{L}^\ast 
\rightarrow \mathrm TN$ over $N$ given by
\[
\operatorname{Ver}(\alpha_b, \beta_b)\colon  = 
\left.\frac{d}{dt}\right|_{t=0}(\alpha_b+ t \beta_b) =\colon 
\operatorname{Ver}_{ \alpha_b}( \beta_b)
\]
has range the total space $V$ of the vertical bundle $V \rightarrow N$ and hence
induces an isomorphism of
vector bundles over $N$
from
$N \times_B  \mathcal L^* \rightarrow N$ to $V\rightarrow N$. 
Notice that, for every $b \in B$,
the map $\operatorname{Ver}_{\alpha_b}$ has the form
 $\operatorname{Ver}_{\alpha_b}\colon  \mathcal{L}_b^\ast\rightarrow 
T_{ \alpha_b}N $. 
Thus we obtain
the following exact sequence of vector bundles with base $N$,
spelled out here for the total spaces:
\begin{equation}
\begin{array}{lllllllll}
0&
\longrightarrow&
N \times _B\mathcal L^*
&\longrightarrow&
\mathrm TN
&\longrightarrow&
N \times _B\mathrm T B
&\longrightarrow&
0\\
&&(\alpha_b, \beta_b)&\longmapsto&
\operatorname{Ver}_{\alpha_b} (\beta_b) &
& \\
&&&&v_n &\longmapsto& (n, \operatorname{T}_n \rho(v_n))&
\end{array}, 
\label{1.3}
\end{equation}
where $n \in N$ and $\rho(n)=b$. 
This is the geometric analog of the sequence \eqref{dddif}.

Dualizing the sequence \eqref{1.3} we obtain an exact sequence
of vector bundles over $N$;  written out in terms of the total spaces,
this sequence has the form 
\begin{equation}
\begin{array}{lllllllll}
0&
\longleftarrow&
N \times _B \mathcal{L}
& \longleftarrow&
\mathrm{T}^\ast N
& \longleftarrow&
N \times _B\operatorname{T} ^\ast B
& \longleftarrow&
0\\
&&(n,\operatorname{Ver}_n ^\ast \Gamma_n)&\longmapsfrom&
\Gamma_n &
& \\
&&&&\operatorname{T}_n^\ast \rho( \alpha_b) &\longmapsfrom& (n, \alpha_b)&
\end{array}
\label{1.13}
\end{equation}
where $n \in N$ and $\rho(n) = b$.
This is the geometric analog of the sequence \eqref{dif2}.

Indeed, over the reals $\mathbb R$ as ground ring,
let $A=C^{\infty}(B)$ and take 
the spaces of sections in the smooth category; 
relative to the resulting sequence
$\mathbb R \subseteq A \subseteq C^{\infty}(\mathcal L^*)$
of real algebras, 
\eqref{1.13} thus yields
precisely an extension of the kind
\eqref{dif2} 
and \eqref{1.3}
one of the kind
 \eqref{dddif} 
except that the algebra 
$C^{\infty}(\mathcal L^*)$
of smooth functions on $\mathcal L^*$ has now been substituted for $\Sigm$.
We note that the two sequences
\eqref{1.3} 
and \eqref{1.13}
are dual to each other.

\subsection{}
Abstracting from the previous example, 
as before,
let $A$ be a
unital commutative $R$-algebra and suppose that, as an $A$-module,
$\operatorname{Der}(A)$ is finitely generated and projective.
Over a field as ground ring, when $A$ is an affine algebra,
this property can be taken as that defining
the notion of regularity of a commutative algebra.
In the present general case,
the
double $A$-dual 
$\mathrm D_A^{**}$ of $\mathrm D_A$, that is, the
 $A$-dual $\operatorname{Der}(A)^*=\mathrm{Hom}_A
(\operatorname{Der}(A),A)$ 
of $\operatorname{Der}(A)$,
is a finitely generated projective $A$-module 
as well. This module represents the 
derivation functor
$\operatorname{Der}(A,\,\cdot\,)$
on the category of finitely generated projective $A$-modules.
In particular, under heading \eqref{2.6}, the $A$-module
$\mathrm D_A^{\mathrm{smooth}}$ is precisely 
the double $A$-dual 
$\mathrm D_A^{**}$ of $\mathrm D_A$.

\section{The tautological Poisson algebra of a Lie-Rinehart
algebra}
\label{sec_2}

\subsection {Motivation}
To motivate the abstract notion of tautological Poisson algebra
associated to a Lie-Rinehart algebra, we will first explain a special case.

Let $Q$ denote a smooth finite dimensional manifold and 
$\tau_Q\colon  \mathrm{T}Q \rightarrow Q$  and 
$\cotan_Q\colon  \mathrm{T}^*Q \to Q$ its tangent and cotangent bundles,
respectively. The 1-form $\vartheta\in \Omega^1( \mathrm{T}^\ast Q)$ defined 
as the composite
\[
\begin{array}{lll}
\mathrm T \mathrm T^* Q &\longrightarrow
\mathrm T Q\times_Q \mathrm T^*Q
&\longrightarrow \mathbb{R} \\
V_{\alpha_q}&\longmapsto
\left(\mathrm T_{\alpha_q} \cotan_Q \left(V_{\alpha_q}\right), 
\alpha_q \right) & \longmapsto 
\left\langle \alpha_q, \mathrm T_{\alpha_q} \cotan_Q \left(V_{\alpha_q}
\right)\right\rangle,\ q \in Q,\ \alpha_q\in  \mathrm T_q^*Q,
\end{array}
\]
is the {\em tautological\/} 1-form on $\mathrm T^*Q$. 
In standard cotangent bundle coordinates $(q^j,p_j)$, 
the 1-form $\vartheta$ takes  the form
$p_j\mathbf d q^j$.
The standard symplectic form on $\mathrm{T}^* Q$ is 
the 2-form $- \mathbf{d}\vartheta $.
For formal reasons, 
we will instead 
work with the symplectic form 
$\mathbf{d}\vartheta$ on $\mathrm{T}^* Q$, and we will denote the associated
Poisson bracket by  $\{\,\cdot\, ,\,\cdot\,\}$.
This is the Poisson bracket which, in Darboux coordinates,
is given by \eqref{poisexp}.

It is well known that, relative to the obvious filtration of the
algebra of globally defined differential
operators on $Q$, the associated graded algebra is canonically
isomorphic to the 
algebra of smooth functions on $\mathrm{T}^* Q$ that 
are polynomial on the fibers.
In terms of the notation $(A,L)=(C^{\infty}(Q),\mathfrak X(Q))$,
this graded algebra is the symmetric $A$-algebra $\Sigm_A[L]$ on $L$.
It is a standard fact that the commutator of differential operators
induces a Poisson bracket
on $\Sigm_A[L]$, in fact,
the restriction of the Poisson bracket $\{\,\cdot\, ,\,\cdot\,\}$ on
 $C^{\infty}(\mathrm{T}^* Q)$ to the subalgebra $\Sigm_A[L]$.

We will now recall some of the details:
The algebra  $C^{\infty}(\mathrm{T}^*Q)$ being endowed
with the Poisson bracket  $\{\,\cdot\, ,\,\cdot\,\}$,
the vector space $\mathrm{Lin}(\mathrm{T}^\ast Q)$ of 
smooth functions on $\mathrm{T}^\ast Q$ that are linear on the fibers 
of the cotangent bundle $\cotan_Q\colon  \mathrm{T}^*Q \to Q$ 
of $Q$ is a Lie (\textit{not} Poisson) 
subalgebra of  $C^{\infty}(\mathrm{T}^*Q)$.
Each smooth vector field $X$ on 
$Q$ defines a smooth \textit{momentum function} 
$\momen_X\colon  \mathrm{T}^* Q \ni \alpha_q  \mapsto 
\left\langle \alpha_q, X(q) \right\rangle\in \mathbb{R}$
($q \in Q$) and 
the assignment to the vector field $X$ on $Q$ of the function
$\momen_X$  is a Lie algebra
homomorphism
$\momen\colon  \mathfrak{X}(Q) \rightarrow  \mathrm{Lin}(\mathrm{T}^\ast Q)$.
Indeed, write the action of the Lie algebra $\mathfrak{X}(Q)$ 
of vector fields on $Q$ on the algebra $C^{\infty}(Q)$ of smooth functions
on $Q$ as
\[
\mathfrak{X}(Q) \times C^{\infty}(Q) 
\longrightarrow  C^{\infty}(Q),
\quad  (X, f)\longmapsto X(f) ;
\]
now,
for any $f, h \in C ^\infty(Q)$ and $X,Y \in \mathfrak{X}(Q)$, the 
Poisson bracket 
on  $\Sigm_A[L]$ is determined by the
following relations: 
\begin{equation}
\label{fundamental_relations}
\{f \circ \cotan_Q, h \circ \cotan_Q\} = 0, \quad
\{\momen_X, f \circ \cotan_Q\} = X(f), \quad
\{\momen_X, \momen_Y\} = \momen_{[X,Y]}.
\end{equation}
We will refer to these relations later when we study the
general situation. 

Via the ring operations in 
$C^{\infty}(\mathrm{T}^* Q)$ and the action of
$\mathfrak{X}(Q)$ on $C^{\infty}(Q)$,
the functions on
$\mathrm{T}^* Q$ of the kind
 $f\circ \cotan_Q$,
where $f$ ranges over functions in $C ^{\infty}(Q)$,
together with the functions  on
$\mathrm{T}^* Q$ of the kind $\momen_X$, 
where $X$ ranges over
smooth vector fields on $Q$,
completely determine the 
Poisson bracket  $\{\,\cdot\, ,\,\cdot\,\}$ on 
$C^{\infty}(\mathrm{T}^* Q)$. 

\begin{Proposition}
The Poisson bracket on $C^{\infty}(\mathrm{T}^* Q)$
determined by \eqref{fundamental_relations}
is the symplectic Poisson bracket associated to
the 1-form $\mathbf{d}\vartheta$.
\end{Proposition}

\begin{proof} This is immediate:
By construction, in Darboux coordinates
$q_1,\ldots,q_m,p_1,\ldots,p_m$, the Poisson bracket
 \eqref{fundamental_relations} is given by
$
\{p_j,q_k\}=  \delta_{j,k}$, $1 \leq j,k\leq m$.
\end{proof}

We note that here and henceforth we do not distinguish in notation between 
vector fields on $Q$ and derivations of  $C^{\infty}(Q)$.

\begin{Remark}\label{rem0}
{\rm For a general smooth manifold $Q$,
with our conventions,
given a smooth function $f$, the assignment to
$f$ of the smooth vector field
$\{f,\,\cdot\,\}$ on $\mathrm{T}^* Q$ 
is a Lie algebra map from
the Lie algebra underlying the Poisson algebra
to the Lie algebra of smooth vector fields 
whereas, in the standard formalism, the ordinary
Hamiltonian vector field map
is an anti Lie morphism;
at the technical level,
the present convention has the advantage that a structure
map to be introduced later 
(the morphism $\pi^{\sharp}$ in Proposition
{\rm\ref{prop2}} below)
is compatible with the Lie-Rinehart structures,
and this will greatly simplify the exposition.} \quad $\blacklozenge$
\end{Remark}

\begin{Remark}\label{rem00}
{\rm Let $G$ be a connected 
Lie group, and let $\mathfrak g$ denote its Lie algebra,
more precisely, the Lie algebra of left-invariant vector fields
on $G$.
These are the fundamental vector fields for the right translation action
of $G$ on itself, and left translation 
\[
G\times \mathfrak g \longrightarrow \mathrm T G
\]
trivializes the tangent bundle of $G$.
Let $X,Y \in \mathfrak g$;
the Poisson bracket $\{\momen_X,\momen_Y\}$ 
of the associated momentum functions on $\mathrm T^* G$
is given by
$
\{\momen_X,\momen_Y\} =\momen_{[X,Y]};
$
this is a mere tautology since
the Lie bracket $[X,Y]\in \mathfrak g$
is defined in terms of the commutator bracket of vector fields
on $G$.

Let $G$ act on $Q$ from the right. Then the fundamental vector field map
$\mathfrak g \to \mathfrak X(Q)$ is a morphism of Lie algebras,
and the composite
\begin{equation}
\delta\colon \mathfrak g \longrightarrow C^{\infty}(\mathrm{T}^* Q)
\label{comp0}
\end{equation}
of the fundamental vector field map with the momentum function $\momen$
yields the {\em comomentum\/}
associated to the corresponding momentum mapping 
$\mu \colon \mathrm{T}^* Q\to \mathfrak g^*$.
The $G$-action being on the right of $Q$,
the $G$-equivariance of the momentum mapping is equivalent
to the momentum mapping being a Poisson map 
relative to the  Poisson bracket $\{\,\cdot\, , \,\cdot\,\}$
on $\mathrm{T}^* Q$ 
chosen above
(associated to $\mathbf d \vartheta$, that is, the negative of the standard bracket) 
and the ordinary Lie-Poisson bracket on $\mathfrak g^*$.
Hence the comomentum \eqref{comp0} is a morphism of Lie algebras.
We will write $\mathfrak g^*_+$ when we take $\mathfrak g^*$ to be 
endowed with  the standard
Lie-Poisson structure and
$\mathfrak g^*_-$
when we think of $\mathfrak g^*$ as being 
endowed with  the negative of the standard
Lie-Poisson structure.

We will now take $Q$ to be the group $G$ itself.
Let $\{\,\cdot\, , \,\cdot\,\}_{+}$
denote  the standard cotangent bundle Poisson bracket
on $\mathrm T^*G$ so that $\{\,\cdot\, , \,\cdot\,\}$
is the negative of the standard cotangent bundle Poisson bracket.
It is well known that
the composite
\begin{equation}
\mathrm T^* G
\longrightarrow \mathfrak g^* \times G
\longrightarrow
 \mathfrak g^*
\end{equation}
of right translation 
\begin{align*}
\mathrm T^* G
&\longrightarrow \mathfrak g^* \times G,\ \alpha_x \longmapsto (\alpha_x \circ \rho_x,x),\  x\in G,
\end{align*}
with the projection to the first component
yields
the momentum mapping
\begin{equation}
\mu_{\ell}\colon \mathrm T^*G \longrightarrow \mathfrak g^*,
\ 
 \alpha_x \longmapsto \alpha_x \circ \rho_x
\end{equation}
for the $G$-action
on $\mathrm T^*G$ by left translation.
The $G$-equivariance of the momentum mapping $\mu_{\ell}$ is 
well known to be equivalent to
$\mu_{\ell}$ being a Poisson map 
relative to  $\{\,\cdot\, , \,\cdot\,\}_{+}$; hence
$\mu_{\ell}$ induces isomorphisms
\begin{equation}
((\mathrm T^*G)\big/G,\{\,\cdot\, , \,\cdot\,\}_{+}) 
\longrightarrow \mathfrak g_+^*,
\quad
((\mathrm T^*G)\big/G,\{\,\cdot\, , \,\cdot\,\}) 
\longrightarrow \mathfrak g_-^*
\end{equation}
of Poisson manifolds.} \quad $\blacklozenge$
\end{Remark}

\subsection{The formal notion}
\label{subsec_formal}
The pair $(A,L):= (C^{\infty}(Q),\mathfrak{X}(Q))$ 
is among 
the standard examples of a {\em Lie-Rinehart\/} 
algebra 
over the real numbers as ground ring
\cite{poiscoho}, \cite{souriau}, \cite{lradq}, 
\cite{rinehart}
and, abstracting from the cotangent bundle situation, we 
are led to the {\em tautological Poisson algebra\/} 
associated to a Lie-Rinehart algebra that we now describe.
\medskip

Fix the commutative ring with identity $R$; the unadorned
tensor product symbol $\otimes$ will always refer to the tensor
product over $R$. Further, denote by $A$ a 
unital commutative $R$-algebra.
An $(R,A)$-{\em Lie algebra\/} \cite{rinehart} is a Lie algebra
$L$ over $R$ which acts on $A$ (from the left) by derivations (the
action being  written as ${(\alpha \otimes a) \mapsto 
\alpha(a)}$)
and is also an $A$-module (the structure map being written as
${(a\otimes \alpha) \mapsto a\alpha}$),
 in such a way that
suitable compatibility conditions are satisfied which generalize
standard properties of the Lie algebra of vector fields on a
smooth manifold viewed as a module over the ring of functions;
these conditions read
\begin{align} 
(a\,\alpha)(b) &= a\,(\alpha(b)), \quad \alpha \in L,
\,a,b \in A, \label{1.1.a}
\\
\lbrack \alpha, a\,\beta \rbrack &= a\, \lbrack \alpha, \beta
\rbrack + \alpha(a)\,\beta,\quad \alpha,\beta \in L, \,a \in A.
\label{1.1.b}
\end{align}
The pair $(A,L)$ is then referred to a {\em Lie-Rinehart
algebra\/}. Given two Lie-Rinehart algebras $(A,L)$ and $(A',L')$,
a {\em  morphism \/} ${ (\phi,\psi) \colon  (A,L) \longrightarrow
(A',L') }$ of {\em  Lie-Rinehart algebras\/} consists of morphisms $\phi\colon A \rightarrow A'$ and $\psi\colon L \rightarrow L'$ in the appropriate
categories that are compatible with the additional structure. With
this notion of morphism, Lie-Rinehart algebras constitute a
category. Apart from the example of smooth functions and smooth
vector fields on a smooth manifold, a related (but more general)
example is the pair consisting of a commutative algebra $A$ and
the $R$-module $\mathrm {Der}(A)$ of derivations of $A$ with the
obvious $A$-module structure and the commutator bracket
on  $\mathrm {Der}(A)$; here the commutativity of $A$ is
crucial.

Take a Lie-Rinehart algebra $(A,L)$, that is, $A$ 
is a commutative algebra and $L$ an $(R,A)$-Lie algebra. 
Denote by $\Sigm_A [L]$ the symmetric $A$-algebra on $L$.
The $L$-action on $A$ and the $R$-Lie bracket operation on $L$ induce a Poisson bracket
\begin{equation}
\{\,\cdot\, ,\, \cdot\, \} \colon  \Sigm_A [L] \otimes \Sigm_A [L]
\longrightarrow \Sigm_A [L]\label{3.16.1}
\end{equation}
on $\Sigm_A [L]$
completely determined by the relations
\begin{align} 
\{\alpha,\beta\} &= [\alpha,\beta],\quad
\alpha,\beta \in L, \label{3.16.2.1}
\\
\{\alpha,a\} &= \alpha(a) \in A,\quad a \in A,\,\alpha \in L,
\label{3.16.2.2}
\\
\{u,vw\}&= \{u,v\}w + v\{u,w\},\quad u,v,w \in \Sigm_A[L].
\label{3.16.2.3}
\end{align}
In particular, taking $u=a, v=b \in A$ 
and $w=\alpha\in L$, $\alpha\neq 0$ in \eqref{3.16.2.3}, 
we obtain $\{a,b\} = 0$.
We refer to the Poisson algebra $(\Sigm_A [L],\{\,\cdot\, ,\,
\cdot\, \} )$ as the {\em tautological Poisson algebra associated
to the Lie-Rinehart algebra\/} $(A,L)$. 
We note that, for the standard example $(C ^\infty(Q), \mathfrak{X}(Q))$
of a Lie-Rinehart algebra,  the tautological Poisson bracket 
just defined 
is the same as the Poisson bracket given in
\eqref{fundamental_relations}. 

The commutative
$R$-algebra $A$ being fixed, the tautological Poisson algebra is
plainly functorial in $(R,A)$-Lie algebras: A morphism 
$\phi\colon  L_1 \to L_2$ of $(R,A)$-Lie algebras induces 
a morphism
\begin{equation}
(\Sigm_A [L_1],\{\,\cdot\, ,\, \cdot\, \} )\longrightarrow
(\Sigm_A [L_2],\{\,\cdot\, ,\, \cdot\, \} ) \label{mor}
\end{equation}
of Poisson algebras.
\medskip

We will now justify the terminology \lq\lq tautological Poisson
algebra\rq\rq. 
To this end, we recall briefly how, for an
arbitrary Poisson algebra, an appropriate Lie-Rinehart algebra
serves as a replacement for the tangent bundle of a smooth
symplectic manifold.
\medskip

Let 
$(\mathcal{P},\{\, \cdot \, ,\, \cdot \, \})$
be a Poisson algebra.
The association 
$(du,dv)\mapsto \pi (du,dv) := \{u,v\}$,
as $u$ and $v$ range over $\mathcal{P}$,
yields a $\mathcal{P}$-valued $\mathcal{P}$-bilinear skew-symmetric 2-form
\begin{equation}
\pi= \pi_{\{\, \cdot \, ,\, \cdot \, \}} \colon  
\mathrm{D}_{\mathcal{P}} \otimes \mathrm D_{\mathcal{P}} \longrightarrow
\mathcal{P}
 \label{3.5.1}
\end{equation}
on $\mathrm D_{\mathcal{P}}$; this 2-form is said to be 
the {\em Poisson\/} 2-{\em  form\/}
of $(\mathcal{P},\{\, \cdot \, ,\, \cdot \, \})$. The adjoint
\begin{equation}
\pi^{\sharp} \colon  \mathrm D_{\mathcal{P}} \longrightarrow \mathrm
{Der}(\mathcal{P}) = \mathrm {Hom}_{\mathcal{P}}(\mathrm D_{\mathcal{P}},\mathcal{P}), \qquad \pi^\sharp(du) (dv) = \{u,v\},
\end{equation} 
of $\pi$ is a morphism of $\mathcal{P}$-modules and the formula
\begin{equation}
\label{bracket_on_d_a}
[a du,b dv] := a\{u,b\} dv + b\{a,v\} du + ab d\{u,v\}, \quad 
a,b,u,v \in \mathcal{P},
\end{equation}
yields a Lie bracket $[\,\cdot \, , \, \cdot \,]$ on $\mathrm D_{\mathcal{P}}$,
viewed as an $R$-module (see \cite{poiscoho} for more 
details). 
For the record we recall the following result,
established in \cite[Theorem 3.8]{poiscoho}.

\begin{Proposition}\label{prop2} The $\mathcal{P}$-module structure on $\mathrm D_{\mathcal{P}}$,
the $R$-bilinear bracket  $[\,\cdot \, , \, \cdot \,]$ on $\mathrm D_{\mathcal{P}}$, 
and the morphism
 $\pi^{\sharp}\colon  \operatorname{D}_ \mathcal{P} \rightarrow \operatorname{Der}( \mathcal{P})$
of $\mathcal{P}$-modules endow the pair
 $(\mathcal{P},\mathrm D_{\mathcal{P}})$
with a Lie-Rinehart algebra structure in such a way that 
$\pi^{\sharp}\colon  \mathrm D_{\mathcal{P}} \rightarrow \operatorname{Der}(\mathcal{P})$ is a
morphism of $(R, \mathcal{P})$-Lie algebras.
\end{Proposition}

Given a $\mathcal{P}$-module $D$,
we will use the notation $\mathrm{Alt}_{\mathcal{P}}(D,\mathcal{P})$
for the graded $\mathcal{P}$-algebra of $\mathcal{P}$-valued $\mathcal{P}$-multilinear alternating
forms on $D$.
When $(\mathcal{P},D)$ is a Lie-Rinehart algebra (over the ground ring $R$),
the standard Cartan-Chevalley-Eilenberg operator $d$ 
(equivalently: the de Rham operator evaluated formally),
endows  $\mathrm{Alt}_{\mathcal{P}}(D,\mathcal{P})$ with a
differential graded $R$-algebra structure, cf. \cite{poiscoho}.
Denote by $\mathrm D_{\{\, \cdot \, ,\, \cdot \,
\}}$ the $(R,\mathcal{P})$-Lie algebra 
$(\mathrm{D}_{\mathcal{P}},[\cdot,\cdot],\pi^{\sharp})$. 
The 2-form $ \pi_{\{\, \cdot \,,\, \cdot \, \}}, $ which 
is defined for {\em every\/} Poisson
algebra, is a 2-cocycle in the Rinehart algebra 
$(\mathrm{Alt}_{\mathcal{P}}(\mathrm D_{\{\, \cdot \, ,\, \cdot \,\}},\mathcal{P}),d)$. 
As opposed to other definitions
that involve a defining chain complex
and only refer to a Poisson algebra of smooth functions
on a smooth manifold,
in 
\cite{poiscoho},
the {\em Poisson cohomology\/} 
$\mathrm H^*_{\mathrm {Poisson}}(\mathcal{P},\mathcal{P})$ 
of the Poisson algebra $(\mathcal{P},
\{\, \cdot \,,\, \cdot \, \})$
was defined  as the cohomology of the Lie-Rinehart 
algebra $\mathrm{D}_{\{\, \cdot \, ,\, \cdot \,\}}$, 
\textit{ i.e.},
\begin{equation}
\mathrm  H^*_{\mathrm {Poisson}}(\mathcal{P},\mathcal{P}) =\mathrm
H^*\left(\mathrm {Alt}_{\mathcal{P}}(\mathrm D_{\{\, \cdot \, ,\,
\cdot \, \}},\mathcal{P}),d\right).
\end{equation}

When $\mathcal{P}$ is the algebra of smooth functions on a smooth
symplectic manifold $M$, endowed with the symplectic Poisson
structure,
the inverse of the symplectic structure induces a
vector bundle isomorphism from $\cotan_M$ onto $\tau_M$ and hence
an isomorphism from the smooth version $\mathrm
D^{\mathrm{smooth}}_{\{\, \cdot \, ,\, \cdot \, \}}\cong
\Omega^1(M)$ (the $C^{\infty}$-module of ordinary 1-forms on $M$)
of $\mathrm D_{\{\, \cdot \, ,\, \cdot \, \}}$ onto the ordinary
$(\mathbb R,\mathcal{P})$-Lie algebra of smooth vector fields on
$M$ which identifies the 2-form $\pi$ on $\operatorname{D}_ \mathcal{P}$ with the symplectic
structure. Furthermore, under the canonical isomorphism between
$\mathrm {Alt}_{\mathcal{P}}(\mathrm D^{\mathrm{smooth}}_{\mathcal{P}},\mathcal{P})$ and the exterior $\mathcal{P}$-algebra
$\wedge_{\mathcal{P}}(\mathfrak{X}(M))$ on $\mathfrak{X}(M)$,
the 2-form $\pi$ corresponds to the familiar Poisson 2-tensor. 
In addition, when
$P$ is a Poisson manifold and $d$  the de Rham operator,
 the bracket \eqref{bracket_on_d_a}
coincides with the Lie bracket on 1-forms on $P$; see
\cite[Proposition 3.3.8]{abramars} for the symplectic case and 
\cite[Subsection 3.12]{poiscoho}
and
\cite[Theorem 4.1]{vaisman} for the general Poisson case.
We emphasize,
cf. \cite{poiscoho},
 that the abstract construction given above
defines the
2-form $\pi$ for an arbitrary Poisson algebra, 
whether or not it arises from a smooth symplectic manifold. 

We return to a general Lie-Rinehart algebra $(A,L)$.
For the algebra
$\Sigm := \Sigm_A [L]$, endowed with the Poisson structure
\eqref{3.16.1}, 
the Poisson 2-form $\pi \colon  \mathrm D_\Sigm 
\otimes _\Sigm \mathrm D_\Sigm \to \Sigm$ given by 
\eqref{3.5.1} is a Poisson coboundary, that is, 
$\{\, \cdot \, , \,\cdot \, \}$ admits a Poisson potential 
\cite{poiscoho}. Indeed,
as an algebra, $\Sigm$ is generated by the elements of 
$A$ and those of $L$. Hence, as an $\Sigm$-module, 
$\mathrm D_\Sigm$ is generated by the formal differentials 
$da,\, a \in A$, and $d\alpha,\,\alpha \in L$. Let 
$\mathrm D_{\mathcal S|A}$ denote the $\mathcal S$-module 
of formal differentials of $\mathcal S$
over $A$; as an $\mathcal S$-module, 
$\mathrm D_{\mathcal S|A}$ is an induced module of 
the form $\mathcal S \otimes_A L$
generated by the formal differentials $d\alpha$, where 
$\alpha$ ranges over $L$.

A straightforward
calculation (see \cite[p.~92]{poiscoho}) shows that the 1-form $\vartheta \colon 
\mathrm D_\Sigm \to \Sigm$ given through the projection to
$\mathrm D_{\mathcal S|A}$ by
\begin{equation}
\vartheta(d\alpha):= \alpha,\quad \alpha  \in L, 
\label{3.16.3}
\end{equation}
is a Poisson potential for $\{\, \cdot  \, , \, \cdot \, \}$, 
that is,
\begin{equation}
\pi =d\vartheta \in \mathrm { Alt}^2_\Sigm(\mathrm D_{\{\, \cdot
\, , \, \cdot \, \}},\Sigm). \label{3.16.4}
\end{equation}

For later reference, we note that the 2-form $\pi$ is given
by the identities
\begin{align}
\pi(d\alpha,d\beta)&=[\alpha,\beta],\quad \alpha,\beta \in L,
\label{dif3}
\\
\pi(d\alpha,db)&=\alpha (b),\quad \alpha \in L,\ b \in A.
\label{dif4}
\end{align}
This description of the 2-form $\pi$ is consistent in the following sense: Given
$a \in A $ and $\alpha \in L $, when $a \alpha$ is viewed as a member of $\mathcal{S}=\mathcal{S}_A[L]$, the differential $d(a\alpha) $ over the ground ring $R $ takes the form
\[
d(a \alpha) = (da)\alpha + a(d\alpha)
\]
and, given  $\beta \in L $, the identities
\begin{align*}
\pi\left(d(a \alpha), d\beta\right)&=[a \alpha, \beta]
= a[ \alpha, \beta] - \beta(a) \alpha\\
\pi\left((da) \alpha + a d \alpha, \beta\right) & = 
\alpha\pi\left(da,\beta\right)+ a[ \alpha, \beta] =
- \beta(a) \alpha+ a [ \alpha, \beta]
\end{align*}
hold. We emphasize that the 2-form $\pi$ is $\mathcal{S}_A[L]$-bilinear. 
In particular, the following observation is immediate; we
spell it out, for later reference.

\begin{Proposition} \label{prop3}
The Poisson structure \eqref{3.16.1} on $\Sigm_A[L]$
recovers, in fact, is equivalent to,
the Lie-Rinehart structure on $(A,L)$:
Given $\alpha,\beta \in L$, the bracket $[\alpha,\beta]$
is determined by the identity
{\rm \eqref{3.16.2.1}} and,
given
$a \in A$ and $\alpha \in L$, the value $\alpha(a)$
of the action of $\alpha$ on $a$
is determined by the identity
{\rm \eqref{3.16.2.2}}.
Thus the
Lie-Rinehart structure can be reconstructed from the Poisson
structure.
\end{Proposition}

\begin{Remark}\label{rem3} {\rm 
Suppose that $\mathcal{P}$ 
is the
algebra of
smooth functions on the total space $N=\mathrm T^*B$ of the
cotangent bundle on a smooth manifold $B$,
endowed with the standard symplectic Poisson structure. 
Then the smooth version of
the 1-form $\vartheta$ given in \eqref{3.16.3}  is identical to the 
 standard 1-form on $N=\mathrm T^*B$, and the resulting symplectic structure
is therefore occasionally referred to as the {\em tautological symplectic
structure\/} on $N=\mathrm T^*B$. Actually, in this case, with the
notation $(A,L)=(C^{\infty}(B),\mathfrak{X}(B))$, in the
appropriate Fr\'echet topology,  $\Sigm_A[L]$ is dense in the
algebra $\mathcal{P}=C^{\infty}(\mathrm T^*B)$ of smooth functions on
$N=\mathrm T^*B$:
Indeed, using a locally affine open cover of $B$, 
in the standard manner,
we endow 
$C^{\infty}(\mathrm T^*B)$
with the locally convex Fr\'echet algebra structure
which, on the total space of the cotangent bundle
of any of the affine open subsets of $B$ in the cover, 
is given by the standard construction
involving the corresponding families of submultiplicative semi-norms.
In the resulting topology,
 $\Sigm_A[L]$ is dense in $C^{\infty}(\mathrm T^*B)$.
Under the canonical isomorphism between
\[
\mathrm
{Alt}_{\mathcal{P}}(\mathrm D^{\mathrm{smooth}}_{\mathcal
{P}},\mathcal{P})= \mathrm {Alt}_{\mathcal{P}}(\Omega^1(N),\mathcal
P)
\]
and the exterior $\mathcal{P}$-algebra $\wedge_{\mathcal
A}(\mathfrak{X}(N))$ on $\mathfrak{X}(N)$ recalled above 
for a general smooth symplectic manifold written there as $M$, 
the formulas \eqref{dif3}
and \eqref{dif4} then yield the Poisson 2-tensor for the cotangent bundle Poisson structure on $\mathrm T^*B$. 
The Poisson structure
\eqref{3.16.1} on the algebra $\Sigm_A [L]$ relative to a general
Lie-Rinehart algebra $(A,L)$ is formally of the same kind.

For an ordinary real Lie algebra $\mathfrak g$, 
viewed as an $(\mathbb R,\mathbb R)$-Lie algebra with trivial
$\mathfrak g$-action on $\mathbb R$, the tautological Poisson structure
plainly comes down to the algebraic version of the ordinary
Lie-Poisson structure on the dual $\mathfrak g^*$ of $\mathfrak g$.

These observations are intended to justify
the terminology {\em tautological Poisson
structure\/}.}
\quad $\blacklozenge$
\end{Remark}

\section{Extensions of Lie-Rinehart algebras and tautological
Poisson algebra}
\label{sec_extensions}

The algebraic analog of an
\lq\lq Atiyah sequence\rq\rq\  
(associated to a principal bundle)
or of a \lq\lq transitive Lie algebroid\rq\rq\  
(see below for details) is an extension of Lie-Rinehart algebras \cite{extensta}. We note that transitive Lie algebroids which do {\em not\/} 
arise from a principal bundle abound; see, e.g., 
\cite{almolone}.
We recall at the present stage that
a Lie algebroid is said to be {\em transitive\/}
when
its anchor map 
is fiberwise
surjective. 

Let $L'$, $L$, $L''$ be $(R,A)$-Lie algebras. An 
{\em extension\/}
of $(R,A)$-Lie algebras is a short exact sequence
\begin{equation}
\begin{CD}
\mathbf{e} \colon  0 @>>> L' @>>> L @>p>> L''
@>>> 0 \end{CD}
\label{2.1}
\end{equation}
in the category of $(R,A)$-Lie algebras; notice, in particular, that
the Lie algebra $L'$ necessarily acts trivially on $A$. 

Recall that 
an  $\mathbf{e}$-\emph{connection\/} associated to the
extension \eqref{2.1} is an $A$-module section
$\omega \colon  L'' \to L$ of 
the projection $p\colon  L \to L''$. The 
$\mathbf{e}$-\emph{curvature\/} of
$\omega$ is the  $A$-bilinear
map $\Omega\colon  L'' \times L'' \rightarrow L'$
of $A$-modules defined by the identity
\begin{equation}
[\omega(\alpha''),\omega(\beta'')] = 
\omega\left([\alpha'',\beta'']_{L''}\right) +
\Omega(\alpha'',\beta''), \ \alpha'', \beta''  \in L'',
\label{2.3}
\end{equation}
where $[\alpha'',\beta'']_{L''}$ 
refers to the Lie bracket in $L''$. 

As discussed below, the notions of 
$\mathbf{e}$-connection and $\mathbf{e}$-curvature
generalize the standard concepts of  principal connection 
and principal curvature. 
For a general extension 
\eqref{2.1} of $(R,A)$-Lie
algebras, a choice of connection $\omega$ induces an $A$-module
decomposition
\[
L= L' \oplus \omega(L'') \cong L' \oplus L''
\]
and, in terms of this decomposition, the Lie bracket 
$[\, \cdot \,, \, \cdot \, ]$ on $L$ takes the form
\begin{equation}
[(\alpha',\alpha''),(\beta',\beta'')]
=\left([\alpha',\beta']_{L'}+[\alpha',\beta'']
+[\alpha'',\beta']+\Omega(\alpha'',\beta''), 
[\alpha'',\beta'']_{L''}\right);
\label{brack}
\end{equation}
here $[\alpha',\beta']_{L'}$ refers to the Lie bracket in $L'$.
When $L''$ is a projective $A$-module,  one can always
find an $\mathbf{e}$-connection.

Since $\omega$ is
a section of  $p$, 
viewed as a morphism of $A$-modules, it
follows from \eqref{2.3} that  the values 
of $\Omega$ lie
in $L'$. Thus $\Omega$ is an $L'$-valued alternating
$A$-bilinear 2-form on $L''$.
We note that, 
similarly as in standard Schreier extension theory for groups,
suitably interpreted, 
$\Omega$ is a non-abelian 
2-cocycle. 

\begin{Example} \label{2.222} {\rm 
Here we take  the ground ring $R$ to be the field 
$\mathbb{R}$ of real numbers. 
Let $B$ be
a smooth finite-dimensional manifold and denote by  
$A$ the algebra $C^{\infty}(B)$ of smooth functions 
on $B$. Let $G$ be a Lie 
group and $\prin \colon  Q \to B$ a right principal 
$G$-bundle.  The vertical subbundle $\tau_Q|V\colon  V \to Q$ 
of the tangent bundle $\tau_Q\colon \mathrm{T}Q \rightarrow Q$ 
is well known to be trivial (beware, not equivariantly
trivial), having as fiber the Lie algebra $\mathfrak{g}$ of $G$, that is, $V \cong Q \times \mathfrak{g}$. 
These data fit into the fundamental exact sequence
\begin{equation}
\label{first_sequence}
\begin{array}{lllllllll}
0&
\longrightarrow&
Q \times \mathfrak{g}
&\longrightarrow&
\qquad \mathrm{T}Q
&\longrightarrow&
Q \times _B \mathrm{T}B
&\longrightarrow&
0\\
&&(q, \lambda) &\longmapsto&
\lambda_Q(q);\; v_q& \longmapsto&(q, T_q \xi(v_q))& 
\end{array}
\end{equation}
where $\lambda_Q(q)\colon  = \left.\frac{d}{dt}\right|_{t=0}
q \cdot \exp(t \lambda) \in \mathrm{T}_qQ$ defines
the fundamental vector field $\lambda_Q \in \mathfrak{X}(Q)$. 
We note that the range of the second arrow is the
total space  $V \subset \mathrm{T}Q$ of 
the subbundle of the tangent bundle whose sections are
vertical 
vector fields.
All arrows are equivariant vector bundle maps relative
to the diagonal action
\[
(q, \lambda)\cdot x\colon  = (q\cdot x, 
\operatorname{Ad}_{x ^{-1}} x ),
\quad  q \in Q , \ \lambda \in \mathfrak{g},\ x \in G,
\]
on $Q \times \mathfrak{g}$ and the tangent lifted action
on $\mathrm{T}Q$.
Dividing out the $G$-actions, we obtain
an extension
\begin{equation}
0 \longrightarrow \mathrm{ad}(\prin) \longrightarrow 
\tau_\QQ/G
\longrightarrow \tau_\BB \longrightarrow 0 \label{2.2.1}
\end{equation}
of vector bundles over $\BB$, where $\tau_\BB$ is the 
tangent bundle of $\BB$, and 
\begin{equation}
\tau_Q\big/G\colon  (\mathrm TQ)\big/G \longrightarrow B
\label{atiy2}
\end{equation}
is a transitive Lie algebroid over $B$.
Here
the notation
$\mathrm{ad}(\prin)$
refers to the vector bundle
 $Q \times_G \mathfrak{g}\to B$ associated to the principal bundle $\xi$
by the adjoint representation of $G$ on its Lie algebra 
$\mathfrak{g}$. The sequence \eqref{2.2.1} was introduced by 
Atiyah \cite{atiyaone} and is now usually referred 
to as the {\em Atiyah sequence\/} of the principal bundle 
$\prin$. 

The spaces of smooth sections 
$\Gamma_B\left(\mathrm{ad}(\xi)\right)$ 
and $\Gamma_B(\tau_Q/G)$ over
$B$  inherit 
obvious Lie
algebra structures, in fact $(\mathbb R,A)$-Lie algebra structures,
and
\begin{equation}
0 \longrightarrow \Gamma_B\left(\mathrm{ad}(\xi)\right) \longrightarrow \Gamma_B(\tau_Q/G)
\longrightarrow \mathfrak{X}(B) \longrightarrow 0 
\label{2.2.2}
\end{equation}
is an extension of  $(\mathbb R,A)$-Lie algebras; here
$\mathfrak{X}(B)$ is the ordinary Lie algebra of vector
fields on $B$ and $\Gamma_B\left(\mathrm{ad}(\xi)\right)$ is in an 
obvious way the Lie
algebra of the {\em group\/} of {\em gauge transformations\/} of $\xi$. 

Pick a connection 
on the principal bundle 
$\xi\colon Q \rightarrow Q/G=B$.
The associated horizontal lift operator
$\operatorname{hor}\colon Q \times _B \mathrm{T}B \rightarrow 
\mathrm{T}Q$, thought of as a $G$-equivariant vector 
bundle morphism over $Q$, is a section of the sequence
\eqref{first_sequence}. Taking the quotient by $G$ we obtain a 
section $[\operatorname{hor}]\colon \mathrm{T}B \rightarrow 
(\mathrm{T}Q)/G$ of the Atiyah sequence \eqref{2.2.1}.
Conversely, any section $\sigma\colon \mathrm{T}B \rightarrow (\mathrm{T}Q)/G$ of the Atiyah sequence 
\eqref{2.2.1} induces a $G$-equivariant 
section of \eqref{first_sequence}.
Thus, taking $L'=\Gamma_B(\operatorname{ad}(\xi))$, 
$L=\Gamma_B( \tau_Q/G)$, and $L''=\mathfrak{X}(B)$, 
we see that the chosen
connection determines a unique $\mathbf{e}$-connection $\omega\colon  L''\rightarrow L$ described above.
We note that we can identify
$L'$ with the $G$-invariant vertical vector
fields on $Q$, $L$ with the $G$-invariant vector fields on $Q$, and $L''$ with the $G$-invariant horizontal
vector fields on $Q$.

Recall that the curvature of the chosen connection,
thought of as a 2-form on $B$ with 
values in the adjoint bundle $\operatorname{ad}(\xi)$,
is the $C ^{\infty}(B)$-bilinear map 
$\Omega\colon  \mathfrak{X}(B)
\times \mathfrak{X}(B) \rightarrow 
\Gamma_B(\operatorname{ad}(\xi))$ of
$C^{\infty}(B)$-modules characterized by the identity
\[ 
\left[
[\operatorname{hor}]X, [\operatorname{hor}]Y\right]
= [\operatorname{hor}]\left([X,Y] \right) + 
\Omega(X, Y).
\quad \blacklozenge
\]
} \end{Example}
 
\begin{Theorem} \label{general} Given an extension $
\mathrm {\mathbf e} \colon  0 \longrightarrow L' \longrightarrow L
\longrightarrow L'' \longrightarrow 0$ of $(R,A)$-Lie algebras 
together with an 
$\mathrm {\mathbf e}$-connection $\omega\colon  L'' \to L$
and with  $\mathrm {\mathbf e}$-curvature
$\Omega\colon  L'' \otimes_AL''\to L'$ of $\omega$, 
the tautological Poisson structure
\begin{equation*}
\{\,\cdot\, ,\, \cdot\, \} \colon  \Sigm_A [L] \otimes \Sigm_A [L]
\longrightarrow \Sigm_A [L] \label{tauto}
\end{equation*}
on $\Sigm_A [L]$ {\rm (}defined in \eqref{3.16.2.1}, \eqref{3.16.2.2}, \eqref{3.16.2.3}{\rm )} is determined by the following identities:
\begin{align} \{\alpha',\beta'\} &= [\alpha',\beta']_{L'},\quad
\alpha',\beta' \in L', \label{3.16.2.11}
\\
\{\alpha'',\beta''\} &= [\alpha'',\beta'']_{L''} +
\Omega(\alpha'',\beta''),\quad \alpha'',\beta'' \in L'',
\label{3.16.2.12}
\\
\{\alpha',\beta''\} &= [\alpha',\beta''],\quad \alpha' \in L',
\,\beta'' \in L'', \label{3.16.2.111}
\\
\{\alpha'',\beta'\} &= [\alpha'',\beta'],\quad \alpha'' \in L'',
\,\beta' \in L', \label{3.16.2.112}
\\
\{\alpha'',a\} &= \alpha''(a) \in A,\quad a \in A,\,\alpha'' \in
L'', \label{3.16.2.22}
\\
\{\alpha',a\} &= 0 \in A,\quad a \in A,\,\alpha' \in L',
\label{3.16.2.222}
\\
\{u,vw\}&= \{u,v\}w + v\{u,w\},\quad u,v,w \in \Sigm_A[L].
\label{3.16.2.33}
\end{align}
\end{Theorem}

\begin{proof} This is an immediate consequence of the
identities \eqref{3.16.2.1}--\eqref{3.16.2.3} and  
\eqref{brack}. \end{proof}

We note that $\{a, b\} = 0$ for every $a,b \in A$.

\begin{Corollary} Under the hypotheses of Theorem {\rm
\ref{general}}, the extension $ \mathrm {\mathbf e}$ of
$(R,A)$-Lie algebras can be reconstructed from the Poisson
structure on $\Sigm_A[L]$ in the sense  that
the extension is determined by the relations
{\rm \eqref{3.16.2.11}}--{\rm \eqref{3.16.2.112}}.
\end{Corollary}

\begin{proof} This is an immediate consequence of Theorem
\ref{general} and Proposition \ref{prop3}.
\end{proof}

In particular, when we specialize  to the case where the ground ring $R$
is the field $\mathbb R$ of real numbers,
under the circumstances of Example {\rm {\ref{2.222}}}, the tautological
Poisson structure is determined by the identities 
given in Theorem \ref{general}.
Thus this description of the tautological Poisson bracket
{\eqref{tauto}} in Theorem {\rm {\ref{general}}} completely recovers
the Poisson structure on $(\mathrm T^*Q)\big/G$.

\section{Description of the Poisson structure in terms of differentials}
\label{desc}

Our ultimate goal is to reconcile our description of the Poisson structure 
on $(\mathrm T^*Q)\big/G$ 
that results from  Theorem \ref{general}
with that given in Theorem IV.1 of
\cite{zaalaone}. 
To this end, we will now develop a description, in terms of differentials, 
of the Poisson bracket given in Theorem \ref{general}.
As before, $A$ denotes
a commutative unital algebra over a general commutative unital ground ring $R$.

\subsection{Algebraic generalities.}
\label{Subsection_4.1}

We will momentarily abstract from the particular situation 
spelled out in 
Theorem \ref{general}
and consider a finitely
generated projective 
$A$-module $\ML$
together with an
$A$-module decomposition $\ML=\ML' \oplus \ML''$
into a direct sum of  two projective
$A$-modules $\ML'$ and $\ML''$.
We will simplify the exposition somewhat and use the notation
\[
\Sigm :=\Sigm_A[\ML],\quad \Sigm' := \Sigm_A[\ML'], \quad  \Sigm'' :=\Sigm_A[\ML''].
\]
The $A$-module decomposition of $\ML$ induces a canonical decomposition
$
\Sigm\cong \Sigm'\otimes_A \Sigm''
$
as the tensor product of the 
commutative $A$-algebras $\Sigm'$ and $\Sigm''$,
and we will also identify
$\Sigm$ with the symmetric
$\Sigm''$-algebra  
on $\Sigm''\otimes_A \ML'$
as well as
with the symmetric
$\Sigm'$-algebra 
on $\Sigm'\otimes_A \ML''$. 

In view of \eqref{fordif},
$\mathrm D_{\mathcal S'|A}\cong \mathcal S'\otimes_A \ML'$,
$\mathrm D_{\mathcal S''|A}\cong \mathcal S'\otimes_A \ML''$,
$\mathrm D_{\mathcal S|A}\cong \mathcal S\otimes_A \ML$,  
and the $\Sigm$-module 
$\mathrm D_{\mathcal S|A}$ of formal differentials
of $\Sigm$ over the algebra $A$
decomposes as the direct sum
\begin{equation}
\mathrm D_{\mathcal S|A}  
=\left(\Sigm \otimes _{\Sigm'}\mathrm D_{\mathcal S'|A}
\right) \oplus \left(\Sigm \otimes _{\Sigm''}
\mathrm D_{\mathcal S''|A}\right)
\cong 
\Sigm\otimes_A \ML' \oplus \Sigm\otimes_A \ML'' .
\label{dirs1}
\end{equation}
Hence
the exact sequence \eqref{dif2}
of $\mathcal S$-modules takes the form
\begin{equation}
0
\longrightarrow
\mathcal S \otimes _A \mathrm D_A \longrightarrow \mathrm
D_{\mathcal S} \longrightarrow 
\mathcal S \otimes _A \ML' \oplus \mathcal S \otimes _A \ML''
\longrightarrow 0 . \label{dif5}
\end{equation}
To make the basic constructions later in the paper more easily 
comprehensible for the reader, we will now 
somewhat refine the exact sequence \eqref{dif5}:

In view of \eqref{fordif},
$\mathrm D_{\mathcal S'|A}\cong \mathcal S'\otimes_A \ML'$
and $\mathrm D_{\mathcal{S}''|A}\cong 
\mathcal{S}''\otimes_A \ML''$, canonically. The inclusions
$R \subseteq A\subseteq \Sigm'$,
$R \subseteq A\subseteq \Sigm''$, 
$R \subseteq \Sigm'\subseteq \mathcal S$,
and $R \subseteq \Sigm''\subseteq \mathcal S$ 
determine exact sequences of $\mathcal{S}$-modules 
of the kind
\eqref{dif2}.

\begin{Proposition}
These exact sequence fit into the following
commutative  
diagram of $\mathcal{S}$-modules having
exact rows and columns:
\begin{equation}
\begin{CD}
@.
0
@.
0
@.
0
@.
\\
@.
@VVV
@VVV
@VVV
@.
\\
0
@>>>
\mathcal S \otimes_A \mathrm D_A 
@>>> 
\mathcal S \otimes_{\Sigm'}\mathrm D_{\mathcal S'} 
@>>> 
\mathcal S \otimes _A \ML'  
@>>> 
0 
\\
@.
@VVV
@VVV
@|
@.
\\
0
@>>>
\mathcal S \otimes _{\Sigm''} \mathrm D_{\Sigm''} @>>> \mathrm
D_{\mathcal S} @>>> \mathcal S \otimes _A \ML'  
@>>> 0
\\
@.
@VVV
@VVV
@VVV
@.
\\
0
@>>>
\mathcal S \otimes _A \ML''
@=
\mathcal S \otimes _A \ML''
@>>>
0
\\
@.
@VVV
@VVV
@.
@.
\\
@.
0
@.
0
@.
@.
\end{CD}
\label{CD1}
\end{equation}
\end{Proposition}

\begin{proof} These sequences are all of the kind \eqref{dif2}.
In view of the naturality of the constructions,
commutativity of the diagram is therefore immediate.
We leave the details to the reader.
\end{proof}

\subsection{The geometric situation.} \label{geomsit}
We 
will illustrate the diagram \eqref{CD1}
in terms of vector bundles.
Consider two vector bundles $p_{E_1}\colon  {E_1} \rightarrow B$ and $p_{E_2}\colon  {E_2} \rightarrow B$. Form the Whitney sum 
\[
p_{E_1}\oplus p_{E_2}\colon  {E_1}\oplus {E_2} \rightarrow B.
\]
 Denote by 
$p_{{E_1}^\ast}\colon  {E_1}^\ast \rightarrow B$, 
$p_{{E_2}^\ast}\colon  {E_2}^\ast \rightarrow B$, and
$p_{{E_1}^\ast}\oplus p_{{E_2}^\ast}\colon  {E_1}^\ast \oplus {E_2}^\ast  
\rightarrow B$ the corresponding dual bundles. When thought
of as the total space of a vector bundle over ${E_1}^\ast$ or 
as the total space of a vector bundle over
${E_2}^\ast$, 
instead of ${E_1} ^\ast \oplus {E_2} ^\ast$ we shall 
write $N:={E_1}^\ast \times_B {E_2}^\ast$ as in the following diagram:
\begin{equation*}
\begin{CD}
N @>{\cotan_{{E_1} ^\ast}}>> {E_1}^*
\\
@VV{\cotan_{{E_2} ^\ast}}V
@VV{p_{{E_1}^\ast}}V
\\
{E_2} ^\ast @>{p_{{E_2} ^\ast}}>> B
\end{CD}
\end{equation*}
Below we will use all three vector
bundle structures of $N$ over $B$, ${E_1}^\ast$, and
${E_2}^\ast$.

\begin{Proposition}
The following commutative diagram of vector bundles with
base $N$ has all rows and columns exact:
\begin{equation}
\begin{CD}
@.
0
@.
0
@.
0
@.
\\
@.
@AAA
@AAA
@AAA
@.
\\
0
@<<<
N \times_B \mathrm{T}B 
@<<<
N \times_{{E_1} ^\ast } \mathrm{T}({E_1} ^\ast ) 
@<<< 
N \times_B {E_1} ^\ast   
@<<< 
0 
\\
@.
@AAA
@AAA
@|
@.
\\
0
@<<<
N \times_{{E_2} ^\ast } \mathrm{T}({E_2} ^\ast) @<<< 
\mathrm{T}N@<<< N \times_B{E_1} ^\ast   
@<<< 0
\\
@.
@AAA
@AAA
@AAA
@.
\\
0
@<<<
N \times_B{E_2}^\ast
@= 
N \times_B{E_2}^\ast
@<<<
0
\\
@.
@AAA
@AAA
@.
@.
\\
@.
0
@.
0
@.
@.
\end{CD}
\label{CD_general_bundle}
\end{equation}
\end{Proposition}

\begin{proof}
We note first that the diagram is symmetric in ${E_1}$ and
${E_2}$; thus the first row and the first column as well as
the second row and the second column coincide when ${E_1}$
and ${E_2}$ are interchanged. Thus, it suffices to give 
the maps for the first two rows. Let 
$\lambda^{E_1}_b, \mu^{E_1}_b \in {E_1}^\ast$, 
$\lambda^{E_2}_b \in {E_2}^\ast$, so 
$(\lambda^{E_1}_b, \lambda^{E_2}_b) \in N$, and
$v_{\lambda^{E_1}_b} \in 
\mathrm{T}_{\lambda^{E_1}_b}{E_1}^\ast$, 
$v_{(\lambda^{E_1}_b, \lambda^{E_2}_b)} \in 
\mathrm{T}_{(\lambda^{E_1}_b, \lambda^{E_2}_b)}N$. 
The maps are: 
\begin{align*}
N\times_B {E_1} ^\ast \ni \left((\lambda^{E_1}_b, \lambda^{E_2}_b),
\mu^{E_1}_b\right)&\longmapsto \left(
(\lambda^{E_1}_b, \lambda^{E_2}_b), 
\operatorname{Ver}_{\lambda^{E_1}_b}(\mu^{E_1}_b) \right) \in 
N \times_{{E_1}^\ast}\mathrm{T}({E_1}^\ast)\\
N \times_{{E_1}^\ast}\mathrm{T}({E_1}^\ast) \ni \left(
(\lambda^{E_1}_b, \lambda^{E_2}_b), v_{\lambda^{E_1}_b} \right) 
&\longmapsto \left((\lambda^{E_1}_b, \lambda^{E_2}_b), 
T_{\lambda^{E_1}_b}p_{{E_1}^\ast}\big(v_{\lambda^{E_1}_b}\big)
\right) \in N \times_B \mathrm{T}B\\
N \times_B {E_1}^\ast \ni \left(
(\lambda^{E_1}_b, \lambda^{E_2}_b),
\mu^{E_1}_b\right)&\longmapsto 
\operatorname{Ver}_{(\lambda^{E_1}_b, \lambda^{E_2}_b)}
\left(\mu^{E_1}_b\right) \in \mathrm{T}N\\
\mathrm{T}N \ni v_{(\lambda^{E_1}_b, \lambda^{E_2}_b)}
& \longmapsto \left((\lambda^{E_1}_b, \lambda^{E_2}_b), 
\mathrm{T}_{(\lambda^{E_1}_b, \lambda^{E_2}_b)}\cotan_{{E_2}^\ast}
\big(v_{(\lambda^{E_1}_b, \lambda^{E_2}_b)}\big)\right) \in 
N \times_{{E_2}^\ast}\mathrm{T}({E_2}^\ast),
\end{align*}
where
\begin{align*}
\operatorname{Ver}_{\lambda^{E_1}_b}(\mu^{E_1}_b) &= 
\left.\frac{d}{dt}\right|_{t=0}\left(\lambda^{E_1}_b + t
\mu^{E_1}_b \right) \in \mathrm{T}_{\lambda^{E_1}_b}({E_1} ^\ast),\\ \operatorname{Ver}_{(\lambda^{E_1}_b, \lambda^{E_2}_b)}
\left(\mu^{E_1}_b\right) &= \left.\frac{d}{dt}\right|_{t=0}
\left(\lambda^{E_1}_b + t\mu^{E_1}_b, \lambda^{E_2}_b \right)
\in \mathrm{T}_{(\lambda^{E_1}_b, \lambda^{E_2}_b)}N.
\end{align*}
The commutativity of the diagram and the exactness of 
the rows and columns is an easy verification.
\end{proof}

Dualizing diagram \eqref{CD_general_bundle} we obtain
the  commutative diagram 
\begin{equation}
\begin{CD}
@.
0
@.
0
@.
0
@.
\\
@.
@VVV
@VVV
@VVV
@.
\\
0
@>>>
N \times_B \mathrm{T}^\ast B 
@>>>
N \times_{{E_1} ^\ast } \mathrm{T}^\ast ({E_1} ^\ast ) 
@>>>
N \times_B {E_1}   
@>>> 
0 
\\
@.
@VVV
@VVV
@|
@.
\\
0
@>>>
N \times_{{E_2} ^\ast } \mathrm{T}^\ast ({E_2} ^\ast) @>>> 
\mathrm{T}^\ast N@>>> N \times_B{E_1}   
@>>> 0
\\
@.
@VVV
@VVV
@VVV
@.
\\
0
@>>>
N \times_B{E_2}
@= 
N \times_B{E_2}
@>>>
0
\\
@.
@VVV
@VVV
@.
@.
\\
@.
0
@.
0
@.
@.
\end{CD}
\label{CD_general_bundle_dual}
\end{equation}
of vector bundles over 
$N$ having exact rows and columns.

We next show that the diagram that we obtain when we take 
sections in \eqref{CD_general_bundle_dual} is the 
$C^{\infty}$-analog of \eqref{CD1}. 
We will use the notation $\Gamma$ for the smooth sections functor.
Recall
the familiar Serre-Swan equivalence:
Given a smooth (paracompact) manifold $M$,
the assignment to a smooth vector bundle $\zeta$ on $M$ of
its $C^{\infty}(M)$-module $\Gamma(\zeta)$
of smooth sections
is an equivalence of categories
between smooth vector bundles on $M$ and
finitely generated projective  $C^{\infty}(M)$-modules. 

\begin{Lemma}
\label{section_lemma}
For vector bundles $\mathscr{{E}}_1 \rightarrow M$ and 
$\mathscr{{E}}_2 \rightarrow M$, 
the fiber product $\mathscr{{E}}_1 \times_M \mathscr{{E}}_2 
\rightarrow \mathscr{{E}}_1$ being
regarded as a vector bundle over
$\mathscr{{E}}_1$,
the Serre-Swan equivalence 
yields a canonical isomorphism 
\[
\Gamma_{\mathscr{{E}}_1}(\mathscr{{E}}_1 \times_M 
\mathscr{{E}}_2) \cong C^{\infty}(\mathscr{{E}}_1)
\otimes_{C^{\infty}(M)}\Gamma_M(\mathscr{{E}}_2)
\]
of  $C^{\infty}(\mathscr{{E}}_1)$-modules.
\end{Lemma}

For a proof, the reader may consult, e.~g.,
\cite{greubetal} (2.26).

Applying the lemma to the exact sequence that we obtain
by taking sections in \eqref{CD_general_bundle_dual}, 
we get the following commutative diagram of 
$C^{\infty}(N)$-modules
\begin{equation}
\label{section_diagram}
\footnotesize{
\begin{CD}
@.
0
@.
0
@.
0
@.
\\
@.
@VVV
@VVV
@VVV
@.
\\
0
@>>>
C^{\infty}(N)\otimes_{C ^{\infty}(B)} 
\Omega^1(B) 
@>>>
C^{\infty}(N) \otimes_{C^{\infty}({E_1}^\ast)} 
\Omega^1({E_1} ^\ast)
@>>>
C^{\infty}(N) \otimes_{C^{\infty}(B)} \Gamma_B({E_1})   
@>>> 
0 
\\
@.
@VVV
@VVV
@|
@.
\\
0
@>>>
C^{\infty}(N) \otimes_{C^{\infty}({E_2}^\ast)} 
\Omega^1({E_2} ^\ast)
@>>> \Omega^1(N)
@>>> C^{\infty}(N) \otimes_{C^{\infty}(B)} 
\Gamma_B({E_1})   
@>>> 0 \\
@.
@VVV
@VVV
@VVV
@.
\\
0
@>>>
C^{\infty}(N)\otimes_{C^{\infty}(B)}\Gamma_B({E_2})
@= 
C^{\infty}(N)\otimes_{C^{\infty}(B)}\Gamma_B({E_2})
@>>>
0
\\
@.
@VVV
@VVV
@.
@.
\\
@.
0
@.
0
@.
@.
\end{CD} }
\end{equation}
with all rows and columns exact. This
is the precise geometric analog \eqref{CD1},
with $(R, A)=(\mathbb{R}, C^{\infty}(B))$,
$\ML' = \Gamma_B({E_1})$, $\ML''= \Gamma_B({E_2})$,
the algebras $C^{\infty}({E_1}^\ast)$,
$C^{\infty}({E_2}^\ast)$, and $C^{\infty}(N)$
being the
$C ^{\infty}$-completions of the respective rings
$\mathcal{S}'$, $\mathcal{S}''$ and $\mathcal{S}$.

\subsection{Basic example}
\label{illust}
We will study the diagram \eqref{CD_general_bundle_dual} 
in the context of Example \ref{2.222}. Let
$\xi\colon  Q \rightarrow B$ be a right principal 
$G$-bundle and $\widetilde \prin\colon  \widetilde{Q} =
Q \times_{B} \mathrm{T}^\ast B\rightarrow \mathrm{T}^*B$
its pull back relative to the cotangent bundle projection
$\cotan_B\colon \mathrm T^*B \rightarrow B$. Thus, in the 
commutative diagram
\begin{equation}
\begin{CD}
\widetilde{Q} @>{\widetilde \prin}>> \mathrm{T}^*B
\\
@VV{\xi^*(\cotan_B)}V
@VV{\cotan_B}V
\\
{Q} @>{\prin}>> B,
\end{CD}
\label{cd9}
\end{equation} 
the vertical projection map 
$\xi^*(\cotan_B)\colon  \widetilde{Q} \to Q$
is a vector bundle projection with fibers $\mathrm{T}^\ast _bB$ 
($b \in B$)
and 
$\widetilde{\xi}\colon   \widetilde{Q} \rightarrow 
\mathrm{T}^\ast B$ is a principal $G$-bundle.

Consider the tangent bundles $\tau_Q\colon \mathrm TQ \to Q$ and $\tau_B\colon  \mathrm TB \to B$
of $Q$ and $B$, respectively, as well as the induced vector bundle $\xi^*(\tau_B)\colon 
Q \times_B \mathrm TB \to Q$ 
on $Q$, the pull back of $\tau_B$ via $\xi$.
Viewed as a map over $Q$, the 
$G$-equivariant surjective map 
$\mathrm{T}Q \ni v_q\longmapsto (q, \mathrm{T}_q\xi(v_q)) 
\in  Q \times _B \mathrm TB$ ($q \in Q$)
is a morphism of vector bundles on $Q$, spelled out here for the total spaces.

\subsubsection{The Sternberg space}
\label{Sternberg_space}
The total space of the coadjoint bundle 
$\operatorname{ad}^\ast(\widetilde{\xi})\colon 
N_{\mathrm{St}}\colon =\widetilde Q \times_G \mathfrak g^*
\rightarrow \mathrm{T}^\ast B$ is occasionally referred 
to in the literature as the {\em Sternberg\/} space 
(associated to the data).

We will now apply the construction in Subsection
\ref{geomsit} above: 
We take the vector bundles 
$p_{E_1}\colon {E_1} \rightarrow B$ 
and $p_{E_2}\colon {E_2} \rightarrow B$
to be the adjoint
bundle $\operatorname{ad}(\xi)\colon  Q \times _G \mathfrak{g}
\rightarrow B$ of $\xi\colon  Q \rightarrow B$ and 
the tangent
bundle $\tau_B\colon \mathrm{T}B \rightarrow B$ of $B$, respectively. 
Then
$N\colon = (Q \times_G\mathfrak{g}^\ast) \oplus \mathrm{T}^\ast 
B$ is not only the total space of a vector bundle over $B$
but the natural projections $N \rightarrow 
Q \times _G\mathfrak{g}^\ast$ and $N \rightarrow 
\mathrm{T}^\ast B$ also define vector bundles. In what
follows we shall use all three vector bundle structures.
We note that
\begin{equation}
\label{second_space}
N\ni ([q, \mu], \alpha_{[q]})\longleftrightarrow
[(\alpha_{[q]}, q), \mu] \in N_{\rm St},
\ 
q \in Q,\ \mu \in \mathfrak g^*,\ \alpha_{[q]} \in \mathrm T^*_{[q]}B,
\end{equation}
yields a canonical isomorphism of vector bundles 
from $N \rightarrow \mathrm{T}^\ast B$ to $N_{\rm St} \rightarrow 
\mathrm{T}^\ast B$. In view of this isomorphism, 
in the situation under discussion,
the diagrams 
\eqref{CD_general_bundle_dual} and 
\eqref{section_diagram} now specialize to
\begin{equation}
\begin{CD}
@.
0
@.
0
@.
0
@.
\\
@.
@VVV
@VVV
@VVV
@.
\\
0
@>>>
N_{\rm St} \times _B \mathrm T^*B
@>>> 
N_{\rm St}\times_{Q\times_G\mathfrak g^*}\mathrm T^*(Q\times_G \mathfrak g^*)
@>>> 
N_{\rm St}\times_{B}(Q\times_G \mathfrak g)  
@>>> 
0 
\\
@.
@VVV
@VVV
@|
@.
\\
0
@>>>
N_{\rm St}\times_{\mathrm T^* B} \mathrm T^* \mathrm T^* B
@>>> 
\mathrm T^* N_{\rm St} 
@>>> 
N_{\rm St}\times_{B}(Q\times_G \mathfrak g)
@>>> 0
\\
@.
@VVV
@VVV
@VVV
@.
\\
0
@>>>
N_{\rm St} \times_B \mathrm TB
@=
N_{\rm St} \times_B \mathrm TB
@>>>
0
\\
@.
@VVV
@VVV
@.
@.
\\
@.
0
@.
0
@.
@.
\end{CD}
\label{CD2}
\end{equation}
and
\begin{equation}
\!\!\!\!\!\!\!\!\!\!\!\!
\label{CD2_sections}
\footnotesize{
\begin{xy}
\xymatrix{
&0\ar[d]&0\ar[d]&0\ar[d]&\\
0\ar[r]&C^{\infty}(N_{\rm St})\otimes_{C ^{\infty}
(B)}\Omega^1(B)\ar[r]\ar[d]&C^{\infty}(N_{\rm St})\otimes_{C ^{\infty}(Q\times_G\mathfrak g^*)}\Omega^1(Q\times_G \mathfrak g^*)
\ar[r]\ar[d]&
C^{\infty}(N_{\rm St})\otimes_{C ^{\infty}(B)}\Gamma_B(Q\times_G \mathfrak{g})\ar[r]\ar@{=}[d]&0\\
0\ar[r]&C^{\infty}(N_{\rm St})\otimes_{C ^{\infty}(\mathrm{T}^* B)}\Omega^1(\mathrm T^* B) \ar[r]\ar[d]&
\Omega^1(N_{\rm St})\ar[r]\ar[d]&
C^{\infty}(N_{\rm St})\otimes_{C ^{\infty}(B)}\Gamma_B(Q\times_G \mathfrak{g}) \ar[r]\ar[d]&0\\
0\ar[r]&C^{\infty}(N_{\rm St})\otimes_{C ^{\infty}(B)}
\mathfrak{X}(B)\ar@{=}[r]\ar[d]&
C^{\infty}(N_{\rm St})\otimes_{C ^{\infty}(B)}
\mathfrak{X}(B)\ar[r]\ar[d]&0&&\\
&0&0&&&&
}
\end{xy} 
}
\end{equation}
The last diagram 
is the $C^{\infty}$-analog of diagram
\eqref{CD1} in the present special case. 

To be precise, $(R, A)=(\mathbb{R}, C^{\infty}(B))$,
the algebra $C^{\infty}(N_{\rm St})$
of ordinary smooth functions on $N_{\rm St}$ is
substituted for the algebra $\mathcal{S}$,
as an $A$-module, $\ML' = \Gamma_B({Q\times_G\mathfrak g})$
is the space of sections of the adjoint bundle
$\operatorname{ad}(\xi)\colon  Q \times_G
\mathfrak{g}\rightarrow B$ on $B$,
as an $A$-module, $\ML''= \Gamma_B({\mathrm TB})= \mathfrak{X}(B)$,
the smooth vector fields on $B$,
the algebra $\mathcal{S}''$ is the algebra of smooth functions
on $\mathrm{T}^* B$ that are polynomial on the fibers of the
cotangent bundle $\mathrm{T}^\ast B \to B$ of $B$.
Since \eqref{cd9} is a pull back diagram, the obvious projection
$\widetilde Q \times_G \mathfrak g \to \mathrm{T}^\ast B$
is a vector bundle on $\mathrm{T}^\ast B$; hence
the $\mathcal{S}''$-module
$\mathcal{S}''\otimes_A \ML'$ is naturally isomorphic to
the space of sections of the adjoint bundle $\mathrm{ad}(\widetilde \xi)
\colon  \widetilde Q \times_G \mathfrak g
\longrightarrow \mathrm{T}^* B$ of $\widetilde{\xi}$
that are polynomial on the fibers of the cotangent bundle
$\mathrm{T}^\ast B\to B$ of $B$. The algebras
$C^{\infty}({Q\times_G\mathfrak g}^\ast)$,
$C^{\infty}(\mathrm{T}^\ast B)$, and $C^{\infty}(N_{\rm St})$
are the $C ^{\infty}$-completions of the respective rings
$\mathcal{S}'$, $\mathcal{S}''$ and $\mathcal{S}$.
Diagram \eqref{CD2_sections} now arises from \eqref{CD2} by
the operation of taking smooth sections
of the corresponding vector bundles on $N_{\rm St}$.

\subsubsection{The Weinstein space}
\label{Weinstein_space}
The smooth manifold $N_{\mathrm W}:=(\mathrm{T}^\ast Q)/G$ is 
occasionally 
referred to in the literature as
the \textit{Weinstein space};
it inherits a Poisson structure from the ordinary symplectic
Poisson structure of $\mathrm{T}^\ast Q$.
 The projection 
$(\mathrm{T}^\ast Q)/G \ni [ \alpha_q] \mapsto [q] 
\in Q/G = B$ ($q \in Q$) turns  $N_{\mathrm W}$ into a vector bundle over $B$.

Viewed as a map over $Q$, the 
obvious 
$G$-equivariant surjective map
$
\mathrm{T}Q \longrightarrow Q \times _B \mathrm TB
$
is a morphism of vector bundles on $Q$, spelled out here for the total spaces.
A principal connection for $\xi$ is precisely a $G$-equivariant vector bundle section for this
surjection. 
Thus, pick a principal connection for $\prin$. Taking duals in the fiber directions, 
we then obtain a
$G$-equivariant map
\begin{equation*}
\mathrm T^*Q \longrightarrow   Q \times _B \mathrm T^*B=\widetilde{Q},
\end{equation*}
indeed a morphism of $G$-vector bundles on $Q$,
spelled out here for the total spaces.
This map combines with the $G$-momentum mapping
$\mathrm T^*Q \to \mathfrak g^*$ to 
yield
a $G$-equivariant 
diffeomorphism
\begin{equation}
\mathrm T^*Q \longrightarrow  \widetilde Q \times \mathfrak g^*;
\label{gequiv}
\end{equation}
in particular, the $G$-equivariant projection $\mathrm T^*Q \to \widetilde Q$
to the first component
is now seen to be the projection of a $G$-vector bundle on $\widetilde Q$
which, as an ordinary vector bundle on $\widetilde Q$, is trivial (not equivariantly trivial).
Thus the diffeomorphism \eqref{gequiv}, in turn, descends to a diffeomorphism
\begin{equation}
N_W=(\mathrm T^*Q)\big/G \longrightarrow  \widetilde Q \times_G \mathfrak g^* =N_{\mathrm{St}}. 
\label{idstw}
\end{equation}
This identifies $N_{\mathrm W}$ with $\mathrm{T}^\ast(Q/G)\oplus
(Q \times_G\mathfrak{g}^\ast)$, and the latter space acquires, as we have
seen in \eqref{second_space},  a vector bundle structure 
over each of the three base spaces $B$, $Q \times _G\mathfrak{g}^\ast$,
and $\mathrm{T}^\ast B$ in a natural way. One of the main goals of this
paper, carried out in 
Subsection \ref{principal},
is to express the natural quotient Poisson bracket
on $(\mathrm{T}^\ast Q)/G$ via this vector bundle 
isomorphism over $B$ on the space $\mathrm{T}^\ast(Q/G)
\oplus(Q \times_G\mathfrak{g}^\ast) \cong N_{\rm St}$.
In order to do this, we step back and consider the general situation presented in Subsection 
\ref{Subsection_4.1}. 

\subsection{Linear connections}
\label{lincon}

Consider a vector bundle $\lambda\colon  \mathcal{L} 
\to B$; let $V=
\mathcal{L}\times_B\mathcal{L} \rightarrow \mathcal{L}$
denote its vertical subbundle. 
Recall that a {\em linear\/} or {\em Ehresmann\/} 
connection on 
$\lambda$  is given by a section
$\mathcal L \times _B\mathrm TB \to \mathrm T\mathcal L$
for the associated extension
\begin{equation}
0
\longrightarrow
V
\longrightarrow
\mathrm T\mathcal L
\longrightarrow
\mathcal L \times _B\mathrm TB
\longrightarrow
0
\label{ehres}
\end{equation}
of vector bundles on $\mathcal L$ or, equivalenty, by a retraction
$\mathrm T^*\mathcal L  \to \mathcal L \times _B\mathrm T^*B $
for the extension
\begin{equation*}
0
\longleftarrow
V^*
\longleftarrow
\mathrm T^*\mathcal L
\longleftarrow
\mathcal L \times _B\mathrm T^*B
\longleftarrow
0
\end{equation*}
of the dual vector bundles on $\mathcal L$.

In the algebraic setting, 
accordingly, 
 let $\Mo$ be a projective $A$-module, and let
$\Sigm_A[\Mo]$ be the symmetric $A$-algebra on $\Mo$.
We define
a {\em linear connection\/} on the $A$-module $\Mo$ to be a
 retraction
\begin{equation}
r \colon
\mathrm D_{\Sigm_A[\Mo]}
\longrightarrow
\Sigm_A[\Mo]
\otimes_A\mathrm D_A
\label{ret}
\end{equation}
for the extension
\begin{equation}
0
\longrightarrow
\Sigm_A[\Mo]
\otimes_A\mathrm D_A
\longrightarrow
\mathrm D_{\Sigm_A[\Mo]}
\longrightarrow
\Sigm_A[\Mo]
\otimes_A\Mo
\longrightarrow
0
\label{ehres2}
\end{equation}
in the category of $\Sigm_A[\Mo]$-modules.
Using the canonical isomorphism
\[
\mathrm{Hom}_{\Sigm_A[\Mo]}(\Sigm_A[\Mo]
\otimes_A\Mo,\Sigm_A[\Mo])
\longrightarrow
\operatorname {Der}({\Sigm_A[\Mo]|A})
\]
of $\Sigm_A[\Mo]$-modules, cf. \eqref{dfordif},
we see that
a linear connection on $\Mo$ determines
a section 
\begin{equation}
s\colon \Sigm_A[\Mo]
\otimes_A\operatorname {Der}(A)
\longrightarrow
\operatorname {Der}({\Sigm_A[\Mo]})
\end{equation}
for the dual extension
\begin{equation}
0
\longrightarrow
\operatorname {Der}({\Sigm_A[\Mo]|A})
\longrightarrow
\operatorname {Der}({\Sigm_A[\Mo]})
\longrightarrow
\Sigm_A[\Mo]
\otimes_A\operatorname {Der}(A)
\longrightarrow
0
\label{ehres3}
\end{equation}
of $\Sigm_A[\Mo]$-modules
(but, beware, the retraction $r$ is not determined
by the section $s$
unless 
$D_A$ is reflexive, i.e., unless
the canonical
$A$-module morphism from $\mathrm D_A$ to its double $A$-dual
 $\mathrm D_A^{**}$
is an isomorphism).
A linear connection,
i.~e. retraction $r$
of the kind \eqref{ret},
 for $\Mo$ induces the 
direct sum decomposition
\begin{equation}
\mathrm D_{\Sigm_A[\Mo]} = \left(\Sigm_A[\Mo]
\otimes_{A}\mathrm D_{A}\right) \oplus
\mathrm{ker}(r) 
\label{dirs121}
\end{equation}
in the category of $\Sigm_A[\Mo]$-modules, and the projection to
$\Sigm_A[\Mo]\otimes_A\Mo$ restricts to an isomorphism
\[
\mathrm{ker}(r) \longrightarrow \Sigm_A[\Mo]\otimes_A\Mo
\]
of 
$\Sigm_A[\Mo]$-modules.
We 
will use the notation
\[
(\mathrm D_{\Sigm_A[\Mo]})_{\mathrm {vert}}
=\mathrm{ker}(r),
\quad
(\mathrm D_{\Sigm_A[\Mo]})_{\mathrm {hor}}=\Sigm_A[\Mo]\otimes_A
\mathrm D_A, 
\]
and we will refer to
$(\mathrm D_{\Sigm_A[\Mo]})_{\mathrm {vert}}$ 
and $(\mathrm D_{\Sigm_A[\Mo]})_{\mathrm {hor}}$ as the 
{\em vertical\/} and {\em hori\-zon\-tal\/}
components, respectively, of  the resulting decomposition \eqref{dirs121}
of $\mathrm D_{\Sigm_A[\Mo]}$, so that
this decomposition takes the form
\begin{equation}
\mathrm D_{\Sigm_A[\Mo]} =(\mathrm D_{\Sigm_A[\Mo]})_{\mathrm {hor}} 
\oplus (\mathrm D_{\Sigm_A[\Mo]})_{\mathrm {vert}} .
\label{dirs91}
\end{equation}
Since $\Mo$ is supposed to be projective as an $A$-module,
 a linear connection on $\Mo$ exists.

We now return to the situation of Subsection 
\ref{Subsection_4.1} above. Thus $\ML'$ and $\ML''$
are projective $A$-modules, $\ML=\ML' \oplus \ML''$,
$\mathcal{S}= \mathcal{S}_A[\ML]$,
$\mathcal{S'}= \mathcal{S}_A[\ML']$,
$\mathcal{S''}= \mathcal{S}_A[\ML'']$.
Let $r$ be a linear connection,
i.~e. retraction $r\colon
\mathrm D_{\Sigm'}
\longrightarrow
\Sigm'
\otimes_A\mathrm D_A
$ 
of the kind \eqref{ret},
for $\ML'$; such a retraction exists since $\ML'$ is projective as 
an $A$-module. Let
\begin{equation}
\mathrm D_{\Sigm'}  
=(\mathrm D_{\Sigm'})_{\mathrm {hor}} 
\oplus (\mathrm D_{\Sigm'})_{\mathrm {vert}}=
\left(\Sigm'
\otimes_{A}\mathrm D_{A}\right) \oplus
\mathrm{ker}(r) 
\label{dirs12}
\end{equation}
be the corresponding 
direct sum decomposition \eqref{dirs121}
in the category of $\Sigm'$-modules.
Then
\[
\Sigm\otimes_{\Sigm'}r\colon
\Sigm\otimes_{\Sigm'}\mathrm D_{\Sigm'}
\longrightarrow
\Sigm
\otimes_A\mathrm D_A
\]
is a retraction for the upper horizontal exact sequence 
\begin{equation}
\begin{CD}
0
@>>>
\mathcal S \otimes_A \mathrm D_A 
@>>> 
\mathcal S \otimes_{\Sigm'}\mathrm D_{\mathcal S'} 
@>>> 
\mathcal S \otimes _A \ML'  
@>>> 
0 
\end{CD}
\end{equation}
in \eqref{CD1},
and 
\[
\mathrm{ker}(\Sigm\otimes_{\Sigm'}r)
=
\Sigm\otimes_{\Sigm'}\mathrm{ker}(r).
\]
Moreover, since the projection 
from $\mathrm D_{\Sigm'}$ to
$\Sigm'\otimes_A\ML'$ restricts to an isomorphism
$
\mathrm{ker}(r) \longrightarrow \Sigm'\otimes_A \ML'
$
of 
$\Sigm'$-modules,
the projection 
from 
$\Sigm\otimes_{\Sigm'}\mathrm D_{\Sigm'}$ to
$\Sigm\otimes_A\ML'$ 
restricts to an isomorphism 
\begin{equation}
\Sigm\otimes_{\Sigm'}\mathrm{ker}(r) \longrightarrow 
\Sigm\otimes_A\ML'
\label{iso1}
\end{equation}
of 
$\Sigm$-modules. Hence, 
when we view
$\Sigm\otimes_{\Sigm'}\mathrm{ker}(r)$
as an $\Sigm$-submodule of
$\mathrm D_{\Sigm}$
via the injection 
$\Sigm\otimes_{\Sigm'} \mathrm D_{\Sigm'}\to \mathrm D_{\Sigm}$
(arising as an arrow of the middle row of \eqref{CD1}),
the isomorphism \eqref{iso1} is also
the restriction of
the projection from
$\mathrm D_{\Sigm}$ to
$\Sigm\otimes_A\ML'$
(arising  as an arrow of the middle column of \eqref{CD1}).
Consequently the chosen linear connection $r$ for $\ML'$ induces corresponding 
direct sum decompositions of
$\Sigm\otimes_{\Sigm'}\mathrm D_{\Sigm'}$ and $\mathrm D_{\Sigm}$
into the following vertical and  horizontal components:
\begin{align*}
(\Sigm\otimes_{\Sigm'}\mathrm D_{\Sigm'})_{\mathrm {vert}}
&=
(\mathrm D_{\Sigm})_{\mathrm {vert}}
=\Sigm\otimes_{\Sigm'}\mathrm{ker}(r),
\\
(\Sigm\otimes_{\Sigm'}\mathrm D_{\Sigm'})_{\mathrm {hor}}
&= \Sigm\otimes_A\mathrm D_A,\ 
(\mathrm D_{\Sigm})_{\mathrm {hor}}=\Sigm \otimes_{\Sigm''}
\mathrm D_{\Sigm''},
\end{align*}
and we will refer to
$(\Sigm\otimes_{\Sigm'}\mathrm D_{\Sigm'})_{\mathrm {vert}}
$ 
as the {\em vertical\/}
component of $(\Sigm\otimes_{\Sigm'}\mathrm D_{\Sigm'})$,
to
$(\Sigm\otimes_{\Sigm'}\mathrm D_{\Sigm'})_{\mathrm {hor}}$
as its {\em horizontal\/}
component, to
$(\mathrm D_{\Sigm})_{\mathrm {vert}}$ 
as the {\em vertical\/}
component of $\mathrm D_{\Sigm}$, and to
$(\mathrm D_{\Sigm})_{\mathrm {hor}}$ as its {\em horizontal\/} component.
Thus the chosen linear connection $r$ for $\ML'$ induces the 
direct sum decompositions
\begin{align}
\Sigm\otimes_{\Sigm'} \mathrm D_{\Sigm'} &= 
\left(\Sigm\otimes_{\Sigm'} \mathrm D_{\Sigm'}\right)_{\mathrm {hor}}
\oplus
\left(\Sigm\otimes_{\Sigm'} \mathrm D_{\Sigm'}\right)_{\mathrm {vert}}
\label{dirs122}
\\
\mathrm D_{\Sigm} &=(\mathrm D_{\Sigm})_{\mathrm {hor}} 
\oplus (\mathrm D_{\Sigm})_{\mathrm {vert}} .
\label{dirs9}
\end{align}
Further, the projection 
from $\Sigm\otimes_{\Sigm'}\mathrm D_{\Sigm'}$ 
to
$\Sigm\otimes_A\ML'$ 
(in the upper horizontal row of \eqref{CD1})
restricts to an isomorphism
$\left(\Sigm\otimes_{\Sigm'} \mathrm D_{\Sigm'}\right)_{\mathrm {vert}}
\to \Sigm\otimes_A\ML'$
of 
$\Sigm$-modules
and
the projection 
from $\mathrm D_{\Sigm}$ 
to
$\Sigm\otimes_A\ML'$ 
(in the middle horizontal row of \eqref{CD1})
restricts to an isomorphism
$\left(\mathrm D_{\Sigm}\right)_{\mathrm {vert}}
\to \Sigm\otimes_A\ML'$
of 
$\Sigm$-modules;
these isomorphisms are precisely
the isomorphism
\eqref{iso1}
of 
$\Sigm$-modules.
Thus, all told, somewhat more formally,
the resulting retraction
\[
\mathrm D_{\Sigm} \longrightarrow (\mathrm D_{\Sigm})_{\mathrm {hor}} 
=\Sigm \otimes_{\Sigm''}
\mathrm D_{\Sigm''}
\]
for the middle horizontal exact sequence 
\begin{equation}
\begin{CD}
0
@>>>
\Sigm \otimes_{\Sigm''}
\mathrm D_{\Sigm''} 
@>>> 
\mathrm D_{\Sigm} 
@>>> 
\mathcal S \otimes _A \ML'  
@>>> 
0 
\end{CD}
\label{sigma_connection}
\end{equation}
in \eqref{CD1}
is 
a linear connection on the induced
$\Sigm''$-module $\Sigm''\otimes_A\ML'$, and the notation is consistent.

By construction,  
in view of the direct sum decomposition $\ML=\ML'\oplus \ML''$,
any $\alpha' \in \ML'$ is a 
\lq\lq vertical\rq\rq\ element of $\Sigm$
in an obvious sense;
however
the differential $d\alpha'$
refers to the differential of
$\alpha'$  viewed as a member of $\Sigm$,
taken as an $R$-algebra, i.~e.  over the ground ring $R$
(rather than as an $A$-algebra).
Hence,
in general, given $\alpha' \in \ML'$,
the differential $d\alpha'\in \mathrm D_{\Sigm}$ has, 
beyond a vertical component
$(d\alpha')_{\mathrm {vert}}\in (\mathrm D_{\Sigm})_{\mathrm {vert}}$,
a horizontal component
$(d\alpha')_{\mathrm {hor}}\in (\mathrm D_{\Sigm})_{\mathrm {hor}}$;
now, 
$d\alpha'$ is actually the image in $\mathrm D_{\Sigm}$ 
of the differential of $\alpha'$ in $\mathrm D_{\Sigm'}$ 
and,
under 
the projection
$\mathrm
D_{\mathcal S} \longrightarrow 
\mathcal S \otimes _A \ML' \oplus \mathcal S \otimes _A \ML''
$ in \eqref{dif5}, 
the horizontal component
$(d\alpha')_{\mathrm {hor}}
$ goes to zero whence
$(d\alpha')_{\mathrm {hor}}$ is a member of
$\Sigm\otimes_A \mathrm D_A$.

\begin{Remark}
\label{rem4}{\rm 
At the present stage, the theory is symmetric in $\ML'$ and $\ML''$,
that is, we could have chosen a linear connection for $\ML''$ as well
and developed the same kind of decomposition of
$\mathrm D_{\Sigm}$ as {\rm\eqref{dirs9}}.} \quad $\blacklozenge$
\end{Remark}

\begin{Example}
\label{extr}
{\rm Consider the trivial
principal
bundle $Q \to B$ where $Q=B \times G$, endowed with the trivial
connection so that, with the notation established earlier,
$N=(\mathrm T^*B) \times \mathfrak g^*$.
Let $f$ be a function on $\mathrm T^*B$, and let $X \in \mathfrak g$, 
viewed as a linear function on $\mathfrak g^*$;
view 
$X$, $f$ and
$f X$ as functions on $N$ in the obvious way.
The ordinary differential $\mathbf{d}(f X)$ 
of the function $fX$
decomposes as 
\[
\mathbf{d}(f X)=(\mathbf{d}f)X + f \mathbf{d}X,\quad (\mathbf{d}(f X))_{\mathrm{vert}}=
f \mathbf{d}X,\ 
(\mathbf{d}(f X))_{\mathrm{hor}}=(\mathbf{d}f)X.
\]
Here the notation 
$\mathbf{d}f$ 
refers to the ordinary differential
of $f$, 
viewed
as a function on $N=(\mathrm T^*B) \times \mathfrak g^*$ that depends only on
the first variable, and
$\mathbf{d}X$ refers to the ordinary differential
of $X$, viewed 
as a function on $N=(\mathrm T^*B) \times \mathfrak g^*$ that depends only on
the second variable.
}\end{Example}

\subsection{The Poisson structure}
We return to the situation of Section 
\ref{sec_extensions} above. 
Thus 
\[
\mathbf e\colon 
0
\longrightarrow
L'
\longrightarrow
L
\longrightarrow
L''
\longrightarrow
0
\]
is an extension of $(R,A)$-Lie algebras.
Choose an $\mathbf e$-connection for the extension
$\mathbf e$, i.e., an $A$-module section 
$\omega\colon  L'' \to L$.
The connection $\omega$ 
determines (and is determined by) an $A$-module decomposition
$L=L' \oplus L''$; here we identify $L''$ with its image
$\omega (L'')$ in $L$, with a slight abuse of  notation.
We also choose a linear connection 
$r\colon \mathrm{D}_{\mathcal{S}'} \rightarrow \mathcal{S}' 
\otimes \mathrm{D}_A$ on the $A$-module $L'$ 
(see \eqref{ret}) and hence, by \eqref{dirs9}, 
$\mathrm{D}_\mathcal{S}=
(\mathrm{D}_{\Sigm})_{\mathrm{hor}} \oplus 
(\mathrm{D}_{\Sigm})_{\mathrm{vert}}$, as $\mathcal{S}$-modules, where
\[
(\mathrm{D}_{\Sigm})_{\mathrm{vert}}
=\Sigm\otimes_{\Sigm'}\mathrm{ker}(r)
\cong
\Sigm\otimes_A L', \qquad
(\mathrm D_{\Sigm})_{\mathrm{hor}}\colon =
\Sigm \otimes_{\Sigm''}
\mathrm D_{\Sigm''}. 
\]
We will now apply the material developed before,
with $\ML$, $\ML'$,$\ML''$ the $A$-modules
that underlie, respectively, $L$, $L'$, and $L''$ and,
to simplify the notation, we will denote these
underlying $A$-modules by $L$, $L'$, and $L''$ as well.

To develop the promised description, in terms of differentials, of the
Poisson bracket given in Theorem \ref{general},
all we have to do now is to evaluate the formulas
\eqref{dif3} and \eqref{dif4} in terms of the direct sum decomposition
\eqref{dirs9}.
The algebras $\Sigm'=\Sigm_A[L']$ and
$\Sigm''=\Sigm_A[L'']$ carry the tautological Poisson structures
that correspond to the 
$A$-Lie algebra structure of $L'$ and to
the Lie-Rinehart structure of 
$(A,L'')$, 
respectively, and the corresponding
Poisson 2-form \eqref{3.5.1}
is defined for each of these algebras.
We denote these 2-forms by
\[
\pi'\colon  \mathrm D_{\Sigm'|A} \otimes_{\Sigm'}\mathrm D_{\Sigm'|A} 
\longrightarrow
\Sigm',
\ 
\pi''\colon  \mathrm D_{\Sigm''} \otimes_{\Sigm''}\mathrm D_{\Sigm''} \longrightarrow
\Sigm''
\]
and, relative to the decomposition
$\mathrm D_{\Sigm} =(\mathrm D_{\Sigm})_{\mathrm{hor}} 
\oplus (\mathrm D_{\Sigm})_{\mathrm{vert}}$
(cf. \eqref{dirs9}),
with a slight abuse of notation, we will likewise denote 
the $\mathcal{S}$-bilinear
extensions to the corresponding $\Sigm$-valued 
forms by
\[
\pi'\colon  
(\mathrm D_{\Sigm})_{\mathrm{vert}}
\otimes_{\Sigm}
(\mathrm D_{\Sigm})_{\mathrm{vert}}
\longrightarrow
\Sigm,
\ 
\pi''\colon  
(\mathrm D_{\Sigm})_{\mathrm{hor}}
\otimes_{\Sigm}
(\mathrm D_{\Sigm})_{\mathrm{hor}}
\longrightarrow
\Sigm
\]
as well.
Furthermore, as before, let $\Omega\colon  L''\otimes_AL'' 
\to L'$ denote the $\mathbf e$-curvature of the  
 $\mathbf e$-connection $\omega$ (see  \eqref{2.3}). 
For $\alpha'',\beta'' \in L''$, let
$d\alpha'',d\beta'' \in \mathrm D_{\Sigm''}$ be
their formal differentials over the ground ring $R$, and denote by
\[
[d\alpha''], [d\beta''] \in \mathrm D_{\Sigm''|A} \cong \Sigm'' \otimes_AL''
\]
the associated formal differentials over $A$, that is,
the images of $d\alpha''$ and $d\beta''$ under the projection
$\mathrm D_{\Sigm''}\to  \mathrm D_{\Sigm''|A}$
 in
\eqref{dif2}, with the notation adjusted suitably.
We will write
\begin{align}
\label{sharp_omega} 
\Omega^{\sharp}([d\alpha''],[d\beta''])&=\Omega(\alpha'',\beta''),
\end{align}
and we extend $\Omega^{\sharp}$ to an $\Sigm$-bilinear alternating
2-form
\[
\Omega^{\sharp} \colon  
\left(\Sigm\otimes_{\Sigm''}\mathrm D_{\Sigm''|A}\right) 
\otimes_{\Sigm}
\left(\Sigm\otimes_{\Sigm''} \mathrm D_{\Sigm''|A}\right) \longrightarrow
\Sigm
\]
in the obvious manner.
We then write the  $\Sigm$-bilinear alternating 2-form on 
$(\mathrm D_{\Sigm})_{\mathrm{hor}}
=\Sigm\otimes_{\Sigm''}\mathrm D_{\Sigm''}$
induced via the projection
\[
(\mathrm D_{\Sigm})_{\mathrm{hor}}
=\Sigm\otimes_{\Sigm''}\mathrm D_{\Sigm''}
\longrightarrow
\Sigm\otimes_{\Sigm''} \mathrm D_{\Sigm''|A}
\cong
\Sigm\otimes_A L''
\]
as
\[
\widetilde\Omega \colon  
(\mathrm D_{\Sigm})_{\mathrm{hor}}
\otimes_{\Sigm}
(\mathrm D_{\Sigm})_{\mathrm{hor}}
\longrightarrow
\Sigm .
\]

\begin{Theorem} 
\label{thm1}
Given $f,h \in \Sigm=\Sigm_A[L]$, the Poisson bracket $\{f,h\}
\in \Sigm$ is given by the expression
\begin{equation}
\begin{aligned}
\{f,h\} &=\pi'((df)_{\mathrm {vert}},(dh)_{\mathrm {vert}})
+
\pi''((df)_{\mathrm{hor} },(dh)_{\mathrm {hor}})
+\widetilde \Omega((df)_{\mathrm{hor} },(dh)_{\mathrm {hor}}) \\
&\quad + \pi ((df)_{\mathrm{vert}},(dh)_{\mathrm{hor}})
+ \pi ((df)_{\mathrm{hor}},(dh)_{\mathrm{vert}}).
\end{aligned}
\label{pois3}
\end{equation}
Furthermore, the 
Poisson 2-form
$\pi$ on $\mathrm D_{\Sigm}$,
 restricted to 
$(\mathrm D_{\Sigm})_{\mathrm{hor}}=S \otimes _{\Sigm ''}\mathrm D_{\Sigm''}$, 
 is given by the sum 
$\pi'' +\widetilde \Omega$ and,
 restricted to 
$(\mathrm D_{\Sigm})_{\mathrm{vert}}$,
by $\pi'$.
\end{Theorem}

\begin{proof} Let $f,h \in \Sigm=\Sigm_A[L]$.
We will show that
\begin{equation}
\label{first_eq}
\pi((df)_{\rm vert}, (dh)_{\rm vert}) = 
\pi'((df)_{\rm vert}, (dh)_{\rm vert})
\end{equation}
and
\begin{equation}
\label{second_eq}
\pi((df)_{\rm hor}, (dh)_{\rm hor}) = 
\pi''((df)_{\rm hor}, (dh)_{\rm hor})
+ \widetilde{\Omega}((df)_{\rm hor}, (dh)_{\rm hor}).
\end{equation}
Since
\begin{align*}
\{f,h\} &= \pi(df, dh) 
= \pi((df)_{\rm hor} + (df)_{\rm vert},
(dh)_{\rm hor} + (dh)_{\rm vert})\\
&= \pi((df)_{\rm vert}, (dh)_{\rm vert}) + 
\pi((df)_{\rm hor}, (dh)_{\rm hor}) + 
\pi((df)_{\rm vert}, (dh)_{\rm hor}) +
\pi( (df)_{\rm hor}, (dh)_{\rm vert}),
\end{align*}
\eqref{first_eq} and \eqref{second_eq}
imply \eqref{pois3}.

As an $R$-algebra, $\Sigm$ is generated by $A$, $L'$, and $L''$.
It suffices to establish the 
identities \eqref{first_eq} and \eqref{second_eq}
(and hence the 
identity \eqref{pois3})
when $f$ and $h$ range over $a \in A$,
$\alpha' \in L'$, and $\alpha'' \in L''$. By construction, 
\begin{equation}
\label{list}
\left\{
\begin{aligned}
da &\in \mathcal S \otimes_A \mathrm D_A \subseteq (\mathrm D_{\Sigm})_{\mathrm{hor}}
\\
d\alpha'' &\in \mathcal{S}\otimes_{ \mathcal{S}''}\mathrm{D}_{\mathcal{S}''} = 
( \mathrm{D}_{ \mathcal{S}})_{\rm hor}\\
d\alpha' &\in \mathcal{S}\otimes_{ \mathcal{S}'} 
\mathrm{D}_{ \mathcal{S}'}  \text{ and hence }
d\alpha' = (d\alpha')_{\rm vert} + ( d\alpha')_{\rm hor}, 
\text{ where}\\
& \qquad (d\alpha')_{\rm vert} \in 
(\Sigm \otimes_{\Sigm'}
\mathrm D_{\Sigm'} )_{\mathrm{vert}} =  
( \mathrm{D}_{ \mathcal{S}})_{\rm vert}\\
& \qquad (d\alpha')_{\rm hor} =
r(d \alpha') \in \mathcal{S}' \otimes_A \mathrm{D}_A
\subset 
\mathcal{S}\otimes_A \mathrm{D}_A \subset 
\mathcal{S}\otimes_{\mathcal{S}''} 
\mathrm{D}_{ \mathcal{S}''} = 
(\mathrm{D}_\mathcal{S})_{\rm hor} .
\end{aligned}
\right.
\end{equation}
To prove \eqref{first_eq} and \eqref{second_eq}, we
need to check them for the pairs
$(a, b)$, $(a, \alpha')$, $(a, \alpha'')$, 
$(\alpha', \beta')$, $( \alpha', \alpha'')$, 
$(\alpha'', \beta'')$, where $a, b \in A$, $\alpha',
\beta' \in L'$, $\alpha'', \beta'' \in L''$.

In view of \eqref{list}, the only non-zero pair on which
either side of \eqref{first_eq} does not vanish is 
$(\alpha', \beta')$. 
In view of
\eqref{3.16.2.1},
 by the definition of the tautological Poisson structure on $L'$,
the right-hand side of \eqref{first_eq} equals 
$[ \alpha', \beta']_{L'}$ whereas, by \eqref{3.16.2.11},
 the left-hand side has the same 
value.

Next, we prove \eqref{second_eq} by verifying it on
each pair listed above. On $(a, b)$ both sides vanish. Indeed, 
in view of
the definition of the tautological Poisson bracket,
$\pi(da, db) = \pi'(da,db)$ is zero. In
addition, $\widetilde{\Omega}(da, db)=0$
since, by construction, the 2-form
$\widetilde \Omega$ factors through the 2-form 
$\Omega^{\sharp}$, and $\Omega^{\sharp}$ vanishes on generators of the
form $da$. 

On the pair $(a, \alpha')$ we proceed in the following
way. Since $(d \alpha')_{\rm hor} \in 
\mathcal{S}\otimes_A \mathrm{D}_A$ we need to check
equality for a pair $(sda, db)$, where $s \in \mathcal{S}$, and $a,b \in A$. 
By construction, $\pi(sda, db) =s\{a,b\}=0$ 
and, similarly, $\pi''(sda, db) =s\{a,b\}=0$. The term 
$\widetilde{\Omega}(s da,db) = 
s\Omega^\sharp(da,db)$ is zero as before. 

On the pair $(a, \alpha'')$, by \eqref{3.16.2.22}, the left-hand side is
$\pi(da, d \alpha'') = \{a, \alpha''\} = \alpha''(a)$. 
By \eqref{3.16.2.2}, the right-hand side equals
$\pi''(da, d \alpha'') + \widetilde{\Omega}(da, d\alpha'')
= \pi''(da, d\alpha'') = \alpha''(a)$.

To compute both sides of \eqref{second_eq} on the pair 
$( \alpha', \beta')$, in view of \eqref{list}, we just need to
evaluate them on elements of the form $\alpha'=s_1da_1$ and 
$\beta'=s_2da_2$ for $s_1, s _2 \in \mathcal{S}$ and $a _1, a _2\in A$. 
Both sides of \eqref{second_eq} vanish.

To compute both sides of \eqref{second_eq} on the pair 
$( \alpha', \alpha'')$, in view of \eqref{list}, we just need to
evaluate them on elements of the form $sda$ and
$d \alpha''$ for $s \in \mathcal{S}$, $a \in\mathcal{S}$,
and $\alpha'' \in L''$. The left-hand side equals
$\pi(sda, d \alpha'') = s \alpha''(a)$ by 
\eqref{3.16.2.22}.  By 
\eqref{3.16.2.2}, the right-hand side equals
\[
\pi''(sda, d \alpha'') + 
\widetilde{\Omega}(sda, d \alpha'') = 
s\pi''(da, d \alpha'') = s \alpha''(a).
\]

On the pair $(\alpha'', \beta'')$, by 
\eqref{3.16.2.12}, the left-hand side
equals 
\[
\pi(d \alpha'', d \beta'') = 
\{\alpha'', \beta''\} = [ \alpha'', \beta'']_{L''}
+ \Omega(\alpha'', \beta'').
\]
By the definition of
the tautological Poisson bracket on $L''$ (see 
\eqref{3.16.2.1} and \eqref{sharp_omega}),
the right-hand side equals 
\[
\pi''(d \alpha'', d \beta'')
+ \widetilde{\Omega}(d \alpha'', d \beta'') = 
\pi''(d \alpha'', d \beta'')+\Omega^\sharp([d \alpha''],
[d \beta'']) = [ \alpha'', \beta'']_{L''} + 
\Omega(\alpha'', \beta'').
\]

Since the differentials of two members $\alpha''$ and $\beta''$ of $L''$,
viewed as elements of $\Sigm$, have zero vertical components,
the 
Poisson 2-form
$\pi$ on $\mathrm D_{\Sigm}$,
 restricted to 
$(\mathrm D_{\Sigm})_{\mathrm{hor}}=S \otimes _{\Sigm ''}\mathrm D_{\Sigm''}$,
 is obviously given by the sum 
$\pi'' +\widetilde \Omega$.
The remaining \lq furthermore\rq\ statement is a consequence of 
Lemma {\rm \ref{apbp}} below.
\end{proof}

\begin{Lemma}\label{apbp}
Given $\alpha',\beta' \in L'$, the value
$\pi((d \alpha')_{\rm vert}, 
(d \beta')_{\rm hor})$ is zero. 
\end{Lemma}

\begin{proof}
Recall from 
\eqref{list} that $(d \beta')_{\rm hor} \in \mathcal{S}
\otimes_A \mathrm{D}_A$ and 
\[
(d \alpha')_{\rm vert} 
\in (\mathrm{D}_{\mathcal{S}})_{\rm vert} = 
(\mathcal{S}\otimes_{ \mathcal{S}'} 
\mathrm{D}_{ \mathcal{S}'})_{\rm vert} \subset 
\mathcal{S}\otimes_{ \mathcal{S}'} 
\mathrm{D}_{ \mathcal{S}'}.
\] 
Therefore, we can
assume that $(d \alpha')_{\rm vert} = s d\beta''$  or 
$(d \alpha')_{\rm vert} = sda$, for
some $a \in A$, $\beta'' \in L'$, and $s \in \mathcal{S}$.
Similarly, we can assume that $(d \beta')_{\rm hor} = 
\bar{s}db$ for some $b \in A$, $\bar{s} \in \mathcal{S}$. 
Therefore $\pi(s d\beta'',\bar{s}db) = 
s\bar{s}\{\beta'',b\}$ vanishes, by \eqref{3.16.2.222}. Likewise,
$\pi(s da,\bar{s}db) = s\bar{s}\{a,b\} = 0$.
\end{proof}

\subsection{The special case $L''=\mathrm{Der}(A)$}
\label{special}

We will show that,
in this case,
the $\mathbf e$-connection $\omega$
determines a unique linear connection 
$r_{\omega}$
for $L'$ and that, furthermore, with
this linear connection, 
the terms 
$\pi((df)_{\mathrm{vert}},(dh)_{\mathrm{hor}})$
and $ \pi ((df)_{\mathrm{hor}},(dh)_{\mathrm{vert}})$
in formula {\rm \eqref{pois3}} vanish. 
This 
amounts to showing that,
given $\alpha'\in L'$ and $\alpha'' \in L''$,
the Poisson bracket
\[
\{\alpha'',\alpha'\} = [\alpha'',\alpha'],
\]
cf.  \eqref{3.16.2.112}, is completely absorbed in the term
$\pi''(d\alpha'',(d\alpha')_{\mathrm {hor}})$, that is,
the component
$\pi(d\alpha'',(d\alpha')_{\mathrm {vert}})$
is zero.
In the case at hand, this observation greatly
simplifies
the formula \eqref{pois3}.

For illustration we note that,
in the circumstances of Example \ref{extr},
given
a function $f$ on $B$, a vector $X \in \mathfrak g$, 
and  $Y \in \mathfrak{X}(B)$, 
when $f X$ and $Y$ are viewed as functions on $N$ in the obvious manner
as in 
Example \ref{extr},
\[
\{Y,fX\}= Y(f) X =\pi''(dY,(d(f X))_{\mathrm{hor}}).
\]
Thus  the Poisson bracket $\{Y,fX\}$ is manifestly absorbed
in the term $\pi''(dY,(d(f X))_{\mathrm{hor}}$.
The reasoning below 
involving the linear connection
induced from the chosen $\mathbf e$-connection
somewhat reduces the argument to 
what corresponds, at the infinitesimal level,
algebraically to
that special case.

We return to the situation of Section \ref{sec_extensions} above. 
The extension of $(R,A)$-Lie algebras under discussion now has the form
\begin{equation}
\label{exact_sequence_der}
\mathbf e\colon 
0
\longrightarrow
L'
\longrightarrow
L
\longrightarrow
\mathrm{Der}(A)
\longrightarrow
0
\end{equation}
but we continue to write $L''=\mathrm{Der}(A)$.
The sequence
of inclusions $R \subseteq A \subseteq \Sigm'$ 
determines the exact sequence
\begin{equation}
0
\longrightarrow
\Sigm' \otimes_A \mathrm D_A
\longrightarrow
\mathrm D_{\Sigm'}
\longrightarrow
\Sigm' \otimes_A L'
\longrightarrow
0;
\label{middle2}
\end{equation}
this is a particular case of the exact sequence
\eqref{dif2}. 
Taking $\Sigm'$-duals,
we obtain the exact sequence
\begin{equation}
\begin{CD}
0
@>>>
\mathrm{Der}(\Sigm'\big|A)
@>>> 
\mathrm{Der}(\Sigm')
@>>> 
\Sigm'\otimes_A L''
@>>> 
0 
\end{CD}
\label{middle9}
\end{equation}
with $L''=\mathrm{Der}(A)\cong \mathrm{Hom}_A(\mathrm D_A,A)$;
the sequence \eqref{middle9} is a particular case
of \eqref{dddif}.
Notice that, 
when we induce up the sequence \eqref{middle2}
to one of $\Sigm$-modules (by
taking the tensor product with $\Sigm$ 
over $\Sigm'$), we obtain the upper horizontal exact sequence 
of $\Sigm$-modules in the commutative diagram 
\eqref{CD1}.
Taking $\Sigm'$-duals once more,
we arrive at a commutative diagram
\begin{equation}
\begin{CD}
0
@>>>
\Sigm' \otimes_A \mathrm D_A
@>>>
\mathrm D_{\Sigm'}
@>>>
\Sigm' \otimes_A L'
@>>>
0
\\
@.
@VVV
@VVV
@VV{\mathrm{Id}}V
@.
\\
0
@>>>
\Sigm' \otimes_A \mathrm D^{**}_A
@>>>
\mathrm D^{**}_{\Sigm'}
@>>>
\Sigm' \otimes_A L'
@>>>
0,
\end{CD}
\label{CD11}
\end{equation}
the unlabeled vertical maps being the canonical maps to the
double duals, the double $A$-dual of $D_A$ being written as
$\mathrm D^{**}_A$ and, likewise,
the double $\Sigm'$-dual of $\mathrm D_{\Sigm'}$ being written as
$\mathrm D^{**}_{\Sigm'}$.

\begin{Proposition}
\label{r_omega_prop}
An $\mathbf{e}$-connection, i.e., an $A$-module section 
$\omega\colon  \operatorname{Der}(A) \rightarrow L$ for the extension 
\eqref{exact_sequence_der} of 
$(R,A)$-Lie algebras, gives rise to a linear connection 
$r_ \omega\colon  \mathrm{D}_{ \mathcal{S}'} \rightarrow 
\mathcal{S}' \otimes_A \mathrm{D}_A$ on the $A$-module $L'$ 
(equivalently: $\mathcal{S}'$-module retraction for \eqref{middle2})
in such a way that 
the $\Sigm'$-module section
\begin{equation}
s_{\omega}\colon \Sigm'\otimes_A L'' \longrightarrow \mathrm{Der}(\Sigm')
\label{cono2}
\end{equation}
for {\em \eqref{middle9}}
dual to  $r_{\omega}$
is uniquely determined by $\omega$.
\end{Proposition}

\begin{proof}
Let $\omega \colon L'' \to L$ be an $\mathbf e$-connection, i.e.,
$A$-module section for the extension
$\mathbf e$
of $(R,A)$-Lie algebras.
Given $\alpha'' \in L''$,
consider the association
\[
L' \longrightarrow L',\ 
(\alpha'',\alpha') \longmapsto [\omega(\alpha''),\alpha'] \in L'
\]
where the notation $[\omega(\alpha''),\alpha']$
refers to the commutator in $L$.
This defines an $R$-linear map
\begin{equation}
L'' \longrightarrow \mathrm{End}(L').
\end{equation}
Since $L'$ acts trivially on $A$, 
given $\alpha'' \in L''$, the formulas
\[
d_{\alpha''}(\alpha') = [\omega(\alpha''),\alpha'],
\
d_{\alpha''}(a) = \alpha''(a),\ 
\alpha' \in L',\ a \in A,
\]
yield an $R$-derivation $d_{\alpha''}\colon \Sigm' \to \Sigm'$, 
the assignment to $\alpha'' \in L''$ of 
$d_{\alpha''}$
yields
an $A$-linear map
from $L''$ to $\mathrm{Der}(\Sigm')$, 
and the extension
\begin{equation}
s_{\omega}\colon \Sigm'\otimes_A L'' \longrightarrow \mathrm{Der}(\Sigm')
\label{cono}
\end{equation}
thereof to an $\Sigm'$-linear map
is an $\Sigm'$-module section for \eqref{middle9},
plainly uniquely determined by $\omega$.
With the notation
\begin{equation}
\label{der_s_prime_split}
(\mathrm{Der}(\Sigm'))_{\mathrm{vert}}=\mathrm{Der}(\Sigm'\big|A),
\
(\mathrm{Der}(\Sigm'))_{\mathrm{hor}}
=s_{\omega}(
\Sigm'\otimes_A L''),
\end{equation}
we obtain an  $\Sigm'$-module direct sum decomposition
\[
\mathrm{Der}(\Sigm')
=(\mathrm{Der}(\Sigm'))_{\mathrm{vert}}
\oplus
(\mathrm{Der}(\Sigm'))_{\mathrm{hor}}
\]
of 
$\mathrm{Der}(\Sigm')$ into a vertical component
$(\mathrm{Der}(\Sigm'))_{\mathrm{vert}}$
and a horizontal component
$(\mathrm{Der}(\Sigm'))_{\mathrm{hor}}$.
The dual of \eqref{cono} is  a retraction
\[
\widetilde r_{\omega} \colon 
\mathrm D^{**}_{\Sigm'}
\longrightarrow
\Sigm'
\otimes_A\mathrm D^{**}_A
\]
for the bottom row extension
in diagram \eqref{CD11}.
This retraction $\widetilde r_{\omega}$ induces the direct sum decomposition
$\mathrm D^{**}_{\Sigm'}=(\Sigm'
\otimes_A\mathrm D^{**}_A)\oplus \mathrm{ker}(\widetilde r_{\omega})$
in the category of $\Sigm'$-modules, and the restriction
to $ \mathrm{ker}(\widetilde r_{\omega})$ 
of the projection to $\Sigm'\otimes_A L'$ is an isomorphism of
$\Sigm'$-modules. Lifting the resulting section
$\Sigm'\otimes_A L' \to \mathrm D^{**}_{\Sigm'}$
to a section $\Sigm'\otimes_A L' \to \mathrm D_{\Sigm'}$
for
\eqref{middle2}, we obtain a direct sum decomposition
\[
\mathrm D_{\Sigm'}\cong(\Sigm'
\otimes_A\mathrm D_A)\oplus \mathrm{ker}(\widetilde r_{\omega})
\cong
(\Sigm'
\otimes_A\mathrm D_A)\oplus\Sigm'\otimes_A L',
\]
and the associated retraction to 
$\Sigm'
\otimes_A\mathrm D_A$ is a retraction
\[
r_{\omega} \colon 
\mathrm D_{\Sigm'}
\longrightarrow
\Sigm'
\otimes_A\mathrm D_A
\]
for the extension \eqref{middle2}
in the category of $\Sigm'$-modules
of the kind we are looking for.
\end{proof}

We will refer to a  connection of the kind $r_{\omega}$ 
in Proposition \ref{r_omega_prop}
as a
{\em linear connection\/} on $L'$ {\em associated to
the\/} $\bf e$-connection $\omega$.

\begin{Theorem} 
\label{thm2}
Suppose that the quotient $L''$
in the extension $\mathbf e$ of $(R,A)$-Lie algebras
under discussion
is the $(R,A)$-Lie algebra
$\mathrm{Der}(A)$ and that
the requisite linear connection on $L'$ is a  connection 
$r_{\omega}$
associated to
the $\mathbf e$-connection 
$\omega\colon  L'' \to L$.
Then the terms 
$\pi ((df)_{\mathrm{vert}},(dh)_{\mathrm{hor}})$
and $ \pi ((df)_{\mathrm{hor}},(dh)_{\mathrm{vert}})$
in the formula {\rm \eqref{pois3}} vanish.
Hence, given  $f, h\in \Sigm=\Sigm_A[L]$,
\begin{equation}
\begin{aligned}
\{f,h\} &=\pi'((df)_{\mathrm {vert}},(dh)_{\mathrm {vert}})
+
\pi''((df)_{\mathrm{hor} },(dh)_{\mathrm {hor}})
+\widetilde \Omega((df)_{\mathrm{hor} },(dh)_{\mathrm {hor}}) .
\end{aligned}
\label{pois2}
\end{equation}
\end{Theorem}

\begin{proof}
We return to
the proof of Theorem \ref{thm1}.
When the members $f$ and $h$ of $\Sigm$
range over
the members of $A$, $L'$, and $L''$, the only cases where the value
$\pi((df)_{\mathrm{vert}},(dh)_{\mathrm{hor}})$ does not
obviously vanish are those where $(f,h)=(\alpha', \beta')$ 
and $(f,h)=(\alpha', \alpha'')$, with 
$\alpha', \beta' \in L'$ and $\alpha'' \in L''$.

In view of Lemma \ref{apbp}, $\pi((d \alpha')_{\rm vert}, 
(d \beta')_{\rm hor}) = 0$. 
Thus, it remains to consider the case $(f,h)=(\alpha', \alpha'')$, i.e., since $d \alpha'' \in 
(\mathrm{D}_\mathcal{S})_{\rm hor}$, we need to show 
that $\pi((d \alpha')_{\rm vert}, d \alpha'') =0$.
\smallskip

As before, we will now identify $L''$ with 
$\omega(L'')$ so that,
as $A$-modules, $L=L'\oplus L''$, and 
$\omega$ is the canonical injection. In view
of Proposition \ref{r_omega_prop}, the
decomposition \eqref{dirs12} now takes the form
\begin{equation}
\mathrm D_{\Sigm'}  
=(\mathrm D_{\Sigm'})_{\mathrm {hor}} 
\oplus (\mathrm D_{\Sigm'})_{\mathrm {vert}}=
\left(\Sigm'
\otimes_{A}\mathrm D_{A}\right) \oplus
\mathrm{ker}(r_{\omega}) .
\end{equation}
Accordingly, with the notation $\widetilde r_{\omega}$
used in the proof of Proposition
\ref{r_omega_prop},
\begin{equation}
\mathrm D^{**}_{\Sigm'}  
=
\left(\Sigm'
\otimes_{A}\mathrm D^{**}_A\right) \oplus
\mathrm{ker}(\widetilde r_{\omega}) .
\end{equation}
With a slight abuse of notation,
given a formal differential $\beta \in \mathrm D_{\Sigm'}$,
we will denote its image in $\mathrm D^{**}_{\Sigm'}$,
i.e., the induced $\Sigm'$-linear
map from  $\mathrm{Der}(\Sigm')$ to $\Sigm'$, by
$\beta\colon \mathrm{Der}(\Sigm') \to \Sigm'$
 as well.   

By the definition of $\widetilde r_ \omega$, the kernel
$\mathrm{ker}(\widetilde r_{\omega})$ 
of $\widetilde r_{\omega}$ 
consists of the 
$\Sigm'$-linear maps
$\beta\colon  \mathrm{Der}(\Sigm') \to \Sigm'$ such that the composition
\begin{equation}
\begin{CD}
\Sigm' \otimes_A L'' 
@>{s_{\omega}}>>
\mathrm{Der}(\Sigm') 
@>{\beta}>>
\Sigm'
\end{CD}
\end{equation}
is zero. By construction, then, 
$d_{\alpha''}=s_{\omega}(\alpha'') \in (\mathrm{Der}(\Sigm'))_{\rm hor}$, 
in view of the decomposition \eqref{der_s_prime_split},
and $(d\alpha')_{\mathrm{vert}}\colon 
\mathrm{Der}(\Sigm') 
\to
\Sigm'
$
has the following property:
\[
(d\alpha')_{\mathrm{vert}}(d_{\alpha''})=
\left((d\alpha')_{\mathrm{vert}} \circ s_\omega\right)(\alpha'')=0.
\]

Given $\beta \in \mathrm D_{\Sigm'}$ and
$\alpha'' \in L''$, so that $d_{\alpha''}\in 
\mathrm{Der}(\Sigm')$, we will show that 
\begin{equation}
\beta(d_{\alpha''}) =\pi(d\alpha'',\beta) .
\label{iden}
\end{equation}
Indeed,
as an $R$-module, $\mathrm D_{\Sigm'}$
is generated by the elements $s d\alpha'$ 
and $tda$,
where $s,t\in \Sigm'$,
$\alpha'\in L'$, and $a\in A$.
Now
\begin{align*}
(s d\alpha')(d_{\alpha''})&=s (d\alpha'(d_{\alpha''}))
= s (d_{\alpha''}( \alpha'))
=s [\alpha'',\alpha'] = s\{\alpha'',\alpha'\}
=s\pi(d\alpha'',d\alpha')
=\pi(d\alpha'',s d\alpha')
\\
(t da)(d_{\alpha''})&=t (da(d_{\alpha''}))
=t(d_{\alpha''}(a))
=t \alpha''(a) = t\{\alpha'',a\}
=t\pi(d\alpha'',da)
=\pi(d\alpha'',tda)
\end{align*}
which proves \eqref{iden}.

By construction, in view of \eqref{iden},
\[
\pi(d\alpha'',(d\alpha')_{\mathrm{hor}})=(d\alpha')_{\mathrm{hor}}
(d_{\alpha''})
=((d\alpha')_{\mathrm{vert}}
+(d\alpha')_{\mathrm{hor}})
(d_{\alpha''})
=
(d\alpha')(d_{\alpha''})
=d_{\alpha''}(\alpha')=
[\alpha'',\alpha'].
\]
On the other hand, by \eqref{3.16.2.112},
\[
[\alpha'', \alpha']=
\{\alpha'',\alpha'\} =\pi (d\alpha'',d \alpha')
=\pi (d\alpha'',(d\alpha')_{\mathrm{hor}})
+\pi (d\alpha'',(d\alpha')_{\mathrm{vert}}).
\]
This shows that $\pi (d\alpha'',
(d\alpha')_{\mathrm{vert}}) =0$.
\end{proof}

\begin{Remark}
\label{rem5}
{\rm 
In Theorem {\rm \ref{thm2}},
relative to the direct sum decomposition
$\mathrm D_{\Sigm} =(\mathrm D_{\Sigm})_{\mathrm {hor}} 
\oplus (\mathrm D_{\Sigm})_{\mathrm {vert}}$,
cf. {\rm \eqref{dirs9}},
the $\Sigm$-valued $\Sigm$-bilinear alternating form $\pi$ on
$\mathrm D_{\Sigm}$---the Poisson form---has the special feature 
that it decomposes into the two 
separate components
$\pi'' + \widetilde \Omega$ on 
$(\mathrm D_{\Sigm})_{\mathrm {hor}}$
and
$\pi'$ on 
$(\mathrm D_{\Sigm})_{\mathrm {vert}}$, 
that is, though defined on all of
$\mathrm D_{\Sigm}$, the 2-form $\pi$ has
no non-zero cross term on
$(\mathrm D_{\Sigm})_{\mathrm {vert}} 
\times (\mathrm D_{\Sigm})_{\mathrm {hor}}$.
}\quad $\blacklozenge$
\end{Remark}

\begin{Remark}
\label{rem6}
{\rm 
In the situation of Theorem {\rm \ref{thm1}},
suppose that the Lie-Rinehart structure map
from $L''$ to  $\mathrm{Der}(A)$ 
is injective,
and let $\mathcal R = A^{L''}\colon =\{a \in A\mid \alpha''(a) = 0 , 
\alpha'' \in L''\}$. By definition, 
$L''=\mathrm{Der}(A|\mathcal R)$.
For example, when
$A$ is the ring of smooth functions on a smooth manifold $B$,
in view of the Frobenius theorem,
an $(\mathbb R,A)$-Lie algebra $L''$ of the kind under discussion
defines a foliation of $B$,
and $\mathcal R$  is then the algebra of functions
that are constant on the leaves.
In the general case,
we can then build the theory
with $\mathcal R$ as ground ring rather than just $R$:
Relative to the sequence of inclusions 
$\mathcal R \subseteq A \subseteq \Sigm$, 
the exact sequence \eqref{dif2} has the form
\[
0
\longrightarrow
\Sigm\otimes_A \mathrm D_{A|\mathcal R}
\longrightarrow
\mathrm D_{\Sigm|\mathcal R}
\longrightarrow
\Sigm \otimes _AL
\longrightarrow
0,
\]
where we keep in mind that the canonical 
$\Sigm$-module
morphism
$ \Sigm \otimes _AL \to \mathrm D_{\Sigm|A}$
is an isomorphism.
That sequence  fits into the commutative diagram
\begin{equation}
\begin{CD}
0
@>>>
\mathcal S \otimes _A \mathrm D_A @>>> \mathrm
D_{\mathcal S} @>>> \mathcal S \otimes _A L
@>>> 0 
\\
@.
@VVV
@VVV
@|
@.
\\
0
@>>>
\Sigm\otimes_A \mathrm D_{A|\mathcal R}
@>>>
\mathrm D_{\Sigm|\mathcal R}
@>>>
\Sigm \otimes _AL 
@>>>
0
\end{CD}
\end{equation}
having the two vertical arrows surjective;
the top row is precisely the exact sequence \eqref{dif2}
relative to the ground ring $R$.
With $\mathcal R$ as ground ring, 
Theorem {\rm \ref{thm2}}
applies, that is, the Poisson bracket is given by the formula
{\rm\eqref{pois2}}, and there are no non-zero cross terms,
even when
the injection of $L''$ into $\mathrm{Der}(A)$
is not an isomorphism.
This observation justifies building the theory over a
ground ring more general than a field.} \quad $\blacklozenge$
\end{Remark}

\begin{Example}{\rm 
Consider the situation in Example \ref{extr}, that is,
$P=B \times G$ and $P \rightarrow B$ is the trivial principal $G$-bundle
endowed with the trivial principal connection.
In this case $R=\mathbb R$, $A = C ^{\infty}(B)$ and $L'' = 
\mathfrak{X}(B) = \operatorname{Der}(C ^{\infty}(B))$, so Theorem \ref{thm2} applies. The algebra $\Sigm=\Sigm_A[L]$ comes down to the tensor product algebra
$\Sigm[\mathfrak g]\otimes \Sigm_A[L'']$.
Let $a\in C ^{\infty}(B)$, $\lambda \in \mathfrak{g}$, 
and  $Y \in \mathfrak{X}(B)$, and view
$a,\lambda, Y$ as functions on $N=
\mathrm{T}^\ast B \times \mathfrak{g}^\ast$. 
The algebra $\Sigm'=\Sigm_A[L']$ decomposes
as the tensor product algebra $A \otimes \Sigm[\mathfrak g]$
and, as an $\Sigm'$-module,
$\mathrm D_{\Sigm'}$ decomposes canonically as the direct sum
$\Sigm[\mathfrak g]\otimes \mathrm D_A \oplus
\Sigm'\otimes \mathfrak g$.
The projection to $\Sigm[\mathfrak g]\otimes \mathrm D_A$
yields the requisite linear connection, and
it is immediate that $({d}\lambda)_{\rm hor} = 0$.
Thus, by \eqref{pois2},
the Poisson bracket on $\mathrm{T}^\ast B \times \mathfrak{g}^\ast $ is the sum of the canonical
(minus) Poisson bracket on $\mathrm{T}^\ast B$ and
 the minus
Lie-Poisson bracket on $\mathfrak{g}^\ast$,
cf.  Remark {\rm \ref{rem00}} for the Lie-Poisson 
bracket on $\mathfrak{g}^\ast$ that comes into play here.

In the next subsection we will determine the bracket for
a general principal bundle.
\quad $\blacklozenge$}
\end{Example}

\subsection{The gauged Lie-Poisson bracket}
\label{principal}

Return to the situation of
Example {\rm {\ref{2.222}}}:
Thus $\prin\colon Q \to \BB$ is a right principal
$G$-bundle.
The ground ring $R$ is now that of the real numbers, $\mathbb R$,
and the algebra $A$ is the algebra $C^{\infty}(\BB)$.

Let $L'$ denote the $A$-Lie algebra of sections 
of the adjoint bundle 
$
\operatorname{ad}({\prin})\colon
Q \times_G \mathfrak g
\rightarrow \BB$
of $\prin$; as before,
we will view $L'$ as an $(\mathbb R,A)$-Lie
algebra with trivial $L'$-action on $A$.
>From the standard (plus) Lie-Poisson structure of
$\mathfrak g^*$,
the total space
of the coadjoint bundle $
\operatorname{ad}^*({\prin})\colon
Q \times_G \mathfrak g^*
\rightarrow \BB$
of $\prin$
acquires a smooth Poisson structure, 
the symmetric $A$-algebra $\Sigm'=\Sigm_A[L']$ on $L'$
carries the tautological Poisson structure introduced in Section 
\ref{sec_2} above, and the Poisson algebra
$\Sigm'$ embeds into $C^{\infty}(Q \times_G \mathfrak g^*)$
as a Poisson algebra in an obvious manner via the dual pairing between
$\mathfrak g$ and $\mathfrak g^*$.
Indeed, $\Sigm'=\Sigm_A[L']$ is precisely the $A$-algebra of 
functions on $Q \times_G \mathfrak g^*$ that are polynomial on the fibers of
the coadjoint bundle projection  $
\operatorname{ad}^*({\prin})$ and
 $\Sigm'$ thus embeds into
$C^{\infty}(Q \times_G \mathfrak g^*)$
as a Fr\'echet dense subalgebra.

As in Subsection \ref{illust}, let
$\widetilde \prin\colon  \widetilde{Q} =
Q \times_{B} \mathrm{T}^\ast B\rightarrow \mathrm{T}^*B$
denote the pull back bundle
relative to the cotangent bundle projection
$\cotan_B\colon \mathrm T^*B \rightarrow B$
and $N_{\mathrm{St}}=\widetilde Q \times_G \mathfrak g^*$ 
the total space of the coadjoint bundle 
$\operatorname{ad}^\ast(\widetilde{\prin})\colon
\widetilde Q \times_G \mathfrak g^*
\rightarrow \mathrm{T}^\ast B$, the {\em Sternberg\/} space 
(associated to the data).
Exactly as before, let
 $\widetilde L'$ denote the $C^{\infty}(\mathrm T^*B)$-Lie algebra of sections 
of the adjoint bundle 
of $\widetilde \prin$;
the Sternberg space
$N_{\mathrm{St}}$, being the
total space
of the coadjoint bundle 
of $\widetilde \prin$,
acquires a smooth Poisson structure, 
the symmetric $C^{\infty}(\mathrm T^*B)$-algebra $\widetilde\Sigm'$
on $\widetilde L'$
carries the tautological Poisson structure, and the Poisson algebra
$\widetilde\Sigm'$ embeds into $N_{\mathrm{St}}$
as a Fr\'echet dense Poisson subalgebra.
We will refer to the present Poisson structure on
$N_{\mathrm{St}}$ as the {\em Lie-Poisson\/} structure
on $N_{\mathrm{St}}$.

Pick a connection on $\prin\colon Q \to \BB$.
This choice of connection induces the corresponding
diffeomorphism 
$N_W \to N_{\mathrm{St}}$ given above as \eqref{idstw}.
From the obvious Poisson structure on $N_W=(\mathrm T^*Q)\big/G$,
the Sternberg space $N_{\mathrm{St}}$ thus inherits a Poisson structure and,
cf. Remark \ref{stb} below, we shall see 
that the new Poisson structure perturbs the Lie-Poisson structure
 on 
$N_{\mathrm{St}}$ in a very precise way.
We are about to show that the new Poisson structure
depends on the Lie-Poisson structure, the cotangent bundle
Poisson structure on $\mathrm T^*B$, and the choice of connection;
we therefore refer to the resulting Poisson
bracket on 
$N_{\mathrm{St}}$
as a {\em gauged\/} Lie-Poisson bracket.

Introduce
the following notation:
\begin{itemize}
\item $\widehat z$ : for a point of $N_{\mathrm{St}}=\widetilde Q\times_G \mathfrak g^*$,
\item $\vartheta_B\colon \mathrm T \mathrm T^*\BB \to \mathbb R$ for the 
tautological 1-form on $\mathrm T^*B$,
\item $\pi_{\mathfrak g^*}\colon N_{\mathrm{St}}=\widetilde Q\times_G \mathfrak g^* \to
\mathrm T^*B$ for the bundle projection,
\item $\beta =\pi_{\mathfrak g^*}(\widehat z)\in \mathrm T^*B$,
\item
$\omega\colon\mathrm T Q \to \mathfrak g$ for the connection form on
$\mathrm T Q$,
\item
$\Omega$ for the curvature thereof,
\item
$\widetilde \omega\colon \mathrm T \widetilde Q \to \mathfrak g$ for
the lifted connection form on $\mathrm T \widetilde Q$,
\item
$\widetilde \Omega$ for the lifted curvature,
\item $\widehat \mu$ for the variable in the fiber
$\mathfrak g^*_{\widehat z}= \pi_{\mathfrak
g^*}^{-1}\left(\pi_{\mathfrak g^*}(\widehat z)\right) \subseteq
N_{\mathrm{St}}$ 
through the point $\widehat z$
of $N_{\mathrm{St}}$, 
\item
$\langle \,\cdot \, , \,
\cdot \, \rangle$ for the canonical pairing between a vector
space and its dual, 
\item
$\sharp\colon \mathrm T^*(\mathrm T^*B) \longrightarrow\mathrm T(\mathrm T^*B)$
for the isomorphism of vector bundles on $\mathrm T^*B$
induced by the  symplectic 
form $\mathbf d\vartheta_B$, the isomorphism being spelled out here on the total spaces. 
\end{itemize}
The principal connection $\omega$ on $\prin$
induces a linear connection on the coadjoint bundle
$\mathrm{ad}^*(\prin)$ 
of $\prin$
and hence on the lifted coadjoint bundle
$\mathrm{ad}^*(\widetilde \xi)\colon 
N_{\mathrm{St}}=\widetilde Q\times_G \mathfrak g^* \to \mathrm T^*\BB$.
Given a function $f$ on $N_{\mathrm{St}}=\widetilde Q\times_G \mathfrak g^*$
and a point $\widehat z$ of  $N_{\mathrm{St}}$,
 the
differential $\mathbf df(\widehat z)$ 
decomposes uniquely in the form
\[
\mathbf df(\widehat z) = \mathbf d_{\widetilde \omega}f(\widehat z) + \mathbf d_{\widehat
\mu}f(\widehat z)
\]
into a vertical summand  $\mathbf d_{\widehat
\mu}f(\widehat z)$  and a horizontal summand 
$\mathbf d_{\widetilde \omega}f(\widehat z)$. 
This decomposition corresponds precisely to the decomposition \eqref{dirs9}
above,
with the same significance of the terms
{\em horizontal\/} and {\em vertical\/}. 

We will now 
substitute the extension 
\eqref{2.2.2} of
$(\mathbb R,A)$-Lie algebras 
for the extension
\eqref{exact_sequence_der}
 of
$(\mathbb R,A)$-Lie algebras and
apply the reasoning in Subsection \ref{special} above.
Thus, as an $(\mathbb R,A)$-Lie algebra, 
$L''=\mathrm{Der}(A)=\mathfrak X(B)$, the algebra
$\Sigm''=\Sigm_A[L'']$ is the algebra of smooth functions
on $\mathrm T^*\BB$ that are polynomial on the fibers of the cotangent bundle
projection $\mathrm T^*\BB \to \BB$,
as an $(\mathbb R,A)$-Lie algebra, $L$ is that of $G$-equivariant
vector fields on $Q$, and the algebra $\Sigm$ is the algebra
of smooth functions on 
$N_{\mathrm{St}}=\widetilde Q\times_G \mathfrak g^*$
that are polynomial on the fibers of the projection
$\widetilde Q\times_G \mathfrak g^*\to \mathrm T^*\BB$.
Formula \eqref{pois2} then yields
formula \eqref{1} 
below; up to sign, this formula 
coincides with the corresponding formula in
\cite{zaalaone}.
We therefore label the subsequent 
result as a proposition rather than as a theorem.

\begin{Proposition}\label{zaal}
The Poisson bracket
$\{f_1,f_2\}$ of two functions $f_1$ and $f_2$ on $\widetilde
Q\times_G \mathfrak g^*$, evaluated at the point $\widehat z$ of
$N_{\mathrm{St}}=\widetilde Q\times_G \mathfrak g^*$, is given by the following expression
where, as before,
 $\beta =\pi_{\mathfrak g^*}(\widehat z)\in \mathrm T^*B$:
\begin{equation}
\begin{aligned}
\{f_1,f_2\}(\widehat z)&= \quad \mathbf d\vartheta_B(\beta)\left(\mathbf d_{\widetilde
\omega}f_1^\sharp(\widehat z), \mathbf d_{\widetilde
\omega}f_2^\sharp(\widehat z)\right) 
\\
&\quad + \left\langle \widehat z, \widetilde
\Omega(\beta)\left(\mathbf d_{\widetilde \omega}f_1^\sharp(\widehat z),
\mathbf d_{\widetilde \omega}f_2^\sharp(\widehat z)\right)\right \rangle
\\
&\quad +\left\langle \widehat z,\left[\mathbf d_{\widehat \mu} f_1(\widehat
z), \mathbf d_{\widehat \mu} f_2(\widehat z)\right ]_{\widehat z} \right
\rangle . 
\end{aligned}
\label{1}
\end{equation}
In this formula, 
$\mathbf d\vartheta_B$ is the resulting symplectic structure on $\mathrm T^*\BB$,
the expression $ \mathbf d_{\widetilde
\omega}f^\sharp(\widehat z)\in \mathrm T_{\beta}(\mathrm T^*B)$
refers to the vector in $\mathrm T_{\beta}(\mathrm T^*B)$ associated
to $\mathbf d_{\widetilde \omega}f(\widehat z)$ under 
$\sharp$ and, in 
$\left\langle \widehat z, \widetilde
\Omega(\beta)\left(\mathbf d_{\widetilde \omega}f_1^\sharp(\widehat z),
\mathbf d_{\widetilde \omega}f_2^\sharp(\widehat z)\right)\right \rangle$
and 
$\left\langle \widehat z,\left[\mathbf d_{\widehat \mu} f_1(\widehat
z), \mathbf d_{\widehat \mu} f_2(\widehat z)\right ]_{\widehat z} \right
\rangle$,
the point $\widehat z$ is viewed as a vector in the vector
space $\mathfrak g^*_{\widehat z}$. 
\end{Proposition}

\begin{Remark}{\rm
In \cite{zaalaone}, the tautological 1-form on $\BB$ is 
written as $\alpha_B$; 
furthermore, in that reference,
the group acts on  the total space $Q$ from the left
(written there as $P$) and,
accordingly, 
the term
$\left\langle \widehat z,\left[\mathbf d_{\widehat \mu} f_1(\widehat
z), \mathbf d_{\widehat \mu} f_2(\widehat z)\right ]_{\widehat z} \right
\rangle$
comes with a minus sign.} \quad $\blacklozenge$
\end{Remark}

\begin{proof}
This is a special case of formula \eqref{pois2}.
The terms
$\pi'((df)_{\mathrm {vert}},(dh)_{\mathrm {vert}})$,
$\pi''((df)_{\mathrm{hor} },(dh)_{\mathrm {hor}})$,
$\widetilde \Omega((df)_{\mathrm{hor} },(dh)_{\mathrm {hor}})$
on the right-hand side of \eqref{pois2} correspond to, respectively,
$\left\langle \widehat z,\left[\mathbf d_{\widehat \mu} f_1(\widehat
z), \mathbf d_{\widehat \mu} f_2(\widehat z)\right ]_{\widehat z} \right
\rangle $,
$\mathbf d\vartheta_B(\beta)\left(\mathbf d_{\widetilde
\omega}f_1^\sharp(\widehat z), \mathbf d_{\widetilde
\omega}f_2^\sharp(\widehat z)\right) $,
and
$\left\langle \widehat z, \widetilde
\Omega(\beta)\left(\mathbf d_{\widetilde \omega}f_1^\sharp(\widehat z),
\mathbf d_{\widetilde \omega}f_2^\sharp(\widehat z)\right)\right \rangle$,
in \eqref{1}.

A straightforward though slightly tedious comparison
shows that the formula \eqref{1} for the value 
of the Poisson bracket
$\{f_1,f_2\}$ at the point $\widehat z$
is  indeed a special case
of the formula \eqref{pois2}. 
\end{proof}

\begin{Remark}
{\rm Our
approach provides a simple explanation for the formula \eqref{1}
for the Poisson structure:
The term $\mathbf d\vartheta_B(\beta)\left(\mathbf d_{\widetilde
\omega}f_1^\sharp(\widehat z), \mathbf d_{\widetilde
\omega}f_2^\sharp(\widehat z)\right)$
comes from the 
symplectic Poisson structure  associated to the 2-form 
$\mathbf d\vartheta_B$
on the base $\mathrm
T^*B$ of the bundle projection
$\pi_{\mathfrak g^*}\colon N_{\mathrm{St}}\to
\mathrm T^*B$
and thus corresponds to the first term on the right-hand side
of \eqref{3.16.2.12} and to the identity \eqref{3.16.2.22} above; the
term $\left\langle \widehat z, \widetilde
\Omega(\beta)\left(\mathbf d_{\widetilde \omega}f_1^\sharp(\widehat z),
\mathbf d_{\widetilde \omega}f_2^\sharp(\widehat z)\right)\right \rangle$
comes from the curvature of the connection and therefore corresponds
to the second term on the right-hand side of \eqref{3.16.2.12}; and
the term $\left\langle \widehat z,\left[\mathbf d_{\widehat \mu}
f_1(\widehat z), \mathbf d_{\widehat \mu} f_2(\widehat z)\right ]_{\widehat
z} \right \rangle$
comes from the Lie-Poisson structure of $\mathfrak g^*$ and
therefore  corresponds to the identity \eqref{3.16.2.11}. 
The two terms in the Poisson tensor beyond that coming from the
cotangent bundle Poisson structure on $\mathrm T^*\BB$ simply
reconstruct the Atiyah sequence \eqref{2.2.1} of the principal
bundle or, equivalently, the corresponding extension of $(\mathbb R,
C^{\infty}(\BB))$-Lie algebras, in terms of the connection. \quad $\blacklozenge$
}\end{Remark} 

\begin{Remark}
\label{stb}
{\rm The summand
$
\left\langle \widehat z,\left[\mathbf d_{\widehat \mu} f_1(\widehat
z), \mathbf d_{\widehat \mu} f_2(\widehat z)\right ]_{\widehat z} \right
\rangle
$
in \eqref{1}
yields the Lie-Poisson bracket on 
the Sternberg space $N_{\mathrm{St}}=\widetilde Q\times_G \mathfrak g^*$.
Thus the formula \eqref{1}
manifestly perturbs the Lie-Poisson structure
of the Sternberg space $N_{\mathrm{St}}$ in a very precise way.
}\quad $\blacklozenge$
\end{Remark}

\subsection{General Lie algebroids}
Let $\lambda\colon  \mathcal{L}\rightarrow B$
be a general Lie algebroid; 
then the pair
$(A,L)=(C^{\infty}(B),\Gamma_B (\lambda))$
acquires a
Lie-Rinehart algebra structure. A particular example is a
Lie algebroid of the kind reproduced in Example \ref{2.222}. For a
general Lie algebroid $\lambda$,  the  algebra $\Sigm_A[L]$ is in an
obvious manner isomorphic to the algebra of smooth functions on the total space
$\mathcal L^*$ of the dual bundle $\lambda^* \colon  \mathcal L^* \to
\BB$ that are polynomial on the fibers of $\lambda^*$, the algebra
$C^{\infty}( \mathcal L^*)$ of ordinary smooth functions on
$\mathcal L^*$ acquires a Poisson structure in an obvious way, and
the injection of $\Sigm_A[L]$ into $C^{\infty}( \mathcal L^*)$ maps
the former algebra onto a Fr\'echet dense subalgebra 
of the latter
and is plainly compatible with the Poisson structures. The special
case of the tangent bundle has been presented in
Subsection \ref{principal}. In the
general case, the image $L''$ of $L$ under the morphism $L \rightarrow \mathfrak{X}(B)$ of $(\mathbb R,A)$-Lie algebras, which is part of
the Lie algebroid structure of $\lambda$, is an $(\mathbb R,A)$-Lie
algebra, and the morphism $L \rightarrow \mathfrak{X}(B)$ of $(\mathbb
R,A)$-Lie algebras induces a morphism $\Sigm_A[L]\to
\Sigm_A[\mathfrak{X}(B)]$ of Poisson algebras, the latter algebra
being the ordinary Poisson algebra of smooth functions on the total
space $\mathrm T^*\BB$ of the cotangent bundle of $\BB$ that are
polynomial on the fibers. Moreover, the projection from $L$ to $L''$
fits into an extension
\begin{equation} 
\label{extension_algebroids}
\mathrm {\mathbf e}
\colon  0 \longrightarrow L' \longrightarrow L \longrightarrow L''
\longrightarrow 0
\end{equation}
of $(\mathbb R,A)$-Lie algebras of the kind \eqref{2.1}.
A special case is the key example \eqref{2.2.2} above.

Under a suitable additional hypothesis, the construction in the
previous section yields the following description of the Poisson
algebra  $\Sigm_A[L]$ in terms of $A$, the extension
$\mathrm{\mathbf e}$,
and an $\mathrm{\mathbf e}$-connection thereof.

\begin{Theorem} \label{alg} Suppose that the anchor
$\rho\colon  \mathcal L \to \mathrm T\BB$ has constant rank, so
that $L''$ defines a foliation $\mathcal F$ on $\BB$, and denote by
$A^{\mathcal F}$ the algebra of smooth functions on $\BB$ that are
constant on the leaves of $\mathcal F$. Then the Poisson algebras
$\Sigm_A[L]$ and $\Sigm_A[L'']$ are Poisson algebras even over
$A^{\mathcal F}$, the functions in $A^{\mathcal F}$ being
Casimir functions, the extension 
$\mathrm{\mathbf e}$
splits in the
category of $A$-modules 
and, once an $\mathrm{\mathbf e}$-connection 
(i.~e., $A$-module splitting of $\mathrm{\mathbf e}$)
has been chosen, the
identities {\rm \eqref{3.16.2.11}--\eqref{3.16.2.33}} yield an
explicit description of the Poisson algebra $\Sigm_A[L]$ of smooth
functions on $\mathcal L^*$. Furthermore, 
in terms of differentials,
the tautological Poisson bracket on $\Sigm_A[L]$ is then as well given by the 
formula {\rm \eqref{pois2}}.
\end{Theorem}

\begin{proof}
Since the anchor $\rho\colon  \mathcal L \to \mathrm T\BB$
has constant rank, $\rho$ maps $\mathcal L$ onto the total space
of a subbundle of the tangent bundle of $\BB$. As an $A$-module,
$L''$ is the space of sections of this vector bundle whence $L''$ is
necessarily projective as an $A$-module. 
Working over 
$A^{\mathcal F}$ as ground ring, cf.
Remark \ref{rem6}, we conclude that
the tautological Poisson bracket on $\Sigm_A[L]$ is still given by the 
formula {\rm \eqref{pois2}}.
This implies the assertion.
\end{proof}

Here is a special case.

\begin{Corollary} \label{cor1}
Suppose that $\lambda\colon  \mathcal L \to \BB$ is a
transitive Lie algebroid. Then $L''$ coincides with the $(\mathbb
R,A)$-Lie algebra $\mathfrak{X}(B)$ of ordinary smooth vector
fields on $\BB$, the extension $\mathrm{\mathbf e}$
splits in the
category of $A$-modules 
and, once an $\mathrm{\mathbf e}$-connection has been chosen,
the
identities {\rm \eqref{3.16.2.11}--\eqref{3.16.2.33}} yield an
explicit description of the Poisson algebra $\Sigm_A[L]$ of smooth
functions on $\mathcal L^*$ and hence a complete description of the
Poisson algebra $\left(C^{\infty}( \mathcal L^*) ,\{\,\cdot\, ,\,
\cdot\, \}\right)$ of ordinary smooth functions on $\mathcal L^*$.
Furthermore, in terms of differentials,
the tautological Poisson bracket on $\Sigm_A[L]$ is then as well given by the 
formula {\rm \eqref{pois2}}
and hence by a suitable variant of {\rm \eqref{1}}.
\end{Corollary}

The Corollary applies, in particular, to the transitive Lie
algebroid \eqref{atiy2}, but that special case has been dealt with
before.

\section{Non-regular quotients of cotangent bundles}
\label{nonreg}
Let $Q$ denote a smooth manifold, $G$ a group acting
smoothly on $Q$ from the right (as before), 
and lift the action to the total space $\mathrm{T}Q$
of the tangent bundle. 
We no longer suppose that the action of $G$ on $Q$ is principal and 
form the {\em non-regular\/} quotient
$(\mathrm{T}^* Q)\big/G$ . Endow the orbit space
$B=Q/G$ with the quotient topology, let
$\pi \colon  Q \to B$ be the canonical projection, and
let $(A_Q,L_Q):=(C^{\infty}(Q),\mathfrak{X}(Q))$. Moreover, let
$A=C^{\infty}(Q)^G$, the algebra of smooth $G$-invariant functions
on $Q$, viewed as an algebra of continuous functions on $B$, and
let $L= \mathfrak{X}(Q)^G$, the Lie algebra of smooth
$G$-invariant vector fields on $Q$. The Lie-Rinehart structure on
$(A_Q,L_Q)$ induces a Lie-Rinehart structure on $(A,L)$, in
particular, an action $L \to \mathrm{Der}(A)$ of $L$ on $A$ by
derivations. In this section we will show that, in the
present more general case, an algebra of the kind $\Sigm_A[L]$
(cf.
Theorem \ref{general} above for the regular case,)
does not suffice
to recover the functions on the quotient $(\mathrm{T}^* Q)\big/G$,
that is, even though $\Sigm_A[L]$ embeds into
$\left(C^{\infty}(\mathrm{T}^* Q)\right)^G$ in an obvious way, this
embedding cannot be onto a Fr\'echet dense subalgebra. Consequently,
a statement of the kind spelled out in Theorem \ref{general} no
longer suffices to recover the Poisson algebra on the quotient
$(\mathrm{T}^* Q)\big/G$.

\subsection{The Poisson structure on $\mathcal{S}_A[L]$.}
Let $L''$ denote the image of $L$ in
$\mathrm{Der}(A)$. Similarly, the pair $(A,L'')$ acquires
a
Lie-Rinehart structure, and the surjection from $L$ to $L''$ fits
into an extension
\begin{equation} \mathrm {\mathbf e}
\colon  0 \longrightarrow L' \longrightarrow L \longrightarrow L''
\longrightarrow 0 \label{ext11}
\def\emphas{}
\end{equation}
of $(\mathbb R,A)$-Lie algebras of the kind \eqref{2.1}. When $G$ is
a Lie group and when the $G$-action on $Q$ is principal, the
extension \eqref{ext11} is the sequence of the spaces of sections of the Atiyah sequence \eqref{2.2.1}.  In the general case, 
the orbit space $(\mathrm{T}^*Q)/G$ need
no longer be a smooth manifold, 
the algebra $\Sigm_A[L]$ embeds into
$\left(C^{\infty}(\mathrm{T}^*Q)\right)^G$
and
the tautological Poisson
algebra $\Sigm_A[L]$ 
yields a Poisson algebra of
continuous 
functions on $(\mathrm{T}^*Q)/G$.
Thus, provided that the
extension \eqref{ext11} splits in the category of $A$-modules,
Theorem \ref{general} applies and furnishes an explicit description
of the induced Poisson structure on the algebra $\Sigm_A[L]$ of
continuous functions on the orbit space $(\mathrm{T}^*Q)/G$.
As for the algebra
$\Sigm_A[L'']$, it is tempting to try to interpret 
it as an algebra of continuous functions on the total space
$\mathrm{T}^*B$ of the cotangent bundle of $B$
but this interpretation is only available when 
the orbit space
$B=Q/G$ is a smooth manifold;
in general we can think of 
$\Sigm_A[L'']$ as a replacement for an algebra
of functions on $\mathrm{T}^*B$.

\begin{Theorem} \label{sec} Suppose that the orbit space $B=Q\big/G$ is a
smooth manifold in such a way that the projection $\pi \colon  Q
\to B$ is a smooth submersion. Then the algebra  $A$ coincides
with the algebra $C^{\infty}(B)$ of ordinary smooth functions on
$B$, the $(\mathbb R,A)$-Lie algebra $L''$ amounts to the
$(\mathbb R,A)$-Lie algebra $\mathfrak{X}(B)$ of ordinary smooth
vector fields on $B$, the extension {\rm \eqref{ext11}} splits in
the category of $A$-modules, that is, admits a connection, and the
identities {\rm \eqref{3.16.2.11}--\eqref{3.16.2.33}} yield an
explicit description of the induced Poisson algebra $\Sigm_A[L]$ of
continuous functions on the orbit space $(\mathrm{T}^*Q)/G$.
\end{Theorem}

\begin{proof}
The algebra  $A$ plainly coincides with the algebra
$C^{\infty}(B)$ of ordinary smooth functions on $B$, and the
regularity assumption implies that the injection from $L''$ into the
$(\mathbb R,A)$-Lie algebra $\mathrm{Der}(A)\cong
\mathfrak{X}(B)$ of ordinary smooth vector fields on $B$ is surjective
as well. Since,  as an $A$-module, $\mathfrak{X}(B)$ is projective, the
extension \eqref{ext11} splits in the category of $A$-modules.
\end{proof}

A special case of Theorem \ref{sec} arises when 
$\pi\colon  Q \rightarrow B$ is a principal 
$G$-bundle. This case has been dealt with in
Theorem \ref{alg} above. In particular, the induced Poisson algebra
$\Sigm_A[L]$ then yields a complete description of the Poisson
algebra $\left(C^{\infty}(\mathrm{T}^*Q)^G,\{\,\cdot\, ,\, \cdot\,
\}\right)$ on the quotient space $(\mathrm{T}^*Q)/G$. In Theorem
\ref{sec}, we do {\em not\/} assert that, for general $G$ and $Q$, 
we obtain
a
\def\emphas{}
{\em complete\/} description of the Poisson algebra
$\left(C^{\infty}(\mathrm{T}^*Q)^G,\{\,\cdot\, ,\, \cdot\,
\}\right)$ on the quotient space $(\mathrm{T}^*Q)/G$; 
we recover 
only the subalgebra $\Sigm_A[L]$ and, in general, the relationship between $\Sigm_A[L]$ and
$C^{\infty}(\mathrm{T}^*Q)^G$ is more subtle than that between the
functions on the total space of a cotangent bundle that are
polynomial on the fibers and all smooth functions on that total
space. In particular, the algebra $\Sigm_A[L]$ is not necessarily a Fr\'echet dense subalgebra of $C^{\infty}(\mathrm
T^*Q)^G$, and invariants that cannot be recovered from
$\Sigm_A[L]$ may show up. We shall illustrate this fact in Subsection
\ref{homog} 
below.

Under the circumstances of  Theorem \ref{sec}, even though the
quotient $(\mathbb R,A)$-Lie algebra $L''$ is the $(\mathbb
R,A)$-Lie algebra of smooth vector fields on an ordinary smooth
manifold and hence arises from an ordinary Lie algebroid over $B$,
the $(\mathbb R,A)$-Lie algebra $L$ is not necessarily projective as
an $A$-module and hence does not necessarily arise from an ordinary
Lie algebroid over $B$ unless the $G$-action is principal.

For illustration, suppose that  the $G$-action on $Q$ 
is proper and has a single orbit type. 
Then the orbit space $B=Q/G$ is a smooth manifold, the
projection $\pi\colon Q \rightarrow B$ is a smooth locally
trivial fiber bundle (cf. \cite{Palais1961}),
and Theorem \ref{sec} applies. This situation has been
explored in \cite{hochrain}.
Moreover, under these circumstances, the assignment to a vector in the Lie algebra
$\mathfrak g$ of its induced fundamental vector field on $Q$
yields an exact sequence
\begin{equation}
Q \times \mathfrak g  \longrightarrow \mathrm{T}Q \longrightarrow
Q\times _B\mathrm{T}B \longrightarrow 0
\end{equation}
of smooth $G$-vector bundles over $Q$, each of the unlabeled
horizontal arrows being of constant rank.

\subsection{Homogeneous spaces}
\label{homog}
We will now justify the claim made in the previous subsection
that, when the $G$-action on $Q$ is no longer free,
the algebra $\Sigm_A[L]$ does not simply come down to a 
Fr\'echet dense subalgebra of $C^{\infty}(\mathrm T^*Q)^G$.

As before, let $G$ be a Lie group.
Let $H$ be a closed subgroup and consider 
the homogeneous space $Q=G/H$,
viewed as the base of the (right) principal
$H$-bundle $\pi\colon G \to G/H$ having $G$ as total space; 
left translation in $G$ turns this bundle into a principal $H$-bundle
in the category of left $G$-spaces.
Lift the left $G$-action to $\mathrm{T}(G/H)$ and 
$\mathrm{T}^*(G/H)$ in the standard manner. 
To adjust the exposition to the standard description of
homogeneous spaces, we work here 
with left $G$-manifolds rather than with right $G$-manifolds, 
and the reader will easily
translate the reasoning to right $G$-manifolds.
The subsequent reasoning is independent of our discussion
of Poisson structures and there is no need to rebuild the
theory with left $G$-actions rather than right $G$-actions.

We will now suppose that 
$G/H$ is reductive (cf. \cite{nomiztwo}), i.e., as 
an $H$-module, $\mathfrak g$ decomposes as a
 direct sum $\mathfrak g = \mathfrak h \oplus
\mathfrak q$  in such a way that 
$[\mathfrak h,\mathfrak q]\subset \mathfrak q$.
From 
the exact sequence
\begin{equation}
0 \longrightarrow G\times \mathfrak h \longrightarrow \mathrm{T}G
\longrightarrow G \times \mathfrak q \longrightarrow 0
\end{equation}
of $H$-vector bundles over $G$ 
(spelled out for the total spaces)
we see that the canonical morphism
\begin{equation}
G\times_H \mathfrak q \longrightarrow \mathrm{T}(G\big/H)
\end{equation}
yields an isomorphism 
of smooth vector bundles over $G\big/H$ whence 
the algebra
$(C^{\infty}(\mathrm{T}^*(G/H)))^G$ of smooth 
$G$-invariant functions on $\mathrm{T}^*(G/H)$ equals 
$(C^{\infty}(\mathfrak q^*))^H$, the algebra of smooth 
$H$-invariant functions on $\mathfrak q^*$; this algebra contains the algebra $(\Sigm[\mathfrak
q])^H$ of $H$-invariant polynomials in $\mathfrak q$ in an obvious
manner and $(\Sigm[\mathfrak q])^H$ is Fr\'echet dense in
$(C^{\infty}(\mathfrak q^*))^H$. When $H$ is compact, in view of a
result in \cite{gwschwar}, the algebra
$(C^{\infty}(\mathfrak q^*))^H$ of smooth $H$-invariant functions on
$\mathfrak q^*$ is the algebra of smooth functions in the generators
of $(\Sigm[\mathfrak q])^H$.

To reconcile this description with that given earlier,
we note first
that the space $B$ of $G$-orbits in $Q = G/H$ is  a
single point; hence the algebra $A$ amounts to the ground field
$\mathbb R$, and the extension \eqref{ext11} becomes the
identity mapping $L\to L$ where $L$ is the ordinary Lie algebra
$(\mathfrak{X}(Q))^G$ of smooth $G$-invariant vector fields on
$Q=G/H$. However, inspection of the morphism
\begin{equation}
\def\emphas{}
\begin{CD}
G \times \mathfrak q @>>> \mathrm{T} Q
\\
@V{p}VV @V{\tau_Q}VV
\def\emphas{}
\\
G @>{\pi}>> Q
\end{CD}
\end{equation}
of vector bundles shows that, as vector spaces,
$(\mathfrak{X}(Q))^G\cong \mathfrak{q}^H$, the vector 
space of $H$-fixed points in $\mathfrak q$.
Plainly,
 the
algebra $(\Sigm[\mathfrak q])^H$ of $H$-invariant polynomials in
$\mathfrak q$  contains the algebra  $\Sigm[L]\cong
\Sigm[\mathfrak q^H]$ of polynomials in $\mathfrak q^H$ but in
general the two will {\em not\/} coincide. In particular, when the
$H$-representation $\mathfrak q$ does not contain the trivial
representation, $\mathfrak q^H$ is zero whereas the algebra 
$(\Sigm[\mathfrak q])^H$ of $H$-invariant polynomials in
$\mathfrak q$ is non-trivial,
unless $\mathfrak q$ is zero,
since
the space of $H$-orbits in $\mathfrak q$ does not reduce to a
point. Thus Theorem \ref{general} cannot recover a Fr\'echet dense
subalgebra of the algebra $(C^{\infty}(\mathrm{T}^*(G/H)))^G$ of
smooth $G$-invariant functions on $\mathrm{T}^*(G/H)$.

\begin{Example} \label{22.222} {\rm 
Let $G$ be a Lie group and view $Q=G$ as a homogeneous space of
the group $G\times G$ in the standard manner, that is, $G$ is
viewed as a $(G\times G)$-space via left and right translation. In
other words, $H=\Delta G (\cong G)$ being the diagonal group in
$G\times G$, $Q$ is identified with the space $(G\times G)\big/H$
of orbits in $G\times G$ relative to the diagonal action. The
decomposition
\[
\mathfrak g \oplus \mathfrak g = \mathfrak h \oplus \mathfrak q, \
\mathfrak h = \{(X,X);\, X\in \mathfrak g\},\ \mathfrak q =
\{(Y,-Y);\, Y\in \mathfrak g\}
\]
is reductive, and the algebra $(\Sigm[\mathfrak q])^H$ of
$H$-invariant polynomials on $\mathfrak q$ coincides with the
algebra $(\Sigm[\mathfrak g])^G$ of invariants on $\mathfrak g$
(under the adjoint action). When $\mathfrak g$ is semisimple,
$\mathfrak g^G$ is zero whereas $(\Sigm[\mathfrak g])^G$ is never
zero unless $\mathfrak g$ is zero. The Poisson bracket on
$(\Sigm[\mathfrak g])^G$ is zero; indeed, the algebra
$(\Sigm[\mathfrak g])^G$ is that of Casimir elements in the
Lie-Poisson algebra $\Sigm[\mathfrak g]$. Plainly, the method
which works well in the case of a principal $G$-action does not
suffice in the present situation, that is, Theorem {\rm
\ref{general}} cannot recover the algebra under discussion.}
\quad $\blacklozenge$
\end{Example}

\begin{Example} \label{23.222} {\rm
Embed $H=\mathrm{SO}(2,\mathbb R)$ into
$G=\mathrm{SO}(3,\mathbb R)$ as rotations about the vertical
axis. Then $G/H$ is the standard {\rm 2-}sphere $S^2$. In the decomposition 
$\mathfrak g = \mathfrak h \oplus \mathfrak q$  of the Lie algebra
$\mathfrak g =\mathfrak{so}(3,\mathbb R)$,  $\dim\mathfrak h=1$,
and $\mathfrak q$ is the defining (irreducible) {\rm
2-}dimensional representation of $\mathrm{SO}(2,\mathbb R)$ via
rotations in the plane $\mathfrak q$; here we write
the elements of $\mathfrak q$ 
as $\mathbf x=(x_1,x_2)$. Since the representation
$\mathfrak q$ is irreducible, $\mathfrak q^H$ is zero whereas the
algebra $(\Sigm[\mathfrak q])^H$ of $H$-invariant polynomials on
$\mathfrak q$ is the polynomial algebra in the single variable
$\mathbf x^2= x_1^2 + x_2^2$. Thus, the space 
$\mathrm{T}^*S^2$
being viewed as the total space of a vector bundle on $S^2$ with a
$\mathrm{O}(3,\mathbb R)$-invariant Riemannian structure, the
assignment to a covector $\alpha_q\in \mathrm{T}^*S^2$ ($q \in
S^2$) of the square $\alpha_q^2$ of its length relative to the
Riemannian structure  induces a {\em diffeomorphism\/} (beware:
the notion of diffeomorphism being suitably interpreted, see
below) between $(\mathrm{T}^*S^2)\big/G$ and $\mathbb R_{\geq 0}$,
and the Poisson bracket on the resulting algebra
$(C^{\infty}(\mathrm{T}^*S^2))^G$ is zero. Here
$C^{\infty}(\mathbb R_{\geq 0})$ is {\em not\/} an algebra of
ordinary smooth functions: the space $\mathbb R_{\geq 0}$ is
endowed with the algebra of {\em Whitney\/} functions relative to
the embedding into $\mathbb R$; these are continuous functions
that are restrictions of smooth functions on some open interval of
the kind $]-\varepsilon,+\infty[$ where $\varepsilon
>0$. Again, the method
which works well in the case of a principal $G$-action does not
suffice in the present situation, that is, Theorem {\rm
\ref{general}} cannot recover this kind of algebra under
discussion.}\quad $\blacklozenge$
\end{Example}

\section*{Appendix}
We
justify the claim that 
the 
left-hand
$\Sigm$-module morphism 
\begin{equation}
\Sigm \otimes _A \mathrm D_A \longrightarrow \mathrm
D_{\Sigm},\ 
s \otimes_A da \longmapsto s da,\ s \in \Sigm,\ a \in A,
\end{equation}
in \eqref{dif2} is injective.

Since $V$ is finitely generated and projective as an $A$-module,
there is a finitely generated projective $A$-module $W$ such that
$V\oplus W$ is free. Consequently there is a free $R$-module
$U$ such that
$V \oplus W= A \otimes U$. Let $\Sigm_U=\Sigm_A[A \otimes U]$.
Since the $R$-algebra 
$\Sigm[U]$ admits a canonical augmentation map $\Sigm[U]\to R$, the canonical morphism  $A \to \Sigm_U$
of $R$-modules is injective and,
since as an $R$-module,
 $\Sigm[U]$ is free,
the canonical morphism $\Sigm[U] \to \Sigm_U$
of $R$-modules is, likewise, injective.
Now
\[
\Sigm_U=\Sigm_A[A \otimes U]\cong A \otimes  \Sigm[U]
\cong \Sigm_A[V]\otimes_A \Sigm_A[W], 
\]
and the injections $A \to \Sigm_U$ and  $\Sigm[U] \to \Sigm_U$
induce a direct sum decomposition
\begin{equation}
\mathrm D_{\Sigm_U}\cong\Sigm[U]\otimes \mathrm D_A \oplus\Sigm_U\otimes U.
\end{equation}
We rewrite this decomposition in the form
\begin{equation}
\mathrm D_{\Sigm_U}\cong
 \Sigm_A[V]\otimes_A \Sigm_A[W]\otimes_A \mathrm D_A
\oplus
 \Sigm_A[V]\otimes_A \Sigm_A[W]\otimes_A(V\oplus W).
\end{equation}
This decomposition 
implies that the canonical morphism
$\Sigm_U \otimes _A \mathrm D_A \to \mathrm
D_{\Sigm_U}$ of $\Sigm_U$-modules is injective; furthermore,
the decomposition splits the resulting extension
\begin{equation}
0
\longrightarrow
\Sigm_U \otimes _A \mathrm D_A \longrightarrow \mathrm
D_{\Sigm_U} \longrightarrow \mathrm D_{\Sigm_U|A}
\longrightarrow 0 \label{dif+}
\end{equation}
of $\Sigm_U$-modules. Notice that 
$\mathrm D_{\Sigm_U|A}\cong  \Sigm_A[V]\otimes_A \Sigm_A[W]\otimes_A(V\oplus W)$
(canonically).

Recall that $\Sigm=\Sigm_A[V]$ and consider the canonical surjection
$\Sigm_U \to \Sigm$ of $A$-algebras, viewed as a morphism of $R$-algebras.
Since \eqref{dif+} is a split extension in the category of  $\Sigm_U$-modules,
applying the functor $\Sigm\otimes_{\Sigm_U}\cdot\,$ to  \eqref{dif+}, we obtain
the extension 
\begin{equation}
0
\longrightarrow
\Sigm\otimes _A \mathrm D_A \longrightarrow 
\Sigm\otimes_{\Sigm_U}\mathrm
D_{\Sigm_U} \longrightarrow 
 \Sigm\otimes_A(V\oplus W)
\longrightarrow 0 \label{dif++}
\end{equation}
of $\Sigm$-modules, necessarily split; the salient feature here is the
exactness at $\Sigm\otimes _A \mathrm D_A$
(injectivity of the second unlabeled arrow).

Let $I_W\subseteq \Sigm_A[W]$ be the kernel of the obvious
algebra surjection $\Sigm_A[W]\to A$, and recall the fundamental exact sequence
\begin{equation}
\begin{CD}
I_W/I_W^2
@>>>
A \otimes_{\Sigm_A[W]}\mathrm D_{\Sigm_A[W]}
@>>>
\mathrm D_A
@>>>
0
\end{CD}
\end{equation}
in the category of $A$-modules. Standard homological algebra yields 
an isomorphism
\[
\mathrm{Tor}^{\Sigm_A[W]}(A,A)\longrightarrow A \otimes_{\Sigm_A[W]} I_W ,
\]
and the canonical surjection $I_W \to  A \otimes_{\Sigm_A[W]} I_W$
descends to an isomorphism  $I_W\big/I_W^2 \to  A \otimes_{\Sigm_A[W]} I_W$.
Likewise, let $J_U$ denote the kernel of
 the $R$-algebra surjection $\Sigm_U \to \Sigm$, and consider
the  fundamental exact sequence 
\begin{equation}
\begin{CD}
J_U\big/J_U^2
@>>>
\Sigm\otimes_{\Sigm_U}
\mathrm
D_{\Sigm_U} 
@>>>
\mathrm D_{\Sigm}
@>>>
0 
\end{CD}
\label{funda102}
\end{equation}
associated to the surjection $\Sigm_U \to \Sigm$.

The obvious map
\[
\Sigm_U\otimes_{\Sigm_A[W]} I_W \longrightarrow \Sigm_U\otimes_{\Sigm_A[W]} \Sigm_A[W] \cong
\Sigm_U
\]
identifies $\Sigm_U\otimes_{\Sigm_A[W]} I_W$ with the kernel $J_U$ of
$\Sigm_U \to \Sigm$, and the obvious morphism
\[
\Sigm\otimes_A\mathrm{Tor}^{\Sigm_A[W]}(A,A)
\cong \Sigm\otimes_A(I_W/I_W^2) \longrightarrow J_U/J_U^2 
\cong \Sigm\otimes_{\Sigm_U}J_U\cong \mathrm{Tor}^{\Sigm_U}(\Sigm,\Sigm)
\]
is an isomorphism of $\Sigm$-modules.
Indeed, this is the canonical morphism
\[
\Sigm\otimes_A\mathrm{Tor}^{\Sigm_A[W]}(A,A)
\longrightarrow
\mathrm{Tor}^{\Sigm\otimes_A\Sigm_A[W] }(\Sigm,\Sigm)
\]
of $\Sigm$-modules, necessarily an isomorphism since $\Sigm=\Sigm_A[V]$
is projective as an $A$-module, the $A$-module
$V$ being projective by assumption.
Consequently the  fundamental exact sequence 
\eqref{funda102}
takes the form
\begin{equation}
\begin{CD}
\Sigm\otimes_A(I_W/I_W^2)
@>>>
\Sigm\otimes_{\Sigm_U}
\mathrm
D_{\Sigm_U} 
@>>>
\mathrm D_{\Sigm}
@>>>
0 .
\end{CD}
\label{funda103}
\end{equation}

The surjection $\Sigm_U \to \Sigm$ 
induces the commutative diagram
\begin{equation}
\begin{CD}
@.
@.
0
@.
0
@.
\\
@.
@.
@VVV
@VVV
@.
\\
@.
@.
\Sigm\otimes_A(I_W/I_W^2)
@>>>
\Sigm\otimes_AW
@.
\\
@.
@.
@VVV
@VVV
@.
\\
0
@>>>
\Sigm\otimes _A \mathrm D_A @>>> 
\Sigm\otimes_{\Sigm_U}
\mathrm
D_{\Sigm_U} @>>> \Sigm\otimes_A(V\oplus W) 
@>>> 0
\\
@.
@|
@VVV
@VVV
@.
\\
@.
\Sigm \otimes _A \mathrm D_A @>>> \mathrm
D_{\Sigm} @>>> \Sigm\otimes_A V 
@>>> 0 
\\
@.
@.
@VVV
@VVV
@.
\\
@.
@.
0
@.
0
@.
\end{CD}
\end{equation}
with exact rows and columns in the category of $\Sigm$-modules.
The middle horizontal row is the (split) extension \eqref{dif++},
necessarily exact as pointed out above. The right-hand vertical column
is obviously exact.
The middle vertical column is the exact sequence \eqref{funda103},
with an arrow added emanating at $0$.
The $A$-module $I_W/I_W^2$ is canonically isomorphic to $W$ whence
$\Sigm\otimes_A(I_W/I_W^2)
\to
\Sigm\otimes_AW
$
is an isomorphism of  $\Sigm$-modules.
Consequently
the middle vertical column is exact at $\Sigm\otimes_A(I_W/I_W^2)$ and hence
exact.

Diagram chasing or the snake lemma implies that the induced morphism
$\Sigm \otimes _A \mathrm D_A \to \mathrm
D_{\Sigm}$
of $\Sigm$-modules is injective.

\section*{Acknowledgements}

We are indebted  to H. Cendra, J. Marsden,  R. Montgomery 
and J. Stasheff
for various comments that helped  us improve the exposition,
to S. Chase, Y. Kosmann-Schwarzbach,
K. Mackenzie, P. Michor, J. Stasheff and A. Weinstein 
for discussions about Lie-Rinehart algebras and, last but not least,
to the referee for his remarks.
T.S. Ratiu  was partially supported by Swiss NSF grant
200021-121512. J. Huebschmann acknowledges partial support
by the CNRS and by the
Labex CEMPI (ANR-11-LABX-0007-01).


\begin{thebibliography}{19}

\bibitem{abramars} 
R. Abraham and J. E. Marsden: \emph{Foundations of Mechanics.}
Benjamin-Cummings Publishing Company, 1978.

\bibitem{almolone} R. Almeida and P. Molino: {\pemphas Suites d'Atiyah et
feuilletages transversalement complets.} {\empha C. R. Acad. Sci.
Paris I} \textbf {300} (1985), 13--15.

\bibitem{arnobook} V. I. Arnold:
\emph{Mathematical Methods of Classical Mechanics.} Graduate
Texts in Mathematics, No. 60, Springer Verlag, Berlin  $\cdot$
Heidelberg  $\cdot$   New York $\cdot$  Tokyo, 1978, 1989 (2nd
edition).

\bibitem{atiyaone} M. F. Atiyah:
{\pemphas Complex analytic connections in fibre bundles.} {\empha
Trans. Amer. Math. Soc.} \textbf{85} (1957), 181--207.


\bibitem{CW} A. Cannas da Silva  and A. Weinstein:
\emph{Geometric models for noncommutative algebras.}
Berkeley Mathematics Lecture Notes \textbf{10},
American Mathematical Society, Providence, RI,
1999.

\bibitem{CMR} H. Cendra, J.E. Marsden, T.S. Ratiu: 
Lagrangian reduction by stages, \emph{Mem. Amer. Math. Soc.}, 
\textbf{152} (2001), no. 722.

\bibitem{CeMaPeRa} H. Cendra, J.E. Marsden, S. Pekarsky,
T.S. Ratiu:
Variational principles for Lie-Poisson and Hamilton-Poincar\'e equations,
\emph{Mosc. Math. J.} \textbf{3}, (2003), 833--867.

\bibitem{courant0}
 {T. Courant}: {\pemphas Dirac manifolds},
{\empha Trans. Amer. Math. Soc.}, \textbf{319}, {(1990)}, {631--661}.

\bibitem{courant}
 {T. Courant}: {\pemphas Tangent {L}ie algebroids},
{\empha J. Phys. A}, \textbf{27}, {(1994)}, {4527--4536}.


\bibitem{greubetal} W. Greub, St. Halperin,  and R. Vanstone:
\emph{Connections, curvature, and cohomology. {V}ol. {I}: {D}e
{R}ham cohomology of manifolds and vector bundles.}
{Pure and Applied Mathematics, Vol. 47},
{Academic Press},
{New York}, {1972}.

\bibitem{hochrain} S. Hochgerner and A. Rainer:
{\pemphas Singular Poisson reduction of cotangent bundles.}
{\empha Rev. Mat. Complut. } \textbf{19} (2006), 431--466.

\bibitem{poiscoho} J. Huebschmann:
{\pemphas Poisson cohomology and quantization.} {\empha J. reine
angew. Math.} \textbf{408} (1990), 57--113.

\bibitem{souriau} J. Huebschmann:
{\pemphas On the quantization of Poisson algebras.} In: Symplectic
Geometry and Mathematical Physics, Actes du colloque en l'honneur
de Jean-Marie Souriau, P. Donato, C. Duval, J. Elhadad, G.~M.
Tuynman, eds.; Progress in Mathematics, Vol. 99, 
Birkh\"auser-Verlag, Boston $\cdot$ Basel $\cdot$ Berlin,  204--233 (1991).

\bibitem{lradq} J. Huebschmann:
{\pemphas Lie-Rinehart algebras, descent, and quantization.}
{\empha Fields Institute Communications}  \textbf{43} (2004),
295--316, {\tt math.SG/0303016}.

\bibitem{bv} J. Huebschmann:
{\pemphas Lie-Rinehart algebras, Gerstenhaber algebras, and
Batalin-Vilkovisky algebras.} {\empha Annales de l'Institut
Fourier}  \textbf{48} (1998), 425--440, {\tt math.DG/9704005}.

\bibitem{extensta} J. Huebschmann:
{\pemphas Extensions of Lie-Rinehart algebras and the Chern-Weil
construction.} In: Festschrift to honor the 60-th birthday of
Jim Stasheff, {\empha Cont. Math.}  \textbf{227} (1999), 145--176,
{\tt math.DG/9706002}.

\bibitem{duality} J. Huebschmann:
{\pemphas Duality for Lie-Rinehart algebras and the modular
class.} {\empha J. reine angew. Math.}  \textbf{510} (1999),
103--159, {\tt math.DG/9702008}.

\bibitem{banach} J. Huebschmann:
{\pemphas  Differential Batalin-Vilkovisky algebras arising from
twilled Lie-Rinehart algebras.} {\empha Banach center
publications}  \textbf{51} (2002), 87--102.

\bibitem{quasi} J. Huebschmann:
{\pemphas Higher homotopies and Maurer-Cartan algebras:
quasi-Lie-Rinehart, Gerstenhaber-, and Batalin-Vilkovisky
algebras.}  In: \textit {The Breadth of Symplectic and Poisson Geometry},
Festschrift in Honor of Alan Weinstein, J. Marsden and T. Ratiu,
eds., Progress in Mathematics, Vol. 232, Birkh\"auser-Verlag,
Boston $\cdot$ Basel $\cdot$ Berlin,  237--302 (2004), {\tt math.DG/0311294}

\bibitem{kaehler} J. Huebschmann: {\pemphas K\"ahler spaces, nilpotent
orbits, and singular reduction.} {\empha Memoirs AMS}
\textbf{172/814} (2004), Amer. Math. Society, Providence, R. I., {\tt math.DG/0104213}.

\bibitem{bedlepro} J. Huebschmann: {\pemphas Singular Poisson-K\"ahler
geometry of certain adjoint quotients.} In: \textit{The Mathematical 
Legacy of C. Ehresmann}, J. Kubarski, and R. Wolak, eds., {\empha Banach Center Publications} \textbf{76} (2007), 325--347, {\tt math.SG/0610614}.

\bibitem{adjoint} J. Huebschmann: {\pemphas Stratified K\"ahler
structures on adjoint quotients.} {\empha Differential Geometry and
its Applications} \textbf{26} (2008), 704-731, {\tt math.DG/0404141}.

\bibitem{qr}
 J. Huebschmann: K\"ahler reduction and quantization
 {\empha J. reine angew.\ Math.} \textbf {591} (2006), 75--109,
 {\tt math.SG/0207166}.

\bibitem{hurusch} J. Huebschmann, G. Rudolph and M. Schmidt:
A gauge model for quantum mechanics on a stratified space,
\emph{Commun. Math. Phys.} \textbf{286} (2009), 459--494,
{\tt hep-th/0702017}.

\bibitem{mackbtwo} K. Mackenzie:
\emph{General Theory of Lie groupoids and Lie algebroids.}
{London Math. Soc. Lecture Note Series, vol. 213},
{Cambridge University Press},
{Cambridge, England},
 {2005}.


\bibitem{marsden81}  J.E. Marsden: \textit{Lectures
on Geometric Methods in Mathematical Physics\/}, Volume 37, SIAM, Philadelphia, 1981. 

\bibitem{MMOPR} J.E. Marsden, G. Misio\l ek, J.-P. Ortega, M. Perlmutter, T.S. Ratiu: \emph{Hamiltonian Reduction by Stages}, Lecture Notes in Mathematics,  \textbf{1913} (2007), Springer-Verlag.

\bibitem{MP} J.E. Marsden and M. Perlmutter: The orbit bundle picture of 
cotangent bundle reduction, \emph{C. R. Math. Acad. Sci. Soc. R. Can.},
\textbf{22}(2)(2000).

\bibitem{MaRa1994}  J.E. Marsden and T.S. Ratiu: \emph{Introduction to Mechanics and Symmetry}, Texts in Applied Mathematics \textbf{17}
(1994), second ed., second printing (2003), Springer-Verlag.

\bibitem{marssche} J. Marsden and J. Scheurle:
{\pemphas The reduced Euler-Lagrange equations.} In: M. J. Enos.
ed. {\em Dynamics and Control of mechanical systems. The falling
cat and related problems.} Papers from the Fields Institute
Workshop held in Waterloo, Ontario, March 1992, {\empha Fields
Institute Communications} \textbf{1} (1993), 139--164,  American
Math. Society, Providence, R. I.

\bibitem{Mont} R. Montgomery: \emph{The Bundle Picture in Mechanics}. Ph.D. Thesis, Department of Mathematics, University of California, Berkeley, 1986.

\bibitem{MMR1984} R. Montgomery, J.E. Marsden, and T.S. Ratiu: Gauged Lie-Poisson Structures, in \emph{Fluids and Plasmas:  Geometry and Dynamics\/} (J. Marsden, ed.) {\it Cont. Math.\/} {\bf 28} (1984), 101--114.


\bibitem{nomiztwo} K. Nomizu: {\pemphas Invariant affine connections on
homogeneous spaces.} {\empha Amer. J. of Math.} \textbf{76}
(1954), 33--65.

\bibitem{or} J.-P. Ortega and T.S. Ratiu: \emph{Momentum Maps and Hamiltonian Reduction\/}, Progress in Mathematics \textbf{222} (2004), 
Birkh\"auser-Verlag, Boston.

\bibitem{Palais1961}
R.~S. Palais: On the existence of slices for actions of non-compact 
Lie groups. {\empha Ann. of 
Math.} \textbf{73}(2) (1961), 295--323.

\bibitem{PeRO2009}
M. Perlmutter and M. Rodr\'{\i}guez-Olmos:
On singular Poisson Sternberg spaces.
\emph{J. Symplectic Geom.} \text{7}(2) (2009), 15--49. 

\bibitem{PeROSD2007}
M. Perlmutter, M. Rodr\'{\i}guez-Olmos, and M. E. 
Sousa-Dias: On the geometry of reduced cotangent bundles at zero momentum. \emph{J. Geom. Phys.}
\textbf{5}(2) (2007),  571--596.

\bibitem{rinehart} G. Rinehart:
{\pemphas Differential forms for general commutative algebras.}
{\empha  Trans. Amer. Math. Soc.}  \textbf{108} (1963), 195--222.

\bibitem{gwschwar}  G. W. Schwarz: {\pemphas Smooth functions invariant
under the action of a compact Lie group.} {\empha Topology}
\textbf{14} (1975), 63--68.

\bibitem{SjLe1991}
R. Sjamaar and E. Lerman: {\pemphas Stratified symplectic spaces and reduction.} {\empha Ann. of Math.} \textbf{134}(2) (1991), 375--422.

\bibitem{S77} S. Sternberg, On minimal coupling and the symplectic mechanics of a classical particle in the presence of a Yang-Mills field, \emph{Proc. Nat. Acad. Sci.}, \textbf{74} (1977), 5253--5254.

\bibitem{vaisman} I. Vaisman:
\emph{Lectures on the Geometry of Poisson Manifolds}.
Progress in Mathematics \textbf{118}, 
Birkh\"auser-Verlag,
Boston, 1994.

\bibitem{Weinstein1978}
A. Weinstein, {\pemphas A universal phase space for particles in Yang-Mills fields.}  {\empha Lett. Math. Phys.} \textbf{2}(5) (1977/78), 417--420

\bibitem{zaalaone} N. Zaalani:
{\pemphas Phase space reduction and Poisson structure.} {\empha J. Math. Phys.} \textbf{40} (1999), 3431--3438.

\end{thebibliography}
\end{document}